\documentclass[12pt,reqno]{amsart}

\def\setRevDate $#1 #2 #3${#2}
\def\TeXdrawId{\setRevDate $Date: 1995/12/19 16:40:42 $ TeXdraw V2R0}

\chardef\catamp=\the\catcode`\@
\catcode`\@=11

\long                              
\def\centertexdraw #1{\hbox to \hsize{\hss
                                      \btexdraw #1\etexdraw
                                      \hss}}

\def\btexdraw {\x@pix=0             \y@pix=0
               \x@segoffpix=\x@pix  \y@segoffpix=\y@pix
               \t@exdrawdef
               \setbox\t@xdbox=\vbox\bgroup\offinterlineskip
                   \global\d@bs=0           
                   \global\t@extonlytrue    
                   \global\p@osinitfalse
                   \s@avemove \x@pix \y@pix 
                   \m@pendingfalse
                   \global\p@osinitfalse    
                   \p@athfalse
                   \the\everytexdraw}

\def\etexdraw {\ift@extonly \else
                 \t@drclose      
               \fi
               \egroup           
               \ifdim \wd\t@xdbox>0pt
                 \t@xderror {TeXdraw box non-zero size,
                             possible extraneous text}%
               \fi
               \vbox {\offinterlineskip
                      \pixtobp \xminpix \l@lxbp  \pixtobp \yminpix \l@lybp
                      \pixtobp \xmaxpix \u@rxbp  \pixtobp \ymaxpix \u@rybp
                      \hbox{\t@xdinclude 
                        [{\l@lxbp},{\l@lybp}][{\u@rxbp},{\u@rybp}]{\p@sfile}}%
                      \pixtodim \xminpix \t@xpos  \pixtodim \yminpix \t@ypos
                      \kern \t@ypos
                      \hbox {\kern -\t@xpos
                             \box\t@xdbox       
                             \kern \t@xpos}%
                      \kern -\t@ypos\relax}}

\def\drawdim #1 {\def\d@dim{#1\relax}}

\def\setunitscale #1 {\edef\u@nitsc{#1}%
                      \realmult \u@nitsc \s@egsc \d@sc}
\def\relunitscale #1 {\realmult {#1}\u@nitsc \u@nitsc
                      \realmult \u@nitsc \s@egsc \d@sc}
\def\setsegscale #1 {\edef\s@egsc {#1}%
                     \realmult \u@nitsc \s@egsc \d@sc}
\def\relsegscale #1 {\realmult {#1}\s@egsc \s@egsc
                     \realmult \u@nitsc \s@egsc \d@sc}

\def\bsegment {\ifp@ath
                 \f@lushbs
                 \f@lushmove
               \fi
               \begingroup
               \x@segoffpix=\x@pix
               \y@segoffpix=\y@pix
               \setsegscale 1
               \global\advance \d@bs by 1\relax}
\def\esegment {\endgroup
               \ifnum \d@bs=0
                 \writetx {es}%
               \else
                 \global\advance \d@bs by -1
               \fi}

\def\savecurrpos (#1 #2){\getsympos (#1 #2)\a@rgx\a@rgy
                         \s@etcsn \a@rgx {\the\x@pix}%
                         \s@etcsn \a@rgy {\the\y@pix}}
\def\savepos (#1 #2)(#3 #4){\getpos (#1 #2)\a@rgx\a@rgy
                            \coordtopix \a@rgx \t@pixa
                            \advance \t@pixa by \x@segoffpix
                            \coordtopix \a@rgy \t@pixb
                            \advance \t@pixb by \y@segoffpix
                            \getsympos (#3 #4)\a@rgx\a@rgy
                            \s@etcsn \a@rgx {\the\t@pixa}%
                            \s@etcsn \a@rgy {\the\t@pixb}}

\def\linewd #1 {\coordtopix {#1}\t@pixa
                \f@lushbs
                \writetx {\the\t@pixa\space sl}}
\def\setgray #1 {\f@lushbs
                 \writetx {#1 sg}}
\def\lpatt (#1){\listtopix (#1)\p@ixlist
                \f@lushbs
                \writetx {[\p@ixlist] sd}}

\def\lvec (#1 #2){\getpos (#1 #2)\a@rgx\a@rgy
                  \s@etpospix \a@rgx \a@rgy
                  \writeps {\the\x@pix\space \the\y@pix\space lv}}
\def\rlvec (#1 #2){\getpos (#1 #2)\a@rgx\a@rgy
                   \r@elpospix \a@rgx \a@rgy
                   \writeps {\the\x@pix\space \the\y@pix\space lv}}
\def\move (#1 #2){\getpos (#1 #2)\a@rgx\a@rgy
                  \s@etpospix \a@rgx \a@rgy
                  \s@avemove \x@pix \y@pix}
\def\rmove (#1 #2){\getpos (#1 #2)\a@rgx\a@rgy
                   \r@elpospix \a@rgx \a@rgy
                   \s@avemove \x@pix \y@pix}

\def\lcir r:#1 {\coordtopix {#1}\t@pixa
                \writeps {\the\t@pixa\space cr}%
                \r@elupd \t@pixa \t@pixa
                \r@elupd {-\t@pixa}{-\t@pixa}}
\def\fcir f:#1 r:#2 {\coordtopix {#2}\t@pixa
                     \writeps {\the\t@pixa\space #1 fc}%
                     \r@elupd \t@pixa \t@pixa
                     \r@elupd {-\t@pixa}{-\t@pixa}}
\def\lellip rx:#1 ry:#2 {\coordtopix {#1}\t@pixa
                     \coordtopix {#2}\t@pixb
                     \writeps {\the\t@pixa\space \the\t@pixb\space el}%
                     \r@elupd \t@pixa \t@pixb
                     \r@elupd {-\t@pixa}{-\t@pixb}}
\def\fellip f:#1 rx:#2 ry:#3 {\coordtopix {#2}\t@pixa
                     \coordtopix {#3}\t@pixb
                     \writeps {\the\t@pixa\space \the\t@pixb\space #1 fe}%
                     \r@elupd \t@pixa \t@pixb
                     \r@elupd {-\t@pixa}{-\t@pixb}}
\def\larc r:#1 sd:#2 ed:#3 {\coordtopix {#1}\t@pixa
                            \writeps {\the\t@pixa\space #2 #3 ar}}

\def\ifill f:#1 {\writeps {#1 fl}}     
\def\lfill f:#1 {\writeps {#1 fp}}     

\def\htext #1{\def\testit {#1}%
              \ifx \testit\l@paren
                \let\next=\h@move
              \else
                \let\next=\h@text
              \fi
              \next {#1}}

\def\rtext td:#1 #2{\def\testit {#2}%
                    \ifx \testit\l@paren
                      \let\next=\r@move
                    \else
                      \let\next=\r@text
                    \fi
                    \next td:#1 {#2}}

\def\textref h:#1 v:#2 {\ifx #1R%
                          \edef\l@stuff {\hss}\edef\r@stuff {}%
                        \else
                          \ifx #1C%
                            \edef\l@stuff {\hss}\edef\r@stuff {\hss}%
                          \else  
                            \edef\l@stuff {}\edef\r@stuff {\hss}%
                          \fi
                        \fi
                        \ifx #2T%
                          \edef\t@stuff {}\edef\b@stuff {\vss}%
                        \else
                          \ifx #2C%
                            \edef\t@stuff {\vss}\edef\b@stuff {\vss}%
                          \else  
                            \edef\t@stuff {\vss}\edef\b@stuff {}%
                          \fi
                        \fi}

\def\avec (#1 #2){\getpos (#1 #2)\a@rgx\a@rgy
                  \s@etpospix \a@rgx \a@rgy
                  \writeps {\the\x@pix\space \the\y@pix\space (\a@type) %
                            \the\a@lenpix\space \the\a@widpix\space av}}

\def\ravec (#1 #2){\getpos (#1 #2)\a@rgx\a@rgy
                   \r@elpospix \a@rgx \a@rgy
                   \writeps {\the\x@pix\space \the\y@pix\space (\a@type) %
                             \the\a@lenpix\space \the\a@widpix\space av}}

\def\arrowheadsize l:#1 w:#2 {\coordtopix{#1}\a@lenpix
                              \coordtopix{#2}\a@widpix}
\def\arrowheadtype t:#1 {\edef\a@type{#1}}

\def\clvec (#1 #2)(#3 #4)(#5 #6)%
           {\getpos (#1 #2)\a@rgx\a@rgy
            \coordtopix \a@rgx\t@pixa
            \advance \t@pixa by \x@segoffpix
            \coordtopix \a@rgy\t@pixb
            \advance \t@pixb by \y@segoffpix
            \getpos (#3 #4)\a@rgx\a@rgy
            \coordtopix \a@rgx\t@pixc
            \advance \t@pixc by \x@segoffpix
            \coordtopix \a@rgy\t@pixd
            \advance \t@pixd by \y@segoffpix
            \getpos (#5 #6)\a@rgx\a@rgy
            \s@etpospix \a@rgx \a@rgy
            \writeps {\the\t@pixa\space \the\t@pixb\space
                      \the\t@pixc\space \the\t@pixd\space
                      \the\x@pix\space \the\y@pix\space cv}}

\def\drawbb {\bsegment
               \drawdim bp
               \linewd 0.24       
               \setunitscale {\p@sfactor}
               \writeps {\the\xminpix\space \the\yminpix\space mv}%
               \writeps {\the\xminpix\space \the\ymaxpix\space lv}%
               \writeps {\the\xmaxpix\space \the\ymaxpix\space lv}%
               \writeps {\the\xmaxpix\space \the\yminpix\space lv}%
               \writeps {\the\xminpix\space \the\yminpix\space lv}%
             \esegment}

\def\getpos (#1 #2)#3#4{\g@etargxy #1 #2 {} \\#3#4%
                        \c@heckast #3%
                        \ifa@st
                          \g@etsympix #3\t@pixa
                          \advance \t@pixa by -\x@segoffpix
                          \pixtocoord \t@pixa #3%
                        \fi
                        \c@heckast #4%
                        \ifa@st
                          \g@etsympix #4\t@pixa
                          \advance \t@pixa by -\y@segoffpix
                          \pixtocoord \t@pixa #4%
                        \fi}

\def\getsympos (#1 #2)#3#4{\g@etargxy #1 #2 {} \\#3#4%
                     \c@heckast #3%
                     \ifa@st \else
                       \t@xderror {TeXdraw: invalid symbolic coordinate}%
                     \fi
                     \c@heckast #4%
                     \ifa@st \else
                       \t@xderror {TeXdraw: invalid symbolic coordinate}%
                     \fi}

\def\listtopix (#1)#2{\def #2{}%
                      \edef\l@ist {#1 }
                      \m@oretrue
                      \loop
                        \expandafter\g@etitem \l@ist \\\a@rgx\l@ist
                        \a@pppix \a@rgx #2%
                        \ifx \l@ist\empty
                          \m@orefalse
                        \fi
                      \ifm@ore
                      \repeat}

\def\realmult #1#2#3{\dimen0=#1pt
                     \dimen2=#2\dimen0
                     \edef #3{\expandafter\c@lean\the\dimen2}}

\def\intdiv #1#2#3{\t@counta=#1
                   \t@countb=#2
                   \ifnum \t@countb<0
                      \t@counta=-\t@counta
                      \t@countb=-\t@countb
                   \fi
                   \t@countd=1                    
                   \ifnum \t@counta<0
                      \t@counta=-\t@counta
                      \t@countd=-1
                   \fi
                   \t@countc=\t@counta  \divide \t@countc by \t@countb
                   \t@counte=\t@countc  \multiply \t@counte by \t@countb
                   \advance \t@counta by -\t@counte
                   \t@counte=-1
                   \loop
                     \advance \t@counte by 1
                   \ifnum \t@counte<16                  
                       \multiply \t@countc by 2           
                       \multiply \t@counta by 2           
                       \ifnum \t@counta<\t@countb \else   
                         \advance \t@countc by 1          
                         \advance \t@counta by -\t@countb 
                       \fi
                   \repeat
                   \divide \t@countb by 2         
                   \ifnum \t@counta<\t@countb     
                     \advance \t@countc by 1
                   \fi
                   \ifnum \t@countd<0             
                     \t@countc=-\t@countc
                   \fi
                   \dimen0=\t@countc sp           
                   \edef #3{\expandafter\c@lean\the\dimen0}}

\def\coordtopix #1#2{\dimen0=#1\d@dim
                     \dimen2=\d@sc\dimen0
                     \t@counta=\dimen2          
                     \t@countb=\s@ppix
                     \divide \t@countb by 2
                     \ifnum \t@counta<0         
                       \advance \t@counta by -\t@countb
                     \else
                       \advance \t@counta by \t@countb
                     \fi
                     \divide \t@counta by \s@ppix
                     #2=\t@counta}

\def\pixtocoord #1#2{\t@counta=#1%
                     \multiply \t@counta by \s@ppix
                     \dimen0=\d@sc\d@dim
                     \t@countb=\dimen0
                     \intdiv \t@counta \t@countb #2}

\def\pixtodim #1#2{\t@countb=#1%
                   \multiply \t@countb by \s@ppix
                   #2=\t@countb sp\relax}

\def\pixtobp #1#2{\dimen0=\p@sfactor pt
                  \t@counta=\dimen0
                  \multiply \t@counta by #1%
                  \ifnum \t@counta < 0             
                    \advance \t@counta by -32768
                  \else
                    \advance \t@counta by 32768
                  \fi
                  \divide \t@counta by 65536
                  \edef #2{\the\t@counta}}
                  
\newcount\t@counta    \newcount\t@countb   
\newcount\t@countc    \newcount\t@countd
\newcount\t@counte
\newcount\t@pixa      \newcount\t@pixb     
\newcount\t@pixc      \newcount\t@pixd

\newdimen\t@xpos      \newdimen\t@ypos

\newcount\xminpix      \newcount\xmaxpix
\newcount\yminpix      \newcount\ymaxpix

\newcount\a@lenpix     \newcount\a@widpix

\newcount\x@pix        \newcount\y@pix
\newcount\x@segoffpix  \newcount\y@segoffpix
\newcount\x@savepix    \newcount\y@savepix

\newcount\s@ppix       

\newcount\d@bs

\newcount\t@xdnum
\global\t@xdnum=0

\newbox\t@xdbox

\newwrite\drawfile

\newif\ifm@pending
\newif\ifp@ath
\newif\ifa@st
\newif\ifm@ore
\newif \ift@extonly
\newif\ifp@osinit

\newtoks\everytexdraw

\def\l@paren{(}
\def\a@st{*}

\catcode`\%=12
  \def\p@b {
\catcode`\%=14
\catcode`\{=12  \catcode`\}=12  \catcode`\u=1 \catcode`\v=2
  \def\l@br u{v  \def\r@br u}v
\catcode `\{=1  \catcode`\}=2   \catcode`\u=11 \catcode`\v=11

{\catcode`\p=12 \catcode`\t=12
 \gdef\c@lean #1pt{#1}}

\def\sppix#1/#2 {\dimen0=1#2 \s@ppix=\dimen0
                 \t@counta=#1%
                 \divide \t@counta by 2
                 \advance \s@ppix by \t@counta
                 \divide \s@ppix by #1
                 \t@counta=\s@ppix
                 \multiply \t@counta by 65536       
                 \advance \t@counta by 32891        
                 \divide \t@counta by 65782         
                 \dimen0=\t@counta sp
                 \edef\p@sfactor {\expandafter\c@lean\the\dimen0}}

\def\g@etargxy #1 #2 #3 #4\\#5#6{\def #5{#1}%
                           \ifx #5\empty
                             \g@etargxy #2 #3 #4 \\#5#6
                           \else
                             \def #6{#2}%
                             \def\next {#3}%
                             \ifx \next\empty \else
                               \t@xderror {TeXdraw: invalid coordinate}%
                             \fi
                           \fi}

\def\c@heckast #1{\expandafter
                  \c@heckastll #1\\}
\def\c@heckastll #1#2\\{\def\testit {#1}%
                        \ifx \testit\a@st
                          \a@sttrue
                        \else
                          \a@stfalse
                        \fi}

\def\g@etsympix #1#2{\expandafter
               \ifx \csname #1\endcsname \relax
                 \t@xderror {TeXdraw: undefined symbolic coordinate}%
               \fi
               #2=\csname #1\endcsname}

\def\s@etcsn #1#2{\expandafter
                  \xdef\csname#1\endcsname {#2}}

\def\g@etitem #1 #2\\#3#4{\edef #4{#2}\edef #3{#1}}
\def\a@pppix #1#2{\edef\next {#1}%
                  \ifx \next\empty \else
                    \coordtopix {#1}\t@pixa
                    \ifx #2\empty
                      \edef #2{\the\t@pixa}%
                    \else
                      \edef #2{#2 \the\t@pixa}%
                    \fi
                  \fi}

\def\s@etpospix #1#2{\coordtopix {#1}\x@pix
                     \advance \x@pix by \x@segoffpix
                     \coordtopix {#2}\y@pix
                     \advance \y@pix by \y@segoffpix
                     \u@pdateminmax \x@pix \y@pix}

\def\r@elpospix #1#2{\coordtopix {#1}\t@pixa
                     \advance \x@pix by \t@pixa
                     \coordtopix {#2}\t@pixa
                     \advance \y@pix by \t@pixa
                     \u@pdateminmax \x@pix \y@pix}

\def\r@elupd #1#2{\t@counta=\x@pix
                  \advance\t@counta by #1%
                  \t@countb=\y@pix
                  \advance\t@countb by #2%
                  \u@pdateminmax \t@counta \t@countb}

\def\u@pdateminmax #1#2{\ifnum #1>\xmaxpix
                          \global\xmaxpix=#1%
                        \fi
                        \ifnum #1<\xminpix
                          \global\xminpix=#1%
                        \fi
                        \ifnum #2>\ymaxpix
                          \global\ymaxpix=#2%
                        \fi
                        \ifnum #2<\yminpix
                          \global\yminpix=#2%
                        \fi}

\def\s@avemove #1#2{\x@savepix=#1\y@savepix=#2%
                    \m@pendingtrue
                    \ifp@osinit \else
                      \global\p@osinittrue
                      \global\xminpix=\x@savepix \global\yminpix=\y@savepix
                      \global\xmaxpix=\x@savepix \global\ymaxpix=\y@savepix
                    \fi}

\def\f@lushmove {\global\p@osinittrue
                 \ifm@pending
                   \writetx {\the\x@savepix\space \the\y@savepix\space mv}%
                   \m@pendingfalse
                   \p@athfalse
                 \fi}

\def\f@lushbs {\loop
                 \ifnum \d@bs>0
                   \writetx {bs}%
                   \global\advance \d@bs by -1
               \repeat}
               
\def\h@move #1#2 #3)#4{\move (#2 #3)%
                       \h@text {#4}}
\def\h@text #1{\pixtodim \x@pix \t@xpos
               \pixtodim \y@pix \t@ypos
               \vbox to 0pt{\normalbaselines
                            \t@stuff
                            \kern -\t@ypos
                            \hbox to 0pt{\l@stuff
                                         \kern \t@xpos
                                         \hbox {#1}%
                                         \kern -\t@xpos
                                         \r@stuff}%
                            \kern \t@ypos
                            \b@stuff\relax}}

\def\r@move td:#1 #2#3 #4)#5{\move (#3 #4)%
                             \r@text td:#1 {#5}}
\def\r@text td:#1 #2{\vbox to 0pt{\pixtodim \x@pix \t@xpos
                                  \pixtodim \y@pix \t@ypos
                                  \kern -\t@ypos
                                  \hbox to 0pt{\kern \t@xpos
                                               \rottxt {#1}{\z@sb {#2}}%
                                               \hss}%
                                  \vss}}
\def\z@sb #1{\vbox to 0pt{\normalbaselines
                          \t@stuff
                          \hbox to 0pt{\l@stuff \hbox {#1}\r@stuff}%
                          \b@stuff}}

\ifx \rotatebox\@undefined
  \def\rottxt #1#2{\bgroup
                     #2%
                   \egroup}
\else
  \let\rottxt=\rotatebox
\fi

\ifx \t@xderror\@undefined
  \let\t@xderror=\errmessage
\fi

\def\t@exdrawdef {\sppix 300/in            
                  \drawdim in              
                  \edef\u@nitsc {1}
                  \setsegscale 1           
                  \arrowheadsize l:0.16 w:0.08
                  \arrowheadtype t:T
                  \textref h:L v:B }

\ifx \includegraphics\@undefined
  \def\t@xdinclude [#1,#2][#3,#4]#5{%
    \begingroup                           
      \message {<#5>}%
      \leavevmode
      \t@counta=-#1
      \t@countb=-#2%
      \setbox0=\hbox{%
        \includegraphics{#5}}%
      \t@ypos=#4 bp%
        \advance \t@ypos by -#2 bp%
      \t@xpos=#3 bp%
        \advance \t@xpos by -#1 bp%
      \dp0=0pt \ht0=\t@ypos  \wd0=\t@xpos
      \box0%
    \endgroup}
\else
  \let\t@xdinclude=\includegraphics
\fi

\def\writeps #1{\f@lushbs
                \f@lushmove
                \p@athtrue
                \writetx {#1}}
\def\writetx #1{\ift@extonly
                  \global\t@extonlyfalse
                  \t@xdpsfn \p@sfile
                  \t@dropen \p@sfile
                \fi
                \w@rps {#1}}
\def\w@rps #1{\immediate\write\drawfile {#1}}

\def\t@xdpsfn #1{%
  \global\advance \t@xdnum by 1
  \ifnum \t@xdnum<10
    \xdef #1{\jobname.ps\the\t@xdnum}
  \else
    \xdef #1{\jobname.p\the\t@xdnum}
  \fi
}
\def\t@dropen #1{%
  \immediate\openout\drawfile=#1%
  \w@rps {\p@b PS-Adobe-3.0 EPSF-3.0}%
  \w@rps {\p@p BoundingBox: (atend)}%
  \w@rps {\p@p Title: TeXdraw drawing: #1}%
  \w@rps {\p@p Pages: 1}%
  \w@rps {\p@p Creator: \TeXdrawId}%
  \w@rps {\p@p CreationDate: \the\year/\the\month/\the\day}%
  \w@rps {50 dict begin}%
  \w@rps {/mv {stroke moveto} def}%
  \w@rps {/lv {lineto} def}%
  \w@rps {/st {currentpoint stroke moveto} def}%
  \w@rps {/sl {st setlinewidth} def}%
  \w@rps {/sd {st 0 setdash} def}%
  \w@rps {/sg {st setgray} def}%
  \w@rps {/bs {gsave} def /es {stroke grestore} def}%
  \w@rps {/fl \l@br gsave setgray fill grestore}%
  \w@rps    { currentpoint newpath moveto\r@br\space def}%
  \w@rps {/fp {gsave setgray fill grestore st} def}%
  \w@rps {/cv {curveto} def}%
  \w@rps {/cr \l@br gsave currentpoint newpath 3 -1 roll 0 360 arc}%
  \w@rps    { stroke grestore\r@br\space def}%
  \w@rps {/fc \l@br gsave setgray currentpoint newpath}%
  \w@rps    { 3 -1 roll 0 360 arc fill grestore\r@br\space def}%
  \w@rps {/ar {gsave currentpoint newpath 5 2 roll arc stroke grestore} def}%
  \w@rps {/el \l@br gsave /svm matrix currentmatrix def}%
  \w@rps    { currentpoint translate scale newpath 0 0 1 0 360 arc}%
  \w@rps    { svm setmatrix stroke grestore\r@br\space def}%
  \w@rps {/fe \l@br gsave setgray currentpoint translate scale newpath}%
  \w@rps    { 0 0 1 0 360 arc fill grestore\r@br\space def}%
  \w@rps {/av \l@br /hhwid exch 2 div def /hlen exch def}%
  \w@rps    { /ah exch def /tipy exch def /tipx exch def}%
  \w@rps    { currentpoint /taily exch def /tailx exch def}%
  \w@rps    { /dx tipx tailx sub def /dy tipy taily sub def}%
  \w@rps    { /alen dx dx mul dy dy mul add sqrt def}%
  \w@rps    { /blen alen hlen sub def}%
  \w@rps    { gsave tailx taily translate dy dx atan rotate}%
  \w@rps    { (V) ah ne {blen 0 gt {blen 0 lineto} if} {alen 0 lineto} ifelse}%
  \w@rps    { stroke blen hhwid neg moveto alen 0 lineto blen hhwid lineto}%
  \w@rps    { (T) ah eq {closepath} if}%
  \w@rps    { (W) ah eq {gsave 1 setgray fill grestore closepath} if}%
  \w@rps    { (F) ah eq {fill} {stroke} ifelse}%
  \w@rps    { grestore tipx tipy moveto\r@br\space def}%
  \w@rps {\p@sfactor\space \p@sfactor\space scale}%
  \w@rps {1 setlinecap 1 setlinejoin}%
  \w@rps {3 setlinewidth [] 0 setdash}%
  \w@rps {0 0 moveto}%
}

\def\t@drclose {%
  \bgroup
    \w@rps {stroke end showpage}%
    \w@rps {\p@p Trailer:}%
    \pixtobp \xminpix \l@lxbp  \pixtobp \yminpix \l@lybp
    \pixtobp \xmaxpix \u@rxbp  \pixtobp \ymaxpix \u@rybp
    \w@rps {\p@p BoundingBox: \l@lxbp\space \l@lybp\space
                              \u@rxbp\space \u@rybp}%
    \w@rps {\p@p EOF}%
  \egroup
  \immediate\closeout\drawfile
}

\catcode`\@=\catamp

\catcode`\@=11
\font\tenln    = line10
\font\tenlnw   = linew10

\newskip\Einheit \Einheit=0.5cm
\newcount\xcoord \newcount\ycoord
\newdimen\xdim \newdimen\ydim \newdimen\PfadD@cke \newdimen\Pfadd@cke

\newcount\@tempcnta
\newcount\@tempcntb

\newdimen\@tempdima
\newdimen\@tempdimb

\newdimen\@wholewidth
\newdimen\@halfwidth

\newcount\@xarg
\newcount\@yarg
\newcount\@yyarg
\newbox\@linechar
\newbox\@tempboxa
\newdimen\@linelen
\newdimen\@clnwd
\newdimen\@clnht

\newif\if@negarg

\def\@whilenoop#1{}
\def\@whiledim#1\do #2{\ifdim #1\relax#2\@iwhiledim{#1\relax#2}\fi}
\def\@iwhiledim#1{\ifdim #1\let\@nextwhile=\@iwhiledim
        \else\let\@nextwhile=\@whilenoop\fi\@nextwhile{#1}}

\def\@whileswnoop#1\fi{}
\def\@whilesw#1\fi#2{#1#2\@iwhilesw{#1#2}\fi\fi}
\def\@iwhilesw#1\fi{#1\let\@nextwhile=\@iwhilesw
         \else\let\@nextwhile=\@whileswnoop\fi\@nextwhile{#1}\fi}

\def\thinlines{\let\@linefnt\tenln \let\@circlefnt\tencirc
  \@wholewidth\fontdimen8\tenln \@halfwidth .5\@wholewidth}
\def\thicklines{\let\@linefnt\tenlnw \let\@circlefnt\tencircw
  \@wholewidth\fontdimen8\tenlnw \@halfwidth .5\@wholewidth}
\thinlines

\PfadD@cke1pt \Pfadd@cke0.5pt
\def\PfadDicke#1{\PfadD@cke#1 \divide\PfadD@cke by2 \Pfadd@cke\PfadD@cke  
\multiply\PfadD@cke by2}
\long\def\LOOP#1\REPEAT{\def\BODY{#1}\ITERATE}
\def\ITERATE{\BODY \let\next\ITERATE \else\let\next\relax\fi \next}
\let\REPEAT=\fi
\def\Punkt{\hbox{\raise-2pt\hbox to0pt{\hss$\ssize\bullet$\hss}}}
\def\DuennPunkt(#1,#2){\unskip
  \raise#2 \Einheit\hbox to0pt{\hskip#1 \Einheit
          \raise-2.5pt\hbox to0pt{\hss$\bullet$\hss}\hss}}
\def\NormalPunkt(#1,#2){\unskip
  \raise#2 \Einheit\hbox to0pt{\hskip#1 \Einheit
          \raise-3pt\hbox to0pt{\hss\twelvepoint$\bullet$\hss}\hss}}
\def\DickPunkt(#1,#2){\unskip
  \raise#2 \Einheit\hbox to0pt{\hskip#1 \Einheit
          \raise-4pt\hbox to0pt{\hss\fourteenpoint$\bullet$\hss}\hss}}
\def\Kreis(#1,#2){\unskip
  \raise#2 \Einheit\hbox to0pt{\hskip#1 \Einheit
          \raise-4pt\hbox to0pt{\hss\fourteenpoint$\circ$\hss}\hss}}

\def\Line@(#1,#2)#3{\@xarg #1\relax \@yarg #2\relax
\@linelen=#3\Einheit
\ifnum\@xarg =0 \@vline
  \else \ifnum\@yarg =0 \@hline \else \@sline\fi
\fi}

\def\@sline{\ifnum\@xarg< 0 \@negargtrue \@xarg -\@xarg \@yyarg -\@yarg
  \else \@negargfalse \@yyarg \@yarg \fi
\ifnum \@yyarg >0 \@tempcnta\@yyarg \else \@tempcnta -\@yyarg \fi
\ifnum\@tempcnta>6 \@badlinearg\@tempcnta0 \fi
\ifnum\@xarg>6 \@badlinearg\@xarg 1 \fi
\setbox\@linechar\hbox{\@linefnt\@getlinechar(\@xarg,\@yyarg)}%
\ifnum \@yarg >0 \let\@upordown\raise \@clnht\z@
   \else\let\@upordown\lower \@clnht \ht\@linechar\fi
\@clnwd=\wd\@linechar
\if@negarg \hskip -\wd\@linechar \def\@tempa{\hskip -2\wd\@linechar}\else
     \let\@tempa\relax \fi
\@whiledim \@clnwd <\@linelen \do
  {\@upordown\@clnht\copy\@linechar
   \@tempa
   \advance\@clnht \ht\@linechar
   \advance\@clnwd \wd\@linechar}%
\advance\@clnht -\ht\@linechar
\advance\@clnwd -\wd\@linechar
\@tempdima\@linelen\advance\@tempdima -\@clnwd
\@tempdimb\@tempdima\advance\@tempdimb -\wd\@linechar
\if@negarg \hskip -\@tempdimb \else \hskip \@tempdimb \fi
\multiply\@tempdima \@m
\@tempcnta \@tempdima \@tempdima \wd\@linechar \divide\@tempcnta \@tempdima
\@tempdima \ht\@linechar \multiply\@tempdima \@tempcnta
\divide\@tempdima \@m
\advance\@clnht \@tempdima
\ifdim \@linelen <\wd\@linechar
   \hskip \wd\@linechar
  \else\@upordown\@clnht\copy\@linechar\fi}

\def\@hline{\ifnum \@xarg <0 \hskip -\@linelen \fi
\vrule height\Pfadd@cke width \@linelen depth\Pfadd@cke
\ifnum \@xarg <0 \hskip -\@linelen \fi}

\def\@getlinechar(#1,#2){\@tempcnta#1\relax\multiply\@tempcnta 8
\advance\@tempcnta -9 \ifnum #2>0 \advance\@tempcnta #2\relax\else
\advance\@tempcnta -#2\relax\advance\@tempcnta 64 \fi
\char\@tempcnta}

\def\Vektor(#1,#2)#3(#4,#5){\unskip\leavevmode
  \xcoord#4\relax \ycoord#5\relax
      \raise\ycoord \Einheit\hbox to0pt{\hskip\xcoord \Einheit
         \Vector@(#1,#2){#3}\hss}}

\def\Vector@(#1,#2)#3{\@xarg #1\relax \@yarg #2\relax
\@tempcnta \ifnum\@xarg<0 -\@xarg\else\@xarg\fi
\ifnum\@tempcnta<5\relax
\@linelen=#3\Einheit
\ifnum\@xarg =0 \@vvector
  \else \ifnum\@yarg =0 \@hvector \else \@svector\fi
\fi
\else\@badlinearg\fi}

\def\@hvector{\@hline\hbox to 0pt{\@linefnt
\ifnum \@xarg <0 \@getlarrow(1,0)\hss\else
    \hss\@getrarrow(1,0)\fi}}

\def\@vvector{\ifnum \@yarg <0 \@downvector \else \@upvector \fi}

\def\@svector{\@sline
\@tempcnta\@yarg \ifnum\@tempcnta <0 \@tempcnta=-\@tempcnta\fi
\ifnum\@tempcnta <5
  \hskip -\wd\@linechar
  \@upordown\@clnht \hbox{\@linefnt  \if@negarg
  \@getlarrow(\@xarg,\@yyarg) \else \@getrarrow(\@xarg,\@yyarg) \fi}%
\else\@badlinearg\fi}

\def\@upline{\hbox to \z@{\hskip -.5\Pfadd@cke \vrule width \Pfadd@cke
   height \@linelen depth \z@\hss}}

\def\@downline{\hbox to \z@{\hskip -.5\Pfadd@cke \vrule width \Pfadd@cke
   height \z@ depth \@linelen \hss}}

\def\@upvector{\@upline\setbox\@tempboxa\hbox{\@linefnt\char'66}\raise
     \@linelen \hbox to\z@{\lower \ht\@tempboxa\box\@tempboxa\hss}}

\def\@downvector{\@downline\lower \@linelen
      \hbox to \z@{\@linefnt\char'77\hss}}

\def\@getlarrow(#1,#2){\ifnum #2 =\z@ \@tempcnta='33\else
\@tempcnta=#1\relax\multiply\@tempcnta \sixt@@n \advance\@tempcnta
-9 \@tempcntb=#2\relax\multiply\@tempcntb \tw@
\ifnum \@tempcntb >0 \advance\@tempcnta \@tempcntb\relax
\else\advance\@tempcnta -\@tempcntb\advance\@tempcnta 64
\fi\fi\char\@tempcnta}

\def\@getrarrow(#1,#2){\@tempcntb=#2\relax
\ifnum\@tempcntb < 0 \@tempcntb=-\@tempcntb\relax\fi
\ifcase \@tempcntb\relax \@tempcnta='55 \or
\ifnum #1<3 \@tempcnta=#1\relax\multiply\@tempcnta
24 \advance\@tempcnta -6 \else \ifnum #1=3 \@tempcnta=49
\else\@tempcnta=58 \fi\fi\or
\ifnum #1<3 \@tempcnta=#1\relax\multiply\@tempcnta
24 \advance\@tempcnta -3 \else \@tempcnta=51\fi\or
\@tempcnta=#1\relax\multiply\@tempcnta
\sixt@@n \advance\@tempcnta -\tw@ \else
\@tempcnta=#1\relax\multiply\@tempcnta
\sixt@@n \advance\@tempcnta 7 \fi\ifnum #2<0 \advance\@tempcnta 64 \fi
\char\@tempcnta}

\def\Diagonale(#1,#2)#3{\unskip\leavevmode
  \xcoord#1\relax \ycoord#2\relax
      \raise\ycoord \Einheit\hbox to0pt{\hskip\xcoord \Einheit
         \Line@(1,1){#3}\hss}}
\def\AntiDiagonale(#1,#2)#3{\unskip\leavevmode
  \xcoord#1\relax \ycoord#2\relax 
      \raise\ycoord \Einheit\hbox to0pt{\hskip\xcoord \Einheit
         \Line@(1,-1){#3}\hss}}
\def\Pfad(#1,#2),#3\endPfad{\unskip\leavevmode
  \xcoord#1 \ycoord#2 \thicklines\ZeichnePfad#3\endPfad\thinlines}
\def\ZeichnePfad#1{\ifx#1\endPfad\let\next\relax
  \else\let\next\ZeichnePfad
    \ifnum#1=1
      \raise\ycoord \Einheit\hbox to0pt{\hskip\xcoord \Einheit
         \vrule height\Pfadd@cke width1 \Einheit depth\Pfadd@cke\hss}%
      \advance\xcoord by 1
    \else\ifnum#1=2
      \raise\ycoord \Einheit\hbox to0pt{\hskip\xcoord \Einheit
        \hbox{\hskip-\PfadD@cke\vrule height1 \Einheit width\PfadD@cke  
depth0pt}\hss}%
      \advance\ycoord by 1
    \else\ifnum#1=3
      \raise\ycoord \Einheit\hbox to0pt{\hskip\xcoord \Einheit
         \Line@(1,1){1}\hss}
      \advance\xcoord by 1
      \advance\ycoord by 1
    \else\ifnum#1=4
      \raise\ycoord \Einheit\hbox to0pt{\hskip\xcoord \Einheit
         \Line@(1,-1){1}\hss}
      \advance\xcoord by 1
      \advance\ycoord by -1
    \fi\fi\fi\fi
  \fi\next}
\def\hSSchritt{\leavevmode\raise-.4pt\hbox to0pt{\hss.\hss}\hskip.2\Einheit
  \raise-.4pt\hbox to0pt{\hss.\hss}\hskip.2\Einheit
  \raise-.4pt\hbox to0pt{\hss.\hss}\hskip.2\Einheit
  \raise-.4pt\hbox to0pt{\hss.\hss}\hskip.2\Einheit
  \raise-.4pt\hbox to0pt{\hss.\hss}\hskip.2\Einheit}
\def\vSSchritt{\vbox{\baselineskip.2\Einheit\lineskiplimit0pt
\hbox{.}\hbox{.}\hbox{.}\hbox{.}\hbox{.}}}
\def\DSSchritt{\leavevmode\raise-.4pt\hbox to0pt{%
  \hbox to0pt{\hss.\hss}\hskip.2\Einheit
  \raise.2\Einheit\hbox to0pt{\hss.\hss}\hskip.2\Einheit
  \raise.4\Einheit\hbox to0pt{\hss.\hss}\hskip.2\Einheit
  \raise.6\Einheit\hbox to0pt{\hss.\hss}\hskip.2\Einheit
  \raise.8\Einheit\hbox to0pt{\hss.\hss}\hss}}
\def\dSSchritt{\leavevmode\raise-.4pt\hbox to0pt{%
  \hbox to0pt{\hss.\hss}\hskip.2\Einheit
  \raise-.2\Einheit\hbox to0pt{\hss.\hss}\hskip.2\Einheit
  \raise-.4\Einheit\hbox to0pt{\hss.\hss}\hskip.2\Einheit
  \raise-.6\Einheit\hbox to0pt{\hss.\hss}\hskip.2\Einheit
  \raise-.8\Einheit\hbox to0pt{\hss.\hss}\hss}}
\def\SPfad(#1,#2),#3\endSPfad{\unskip\leavevmode
  \xcoord#1 \ycoord#2 \ZeichneSPfad#3\endSPfad}
\def\ZeichneSPfad#1{\ifx#1\endSPfad\let\next\relax
  \else\let\next\ZeichneSPfad
    \ifnum#1=1
      \raise\ycoord \Einheit\hbox to0pt{\hskip\xcoord \Einheit
         \hSSchritt\hss}%
      \advance\xcoord by 1
    \else\ifnum#1=2
      \raise\ycoord \Einheit\hbox to0pt{\hskip\xcoord \Einheit
        \hbox{\hskip-2pt \vSSchritt}\hss}%
      \advance\ycoord by 1
    \else\ifnum#1=3
      \raise\ycoord \Einheit\hbox to0pt{\hskip\xcoord \Einheit
         \DSSchritt\hss}
      \advance\xcoord by 1
      \advance\ycoord by 1
    \else\ifnum#1=4
      \raise\ycoord \Einheit\hbox to0pt{\hskip\xcoord \Einheit
         \dSSchritt\hss}
      \advance\xcoord by 1
      \advance\ycoord by -1
    \fi\fi\fi\fi
  \fi\next}
\def\Koordinatenachsen(#1,#2){\unskip
 \hbox to0pt{\hskip-.5pt\vrule height#2 \Einheit width.5pt depth1 \Einheit}%
 \hbox to0pt{\hskip-1 \Einheit \xcoord#1 \advance\xcoord by1
    \vrule height0.25pt width\xcoord \Einheit depth0.25pt\hss}}
\def\Koordinatenachsen(#1,#2)(#3,#4){\unskip
 \hbox to0pt{\hskip-.5pt \ycoord-#4 \advance\ycoord by1
    \vrule height#2 \Einheit width.5pt depth\ycoord \Einheit}%
 \hbox to0pt{\hskip-1 \Einheit \hskip#3\Einheit
    \xcoord#1 \advance\xcoord by1 \advance\xcoord by-#3
    \vrule height0.25pt width\xcoord \Einheit depth0.25pt\hss}}
\def\Gitter(#1,#2){\unskip \xcoord0 \ycoord0 \leavevmode
  \LOOP\ifnum\ycoord<#2
    \loop\ifnum\xcoord<#1
      \raise\ycoord \Einheit\hbox to0pt{\hskip\xcoord \Einheit\Punkt\hss}%
      \advance\xcoord by1
    \repeat
    \xcoord0
    \advance\ycoord by1
  \REPEAT}
\def\Gitter(#1,#2)(#3,#4){\unskip \xcoord#3 \ycoord#4 \leavevmode
  \LOOP\ifnum\ycoord<#2
    \loop\ifnum\xcoord<#1
      \raise\ycoord \Einheit\hbox to0pt{\hskip\xcoord \Einheit\Punkt\hss}%
      \advance\xcoord by1
    \repeat
    \xcoord#3
    \advance\ycoord by1
  \REPEAT}
\def\Label#1#2(#3,#4){\unskip \xdim#3 \Einheit \ydim#4 \Einheit
  \def\lo{\advance\xdim by-.5 \Einheit \advance\ydim by.5 \Einheit}%
  \def\llo{\advance\xdim by-.25cm \advance\ydim by.5 \Einheit}%
  \def\loo{\advance\xdim by-.5 \Einheit \advance\ydim by.25cm}%
  \def\o{\advance\ydim by.25cm}%
  \def\ro{\advance\xdim by.5 \Einheit \advance\ydim by.5 \Einheit}%
  \def\rro{\advance\xdim by.25cm \advance\ydim by.5 \Einheit}%
  \def\roo{\advance\xdim by.5 \Einheit \advance\ydim by.25cm}%
  \def\l{\advance\xdim by-.30cm}%
  \def\r{\advance\xdim by.30cm}%
  \def\lu{\advance\xdim by-.5 \Einheit \advance\ydim by-.6 \Einheit}%
  \def\llu{\advance\xdim by-.25cm \advance\ydim by-.6 \Einheit}%
  \def\luu{\advance\xdim by-.5 \Einheit \advance\ydim by-.30cm}%
  \def\u{\advance\ydim by-.30cm}%
  \def\ru{\advance\xdim by.5 \Einheit \advance\ydim by-.6 \Einheit}%
  \def\rru{\advance\xdim by.25cm \advance\ydim by-.6 \Einheit}%
  \def\ruu{\advance\xdim by.5 \Einheit \advance\ydim by-.30cm}%
  #1\raise\ydim\hbox to0pt{\hskip\xdim
     \vbox to0pt{\vss\hbox to0pt{\hss$#2$\hss}\vss}\hss}%
}
\catcode`\@=13

\def\ldreieck{\bsegment
  \rlvec(0.866025403784439 .5) \rlvec(0 -1)
  \rlvec(-0.866025403784439 .5)
  \savepos(0.866025403784439 -.5)(*ex *ey)
        \esegment
  \move(*ex *ey)
        }
\def\rdreieck{\bsegment
  \rlvec(0.866025403784439 -.5) \rlvec(-0.866025403784439 -.5)  \rlvec(0 1)
  \savepos(0 -1)(*ex *ey)
        \esegment
  \move(*ex *ey)
        }
\def\rhombus{\bsegment
  \rlvec(0.866025403784439 .5) \rlvec(0.866025403784439 -.5)
  \rlvec(-0.866025403784439 -.5)  \rlvec(0 1)
  \rmove(0 -1)  \rlvec(-0.866025403784439 .5)
  \savepos(0.866025403784439 -.5)(*ex *ey)
        \esegment
  \move(*ex *ey)
        }
\def\RhombusA{\bsegment
  \rlvec(0.866025403784439 .5) \rlvec(0.866025403784439 -.5)
  \rlvec(-0.866025403784439 -.5) \rlvec(-0.866025403784439 .5)
  \savepos(0.866025403784439 -.5)(*ex *ey)
        \esegment
  \move(*ex *ey)
        }
\def\RhombusB{\bsegment
  \rlvec(0.866025403784439 .5) \rlvec(0 -1)
  \rlvec(-0.866025403784439 -.5) \rlvec(0 1)
  \savepos(0 -1)(*ex *ey)
        \esegment
  \move(*ex *ey)
        }
\def\RhombusC{\bsegment
  \rlvec(0.866025403784439 -.5) \rlvec(0 -1)
  \rlvec(-0.866025403784439 .5) \rlvec(0 1)
  \savepos(0.866025403784439 -.5)(*ex *ey)
        \esegment
  \move(*ex *ey)
        }
\def\hdSchritt{\bsegment
  \lpatt(.05 .13)
  \rlvec(0.866025403784439 -.5)
  \savepos(0.866025403784439 -.5)(*ex *ey)
        \esegment
  \move(*ex *ey)
        }
\def\vdSchritt{\bsegment
  \lpatt(.05 .13)
  \rlvec(0 -1)
  \savepos(0 -1)(*ex *ey)
        \esegment
  \move(*ex *ey)
        }

\newcommand{\arccot}{\operatorname{arccot}}			
\newcommand{\HypFsimple}[2]{\sideset{_#1}{_#2}{F}}
\newcommand{\HypF}[5]{\sideset{_#1}{_#2}{F}\!\left[\matrix  
#3\\#4\endmatrix;{\displaystyle #5}\right]}
\def\po#1#2{(#1)_{#2}}
\newcommand{\card}[1]{\left|#1\right|}
\newcommand{\abs}[1]{\left|#1\right|}
\newcommand{\floor}[1]{\left\lfloor#1\right\rfloor}
\newcommand{\ceil}[1]{\left\lceil#1\right\rceil}
\newcommand{\arcsinh}{\operatorname{arcsinh}}
\newcommand{\sgn}{\operatorname{sgn}}
\newcommand{\mydet}[1]{\det\left(#1\right)}
\newcommand{\dx}{{\,dx}}
\def\ep{\varepsilon}
\def\la{\lambda}

\newtheorem{thm}{Theorem}

\newtheorem{cor}[thm]{Corollary}
\newtheorem{prop}[thm]{Proposition}
\newtheorem{lem}[thm]{Lemma}
\newtheorem*{conj}{Conjecture}

\theoremstyle{remark}
\newtheorem*{rem}{Remark}

\catcode`\@=11
\def\sideset#1#2#3{%
  \@mathmeasure\z@\displaystyle{#3}%
  \global\setbox\@ne\vbox to\ht\z@{}\dp\@ne\dp\z@
  \setbox\tw@\box\@ne
  \@mathmeasure4\displaystyle{\copy\tw@#1}%
  \@mathmeasure6\displaystyle{#3{#2}}%
  \dimen@-\wd6 \advance\dimen@\wd4 \advance\dimen@\wd\z@
  \hbox to\dimen@{}\mathop{\kern-\dimen@\box4\box6}%
}
\catcode`\@=12

\numberwithin{equation}{section}

\def\newstep#1{\par\vskip6pt\noindent{\it #1.\/ }}

\def\complexi{{\sqrt{-1}}}
\def\complexikl{{(\sqrt{-1})}}

\begin{document}

\title[The number of rhombus tilings]{The number of rhombus
tilings of a symmetric hexagon which contain a fixed rhombus on the
symmetry axis, II}

\author[M.~Fulmek and C.~Krattenthaler]
{M.~Fulmek and C.~Krattenthaler$^\dagger$}
\address{
Institut f\"ur Mathematik der Universit\"at Wien\\
Strudlhofgasse 4, A-1090 Wien, Austria}
\email{{\tt Mfulmek@Mat.Univie.Ac.At, Kratt@Pap.Univie.Ac.At}\newline\leavevmode\indent
{\it WWW}: {\tt http://radon.mat.univie.ac.at/\~{}mfulmek}\newline\leavevmode\indent
\hphantom{{\it WWW}: }{\tt http://radon.mat.univie.ac.at/People/kratt}
}

\thanks{$^\dagger$ Research partially supported by the Austrian
Science Foundation FWF, grant P13190-MAT}

\begin{abstract}
We compute the number
of rhombus tilings of a hexagon with side lengths $N,M,N,N,M,N$,
with $N$ and $M$ having the same parity, which
contain a particular rhombus next to the center of the hexagon.
The special case $N=M$ of one of our results solves a problem posed by Propp.
In the proofs, Hankel determinants featuring Bernoulli numbers play
an important role.
\end{abstract}

\maketitle

\section{Introduction}

Let $a$, $b$ and $c$ be positive integers, and consider a 
hexagon with side lengths $a,b,c,a,b,c$ whose angles are $120^\circ$
(see Figure~\ref{fig:0}.a). 
The subject of enumerating rhombus tilings 
of this hexagon (cf\@. Figure~\ref{fig:0}.b;
here, and in the sequel, by a
rhombus we always mean a rhombus with side lengths 1 and angles of
$60^\circ$ and $120^\circ$) gained a lot of interest recently.
This interest comes from two facts. First, it is a rich source of
non-trivial enumeration problems which have (or appear to have)
beautiful solutions (see e.g\@. 
\cite{CiucAH,CiEKAA,CiucKratAB,CiKrAB,CiKrAC,FuKrAC,KupeAA,Propp2,RobbAA}). 
Second, these problems are very often related to 
the theory of symmetric functions and/or the representation theory of
classical and quantum Lie algebras, and to statistical physics
(sometimes in disguise; see e.g\@. 
\cite{CoLPAA,GuOVAA,JohaAA,KupeAF,KupeAG,OkKrAA,ProcAE,ProcAD,ProcAB}).

\begin{figure}
\centertexdraw{
  \drawdim truecm  \linewd.02
  \rhombus \rhombus \rhombus \rhombus \ldreieck
  \move (-0.866025403784439 -.5)
  \rhombus \rhombus \rhombus \rhombus \rhombus \ldreieck
  \move (-1.732050807568877 -1)
  \rhombus \rhombus \rhombus \rhombus \rhombus \rhombus \ldreieck
  \move (-1.732050807568877 -1)
  \rdreieck
  \rhombus \rhombus \rhombus \rhombus \rhombus \rhombus \ldreieck
  \move (-1.732050807568877 -2)
  \rdreieck
  \rhombus \rhombus \rhombus \rhombus \rhombus \rhombus \ldreieck
  \move (-1.732050807568877 -3)
  \rdreieck
  \rhombus \rhombus \rhombus \rhombus \rhombus \rhombus 
  \move (-1.732050807568877 -4)
  \rdreieck
  \rhombus \rhombus \rhombus \rhombus \rhombus 
  \move (-1.732050807568877 -5)
  \rdreieck
  \rhombus \rhombus \rhombus \rhombus 
\move(8 0)
\bsegment
  \drawdim truecm  \linewd.02
  \rhombus \rhombus \rhombus \rhombus \ldreieck
  \move (-0.866025403784439 -.5)
  \rhombus \rhombus \rhombus \rhombus \rhombus \ldreieck
  \move (-1.732050807568877 -1)
  \rhombus \rhombus \rhombus \rhombus \rhombus \rhombus \ldreieck
  \move (-1.732050807568877 -1)
  \rdreieck
  \rhombus \rhombus \rhombus \rhombus \rhombus \rhombus \ldreieck
  \move (-1.732050807568877 -2)
  \rdreieck
  \rhombus \rhombus \rhombus \rhombus \rhombus \rhombus \ldreieck
  \move (-1.732050807568877 -3)
  \rdreieck
  \rhombus \rhombus \rhombus \rhombus \rhombus \rhombus 
  \move (-1.732050807568877 -4)
  \rdreieck
  \rhombus \rhombus \rhombus \rhombus \rhombus 
  \move (-1.732050807568877 -5)
  \rdreieck
  \rhombus \rhombus \rhombus \rhombus 
  \linewd.12
  \move(0 0)
  \RhombusA \RhombusB \RhombusB 
  \RhombusA \RhombusA \RhombusB \RhombusA \RhombusB \RhombusB
  \move (-0.866025403784439 -.5)
  \RhombusA \RhombusB \RhombusB \RhombusB \RhombusB
  \RhombusA \RhombusA \RhombusB \RhombusA 
  \move (-1.732050807568877 -1)
  \RhombusB \RhombusB \RhombusA \RhombusB \RhombusB \RhombusA
  \RhombusB \RhombusA \RhombusA 
  \move (1.732050807568877 0)
  \RhombusC \RhombusC \RhombusC 
  \move (1.732050807568877 -1)
  \RhombusC \RhombusC \RhombusC 
  \move (3.464101615137755 -3)
  \RhombusC 
  \move (-0.866025403784439 -.5)
  \RhombusC
  \move (-0.866025403784439 -1.5)
  \RhombusC
  \move (0.866025403784439 -2.5)
  \RhombusC \RhombusC 
  \move (0.866025403784439 -3.5)
  \RhombusC \RhombusC \RhombusC 
  \move (2.598076211353316 -5.5)
  \RhombusC 
  \move (0.866025403784439 -5.5)
  \RhombusC 
  \move (-1.732050807568877 -3)
  \RhombusC 
  \move (-1.732050807568877 -4)
  \RhombusC 
  \move (-1.732050807568877 -5)
  \RhombusC \RhombusC 
\esegment
\htext (-1.5 -9){\small a. A hexagon with sides $a,b,c,a,b,c$,}
\htext (-1.5 -9.5){\small \hphantom{a. }where $a=3$, $b=4$, $c=5$}
\htext (6.8 -9){\small b. A rhombus tiling of a hexagon}
\htext (6.8 -9.5){\small \hphantom{b. }with sides $a,b,c,a,b,c$}
\rtext td:0 (4.3 -4){$\sideset {}  c 
    {\left.\vbox{\vskip2.2cm}\right\}}$}
\rtext td:60 (2.6 -.5){$\sideset {} {} 
    {\left.\vbox{\vskip1.7cm}\right\}}$}
\rtext td:120 (-.44 -.25){$\sideset {}  {}  
    {\left.\vbox{\vskip1.3cm}\right\}}$}
\rtext td:0 (-2.3 -3.6){$\sideset {c}  {} 
    {\left\{\vbox{\vskip2.2cm}\right.}$}
\rtext td:240 (0 -7){$\sideset {}  {}  
    {\left.\vbox{\vskip1.7cm}\right\}}$}
\rtext td:300 (3.03 -7.25){$\sideset {}  {}  
    {\left.\vbox{\vskip1.4cm}\right\}}$}
\htext (-.9 0.1){$a$}
\htext (2.8 -.1){$b$}
\htext (3.2 -7.8){$a$}
\htext (-0.3 -7.65){$b$}
}
\caption{}
\label{fig:0}
\end{figure}

As is well-known, the total number of rhombus tilings of a hexagon
with side lengths $a,b,c,a,b,c$ equals
\begin{equation}
\label{eq:MacMahon}
\prod_{i=1}^a\prod_{j=1}^b\prod_{k=1}^c\frac{i+j+k-1}{i+j+k-2}.
\end{equation}
(This follows from MacMahon's enumeration 
\cite[Sec.~429, $q\rightarrow 1$; proof in Sec.~494]{MacMahon}
of all plane partitions contained in an $a\times b\times c$ box,
as these are in bijection with rhombus tilings of a hexagon
with side lengths $a,b,c,a,b,c$, as explained e.g\@. in \cite{DT}.)

A natural question to be asked is what the distribution of the
rhombi in a random tiling is. 
On an {\it asymptotic} level, this question was answered by Cohn, 
Larsen and Propp \cite{CoLPAA}. On the {\it exact\/} (enumerative) level, 
the significant contributions are \cite{CiucKratAB,FuKrAC,HeGeAA}.
The most general result was obtained in \cite{FuKrAC} by the authors, where
the number of all rhombus
tilings of a hexagon with side lengths $N,M,N,N,M,N$ was computed,
which contain an arbitrary fixed rhombus
on the symmetry axis which cuts through the sides of length $M$.

The purpose of this paper is to add other results in this direction.
We compute 
the number of all rhombus tilings of a hexagon 
which contain a particular rhombus on the ``other" symmetry axis,
i.e., the symmetry axis which runs {\em in parallel\/} to the sides
of length $M$. In difference to \cite{FuKrAC}, we are not able to
solve this problem for an {\em arbitrary} rhombus on this symmetry
axis, but only for rhombi which are close to the center. In fact, the
case of the central rhombus is already covered by the papers 
\cite{CiucKratAB,FuKrAC,HeGeAA}. We provide results for the ``next"
three cases, i.e., for a rhombus which is by one, two, or three
``half units" off the center, see
Theorems~\ref{thm:MEven}--\ref{thm:NOdd2} below. In fact, as our proofs
show, the computations become increasingly harder, and more
elaborate, as we move the
rhombus gradually farther away from the center, so that it seems highly
unlikely that a uniform formula, similar to the one in \cite{FuKrAC}
for the other symmetry axis,
can be found for an arbitrary rhombus. (For further comments
on this issue see Section~\ref{sec:open}, (2).)

Here are our results.
(Each of the following formulas has to be interpreted as the
appropriate limit if singularities are encountered. For example, if
we directly set $m=1$ in \eqref{eq:NEven2} then the term $(m-1)$
in the second line becomes 0, and, on the other hand, we have a
singularity in the sum in the third line caused by the term $(3-n-m)_h$
in the summand for $h=n-1$. The correct way to interpret the
expression is as the limit as $m$ goes to 1.)

\begin{figure}
\bigskip
\centertexdraw{
  \drawdim truecm  \linewd.02
  \rhombus \rhombus \rhombus \rhombus \ldreieck
  \move (-0.866025403784439 -.5)
  \rhombus \rhombus \rhombus \rhombus \rhombus \ldreieck
  \move (-1.732050807568877 -1)
  \rhombus \rhombus \rhombus \rhombus \rhombus \rhombus
  \move (-2.598076211353316 -1.5)
  \rhombus \rhombus \rhombus \rhombus \rhombus \rhombus
  \move (-2.598076211353316 -1.5)
  \rdreieck
  \rhombus \rhombus \rhombus \rhombus \rhombus
  \move (-2.598076211353316 -2.5)
  \rdreieck
  \rhombus \rhombus \rhombus \rhombus
  \move (0.866025403784439 -1.5)
  \bsegment
    \rlvec(0.866025403784439 -.5) \rlvec(-0.866025403784439 -.5)
	   \rlvec(-0.866025403784439 .5)
    \lfill f:0.8
  \esegment
  \linewd.05
  \move (0.866025403784439 -5.75)
  \rlvec(0 6.5)
\move(8 0)
\bsegment
  \drawdim truecm  \linewd.02
  \rhombus \rhombus \rhombus \ldreieck
  \move (-0.866025403784439 -.5)
  \rhombus \rhombus \rhombus \rhombus \ldreieck
  \move (-1.732050807568877 -1)
  \rhombus \rhombus \rhombus \rhombus \rhombus \ldreieck
  \move (-1.732050807568877 -1)
  \rdreieck \rhombus \rhombus \rhombus \rhombus \rhombus
  \move (-1.732050807568877 -2)
  \rdreieck \rhombus \rhombus \rhombus \rhombus
  \move (-1.732050807568877 -3)
  \rdreieck \rhombus \rhombus \rhombus
  \move (0.866025403784439 -1.5)
  \bsegment
    \rlvec(0.866025403784439 -.5) \rlvec(-0.866025403784439 -.5)
	   \rlvec(-0.866025403784439 .5)
    \lfill f:0.8
  \esegment
  \linewd.05
  \move (0.866025403784439 -5.75)
  \rlvec(0 6.5)
\esegment
\htext (-2.5 -6.4){\small a.}
\htext (-2 -6.5){\small A hexagon with sides $N,M,N,$}
\htext (-2 -7){\small $N,M,N$, where $N=4$, $M=2$.}
\htext (5.3 -6.4){\small b.}
\htext (5.8 -6.5){\small A hexagon with sides $N,M,N,$}
\htext (5.8 -7){\small $N,M,N$, where $N=3$, $M=3$.}
\htext (-1.4 -0.1){$N$}
\htext (2.9 -0.1){$N$}
\htext (4.5 -2.6){$M$}
\htext (7.1 0.1){$N$}
\htext (10.4 0.1){$N$}
\htext (5.7 -2.6){$M$}
}
\caption{}
\label{fig:1}
\end{figure}

\begin{thm}
\label{thm:MEven}
Let $n$ and $m$ be positive integers. The number
of rhombus tilings of a hexagon with side lengths $2n,2m,2n,2n,2m,2n$, which
contain the rhombus above and next to the center of the hexagon
{\em (}see
Figure~\ref{fig:1}.a; the rhombus which is contained in every
tiling is shaded{\em)}, equals
\begin{multline}
\label{eq:MEven}
\frac{
n m \binom{2n}{n}\binom{2n-1}{n}\binom{2m}{m}
}{
\binom{4n+2m-1}{2n+m}
}
\Biggl(
	-\frac{1}{(n+m)^2}
	+\frac{4n+2}{(n+1)(2n-1)(n+m-1)(n+m+1)}\\
	\times\sum_{h=0}^{n-1}
		\frac{
			\po{2}{h}\,\po{1-n}{h}\,\po{\frac {3} {2}+n}{h}\,\po{1-n-m}{h}\,\po{1+n+m}{h}
		}{
			\po{1}{h}\,\po{2+n}{h}\,\po{\frac {3} {2}-n}{h}\,\po{2+n+m}{h}\,\po{2-n-m}{h}
		}
\Biggr)\\
\times
\prod_{i=1}^{2n}\prod_{j=1}^{2m}\prod_{k=1}^{2n}\frac{i+j+k-1}{i+j+k-2}.
\end{multline}
where the shifted factorial
$(a)_k$ is defined by $(a)_k:=a(a+1)\cdots(a+k-1)$,
$k\ge1$, $(a)_0:=1$.
\end{thm}

\begin{thm}
\label{thm:MOdd}
Let $n$ be a nonnegative integer and $m$ be a positive integer. The number
of rhombus tilings of a hexagon with side lengths
$2n+1,2m-1,2n+1,2n+1,2m-1,2n+1$, which
contain the rhombus above and next to the center of the hexagon
{\em (}see
Figure~\ref{fig:1}.b; the rhombus which is contained in every
tiling is shaded{\em)}, equals
\begin{multline}
\label{eq:MOdd}
\frac{
(n+1) m \binom{2n}{n} \binom{2n+1}{n} \binom{2m-1}{m}
}{
\binom{4n+2m}{2n+m}
}
\Biggl(
	\frac{1}{(n+m)^2}
	+\frac{4n}{(n+1)(2n-1)(n+m-1)(n+m+1)}\\
\times	\sum_{h=0}^{n-1}
		\frac{
			\po{2}{h}\,\po{1-n}{h}\,\po{\frac {3} {2}+n}{h}\,\po{1-n-m}{h}\,\po{1+n+m}{h}
		}{
			\po{1}{h}\,\po{2+n}{h}\,\po{\frac {3} {2}-n}{h}\,\po{2+n+m}{h}\,\po{2-n-m}{h}
		}
\Biggr)\\
\times
\prod_{i=1}^{2n+1}\prod_{j=1}^{2m-1}\prod_{k=1}^{2n+1}\frac{i+j+k-1}{i+j+k-2}.
\end{multline}
\end{thm}

\begin{thm}
\label{thm:NEven}
Let $n$ and $m$ be positive integers. The number
of rhombus tilings of a hexagon with side lengths $2n,2m-1,2n,2n,2m-1,2n$,
which contain the rhombus above and next to the central rhombus
(see Figure~\ref{fig:1a}.a; the rhombus which is contained in
every tiling is shaded), equals
\begin{multline}
\label{eq:NEven}
\frac{
	(2 m-1) \binom{2 m-2}{m-1} \binom{2 n-4}{n-2} 
	\binom{2 n+2}{n+1}
}{
	\left(n + m - 1\right)  
	\left(n + m \right)  \binom{4 n + 2 m - 3}{2 n + m - 2}
}
\\
\times
\Biggl(
	\frac{
		n\,\left(n + 1\right)\,\left(2 n - 3\right) \,
		\left(2 n - 1\right) \,
		\left({m^2} - m + n + 2 m\,n + {n^2} + 1\right)
	}{
		(n - 1)\,\left(n + m - 1\right) \,\left(n + m \right)
		\,\left(2 n + 1\right)
	} + \kern3.5cm\\
        \frac {6} {(n+m-2)(n+m+1)}
		\sum_{h = 0}^{n - 2}
			\frac{
				\po{3}{h}\,
				\po{{\frac{5}{2}}}{h}\,
				\po{2 - n}{h}\,
				\po{{\frac{3}{2}} + n}{h}\,
				\po{2 - n - m}{h}\,
				\po{1 + n + m}{h}
			}{
				\po{1}{h}\,
				\po{{\frac{3}{2}}}{h}\,
				\po{2 + n}{h}\,
				\po{{\frac{5}{2}} - n}{h}\,
				\po{2 + n + m}{h}\,
				\po{3 - n - m}{h}
			}
\Biggr)\\
\times
\prod_{i=1}^{2n}\prod_{j=1}^{2m-1}\prod_{k=1}^{2n}\frac{i+j+k-1}{i+j+k-2}.
\end{multline}

\end{thm}

\begin{thm}
\label{thm:NOdd}
Let $n$ and $m$ be positive integers. The number
of rhombus tilings of a hexagon with side lengths
$2n-1,2m,2n-1,2n-1,2m,2n-1$, which
contain the rhombus above and next to the central rhombus 
(see Figure~\ref{fig:1a}.b; the rhombus which is contained in
every tiling is shaded), equals
\begin{multline}
\label{eq:NOdd}
\frac{
	(2 m - 1) \binom{2 m - 2}{m - 1} \binom{2 n - 4}{n - 2} 
	\binom{2 n + 2}{n + 1}
}{
	\left(n + m - 1\right)  
	\left( n + m \right)  \binom{4 n + 2 m - 3}{2 n + m - 2}
}
\\
\times
\Biggl(
	\frac{
		n\,\left(n + 1\right)\,\left(2 n - 3\right) \,
		\left(2 n - 1\right) \,
		\left({m^2} - m - 3n + 2 m\,n + {n^2} + 2\right)
	}{
		(n - 1)\,\left(n + m - 1\right) \,\left(n + m \right)
		\,\left(2 n + 1\right)
	} + \kern3.5cm\\
        \frac {6} {(n+m-2)(n+m+1)}
		\sum_{h = 0}^{n - 2}
			\frac{
				\po{3}{h}\,
				\po{{\frac{5}{2}}}{h}\,
				\po{2 - n}{h}\,
				\po{{\frac{3}{2}} + n}{h}\,
				\po{2 - n - m}{h}\,
				\po{1 + n + m}{h}
			}{
				\po{1}{h}\,
				\po{{\frac{3}{2}}}{h}\,
				\po{2 + n}{h}\,
				\po{{\frac{5}{2}} - n}{h}\,
				\po{2 + n + m}{h}\,
				\po{3 - n - m}{h}
			}
\Biggr)\\
\times
\prod_{i=1}^{2n-1}\prod_{j=1}^{2m}\prod_{k=1}^{2n-1}\frac{i+j+k-1}{i+j+k-2}.
\end{multline}
\end{thm}

\begin{figure}
\bigskip
\centertexdraw{
  \drawdim truecm  \linewd.02
  \rhombus \rhombus \ldreieck
  \move (-0.866025403784439 -.5)
  \rhombus \rhombus \rhombus \ldreieck
  \move (-0.866025403784439 -.5)
  \rdreieck \rhombus \rhombus \rhombus \ldreieck
  \move (-0.866025403784439 -1.5)
  \rdreieck \rhombus \rhombus \rhombus
  \move (-0.866025403784439  -2.5)
  \rdreieck \rhombus \rhombus
  \move (0.866025403784439 -.5)
  \bsegment
    \rlvec(0.866025403784439 -.5) \rlvec(-0.866025403784439 -.5)
	   \rlvec(-0.866025403784439 .5)
    \lfill f:0.8
  \esegment
  \linewd.05
  \move (0.866025403784439 -4.75)
  \rlvec(0 5.5)
\move(8 0)
\bsegment
  \drawdim truecm  \linewd.02
  \rhombus \rhombus \rhombus \ldreieck
  \move (-0.866025403784439 -.5)
  \rhombus \rhombus \rhombus \rhombus \ldreieck
  \move (-1.732050807568877 -1)
  \rhombus \rhombus \rhombus \rhombus \rhombus
  \move (-1.732050807568877 -1)
  \rdreieck
  \rhombus \rhombus \rhombus \rhombus
  \move (-1.732050807568877 -2)
  \rdreieck
  \rhombus \rhombus \rhombus
  \move (0.866025403784439 -.5)
  \bsegment
    \rlvec(0.866025403784439 -.5) \rlvec(-0.866025403784439 -.5)
	   \rlvec(-0.866025403784439 .5)
    \lfill f:0.8
  \esegment
  \linewd.05
  \move (0.866025403784439 -4.75)
  \rlvec(0 5.5)
\esegment
\htext (-2.5 -5.4){\small a.}
\htext (-2 -5.5){\small A hexagon with sides $N,M,N,$}
\htext (-2 -6){\small $N,M,N$, where $N=2$, $M=3$.}
\htext (5.3 -5.4){\small b.}
\htext (5.8 -5.5){\small A hexagon with sides $N,M,N,$}
\htext (5.8 -6){\small $N,M,N$, where $N=3$, $M=2$.}
%
\htext (-0.6 0.1){$N$}
\htext (1.9 0.1){$N$}
\htext (3.0 -2.1){$M$}
\htext (7.0 -0.2){$N$}
\htext (10.4 -0.2){$N$}
\htext (5.5 -2.1){$M$}
}
\caption{}
\label{fig:1a}
\end{figure}

\begin{thm}
\label{thm:NEven2}
Let $n$ and $m$ be positive integers. The number
of rhombus tilings of a hexagon with side lengths $2n,2m,2n,2n,2m,2n$,
which contain the rhombus above and next to the rhombus which is adjacent
to the center of the hexagon
(see Figure~\ref{fig:1b}.a; the rhombus which is contained in
every tiling is shaded), equals
\begin{multline}
\label{eq:NEven2}
\frac{
	\binom{n+m-1}{m} 
	\binom{2n+2}{n-1}
	\binom{2n+m-1}{n}
	\binom{2n+m}{2n+1}
}{
	2(2n-3)(2n-1)(2n+2)(n+m-1)(n+m+1)
	\binom{4n-1}{2n}
	\binom{4n+2m-1}{2m}
}
\\
\times\Biggl( 
{\frac {2   \left( n+2 \right)  \left( n+3 \right)  
     \left( 2 n -1\right) \left( 2 n-3 \right) X(m,n) } 
   {\left( n+m -1\right)  {{\left( n+m \right) }^2} 
     \left( n+m+1 \right) }}
	-
	\frac{
		24(m-1)(2n+m+1)(2n+1)(2n+3)
	}{
		(n+m-2)(n+m+2)
	}\\
	\cdot\sum_{h = 0}^{n - 1}
			\frac{
				\po{4}{h}\,
				\po{1 - n}{h}\,
				\po{{\frac{5}{2}} + n}{h}\,
				\po{2 - n - m}{h}\,
				\po{2 + n + m}{h}
			}{
				\po{1}{h}\,
				\po{4 + n}{h}\,
				\po{{\frac{5}{2}} - n}{h}\,
				\po{3 + n + m}{h}\,
				\po{3 - n - m}{h}
			}
\Biggr)
\prod_{i=1}^{2n}\prod_{j=1}^{2m}\prod_{k=1}^{2n}\frac{i+j+k-1}{i+j+k-2},
\end{multline}
where 

\vskip-10pt
{\small
\begin{multline*} 
X(m,n)= -1 + 4 {m^2} - 3 {m^4} - n + 8 m n + 11 {m^2} n - 12 {m^3} n - 
  6 {m^4} n + 8 {n^2} + 22 m {n^2} - 4 {m^2} {n^2} - 
  24 {m^3} {n^2}\\
 - 2 {m^4} {n^2} + 15 {n^3} + 16 m {n^3} - 
  28 {m^2} {n^3} - 8 {m^3} {n^3} - 5 {n^4} - 8 m {n^4} - 
  12 {m^2} {n^4} - 14 {n^5} - 8 m {n^5} - 2 {n^6}.
\end{multline*}
}
\end{thm}

\begin{thm}
\label{thm:NOdd2}
Let $n$ and $m$ be positive integers. The number
of rhombus tilings of a hexagon with side lengths
$2n+1,2m-1,2n+1,2n+1,2m-1,2n+1$, which
contain the rhombus above and next to the rhombus which is adjacent
to the center of the hexagon
(see Figure~\ref{fig:1b}.b; the rhombus which is contained in
every tiling is shaded), equals
\begin{multline}
\label{eq:NOdd2}
\frac{
	\binom{2m-1}{m}
	\binom{n+m-1}{m-1}
	\binom{2n+2}{n-1}
	\binom{2n+m}{n+1}
}{
	(m+1)(m+2)(2n-3)(n+m-1)(n+m+1)
	\binom{2n+m}{m+2}
	\binom{4n+2m}{2n+m}
}
\\
\times\Biggl(
{\frac { \left( n+2 \right)  \left( n+3 \right)  
     \left( 2 n-1 \right) \left( 2 n-3 \right) Y(m,n) }
   {\left( n + m-1 \right)  {{\left( n + m \right) }^2} 
     \left( n + m+1 \right) }}
	-
	\frac{
		24(m-1)n(2n+m+1)(2n+3)
	}{
		(n+m-2)(n+m+2)
	}\\
	\cdot\sum_{h = 0}^{n - 1}
			\frac{
				\po{4}{h}\,
				\po{1 - n}{h}\,
				\po{{\frac{5}{2}} + n}{h}\,
				\po{2 - n - m}{h}\,
				\po{2 + n + m}{h}
			}{
				\po{1}{h}\,
				\po{4 + n}{h}\,
				\po{{\frac{5}{2}} - n}{h}\,
				\po{3 + n + m}{h}\,
				\po{3 - n - m}{h}
			}
\Biggr)
\prod_{i=1}^{2n+1}\prod_{j=1}^{2m-1}\prod_{k=1}^{2n+1}\frac{i+j+k-1}{i+j+k-2},
\end{multline}
where 

\vskip-10pt
{\small
\begin{multline*} 
Y(m,n)= 1 - 2 {m^2} + {m^4} + 5 n - 4 m n - 2 {m^2} n + 
       4 {m^3} n - 3 {m^4} n + 22 {n^2} - 4 m {n^2} - 
       4 {m^2} {n^2} \\
- 12 {m^3} {n^2} - 2 {m^4} {n^2} + 50 {n^3} - 
       16 m {n^3} - 26 {m^2} {n^3} - 8 {m^3} {n^3} + 39 {n^4} - 
       28 m {n^4} - 12 {m^2} {n^4} + 5 {n^5} - 8 m {n^5} - 2 {n^6}.
\end{multline*}
}
\end{thm}

\begin{figure}
\caption{}
\label{fig:1b}
\bigskip
\centertexdraw{
  \drawdim truecm  \linewd.02
  \rhombus \rhombus \rhombus \rhombus \ldreieck
  \move (-0.866025403784439 -.5)
  \rhombus \rhombus \rhombus \rhombus \rhombus \ldreieck
  \move (-1.732050807568877 -1)
  \rhombus \rhombus \rhombus \rhombus \rhombus \rhombus
  \move (-2.598076211353316 -1.5)
  \rhombus \rhombus \rhombus \rhombus \rhombus \rhombus
  \move (-2.598076211353316 -1.5)
  \rdreieck
  \rhombus \rhombus \rhombus \rhombus \rhombus
  \move (-2.598076211353316 -2.5)
  \rdreieck
  \rhombus \rhombus \rhombus \rhombus
  \move (0.866025403784439 -0.5)
  \bsegment
    \rlvec(0.866025403784439 -.5) \rlvec(-0.866025403784439 -.5)
	   \rlvec(-0.866025403784439 .5)
    \lfill f:0.8
  \esegment
  \linewd.05
  \move (0.866025403784439 -5.75)
  \rlvec(0 6.5)
\move(8 0)
\bsegment
  \drawdim truecm  \linewd.02
  \rhombus \rhombus \rhombus \ldreieck
  \move (-0.866025403784439 -.5)
  \rhombus \rhombus \rhombus \rhombus \ldreieck
  \move (-1.732050807568877 -1)
  \rhombus \rhombus \rhombus \rhombus \rhombus \ldreieck
  \move (-1.732050807568877 -1)
  \rdreieck \rhombus \rhombus \rhombus \rhombus \rhombus
  \move (-1.732050807568877 -2)
  \rdreieck \rhombus \rhombus \rhombus \rhombus
  \move (-1.732050807568877 -3)
  \rdreieck \rhombus \rhombus \rhombus
  \move (0.866025403784439 -0.5)
  \bsegment
    \rlvec(0.866025403784439 -.5) \rlvec(-0.866025403784439 -.5)
	   \rlvec(-0.866025403784439 .5)
    \lfill f:0.8
  \esegment
  \linewd.05
  \move (0.866025403784439 -5.75)
  \rlvec(0 6.5)
\esegment
\htext (-2.5 -6.4){\small a.}
\htext (-2 -6.5){\small A hexagon with sides $N,M,N,$}
\htext (-2 -7){\small $N,M,N$, where $N=4$, $M=2$.}
\htext (5.3 -6.4){\small b.}
\htext (5.8 -6.5){\small A hexagon with sides $N,M,N,$}
\htext (5.8 -7){\small $N,M,N$, where $N=3$, $M=3$.}
\htext (-1.4 -0.1){$N$}
\htext (2.9 -0.1){$N$}
\htext (4.5 -2.6){$M$}
\htext (7.1 0.1){$N$}
\htext (10.4 0.1){$N$}
\htext (5.7 -2.6){$M$}
}
\end{figure}

In general, the sum in \eqref{eq:MEven} and \eqref{eq:MOdd} (note
that it is indeed exactly the same sum), the sum in 
\eqref{eq:NEven} and \eqref{eq:NOdd} 
(it is indeed exactly the same sum), and the sum in 
\eqref{eq:NEven2} and \eqref{eq:NOdd2} (again, 
it is indeed exactly the same sum), does not simplify. 
However, in the case that $n$ and $m$ are roughly of the same size,
the sum does simplify. For the sake of brevity, we give here just a
sample of corollaries to Theorems~\ref{thm:MEven} and \ref{thm:MOdd},
the first two statements of which solve the second part of Problem~4 in
Propp's list \cite{PropAA}. We wish to emphasize that there are
similar corollaries to Theorems~\ref{thm:NEven}--\ref{thm:NOdd2}.

\begin{cor}\label{thm:Propp}
Let $n$ be a nonnegative integer. The number of 
rhombus tilings of a hexagon with all sides of length $2n$,
which contain the rhombus above and next to the center of the hexagon, equals
\begin{equation}
\label{eq:ProppEven}
\left(
	\frac{1}{3}-
	\frac{1}{12}
	\frac{
		\binom{2n}{n}^3
	}{
		\binom{6n}{3n}
	}
\right)
\prod_{i=1}^{2n}\prod_{j=1}^{2n}\prod_{k=1}^{2n}\frac{i+j+k-1}{i+j+k-2}.
\end{equation}
The number of 
rhombus tilings of a hexagon with all sides of length $2n+1$,
which contain the rhombus above and next to the center of the hexagon, equals
\begin{equation}
\label{eq:ProppOdd}
\left(
	\frac{1}{3}+
	\frac{1}{3}
	\frac{
		\binom{2n}{n}^3
	}{
		\binom{6n+2}{3n+1}
	}
\right)
\prod_{i=1}^{2n+1}\prod_{j=1}^{2n+1}\prod_{k=1}^{2n+1}\frac{i+j+k-1}{i+j+k-2}.
\end{equation}
The number of 
rhombus tilings of a hexagon with side lengths $2n,2n+2,2n,2n,2n+2,2n$,
which contain the rhombus above and next to the center of the hexagon, equals
\begin{equation}
\label{eq:ProppEven1}
\left(
	\frac{1}{3}-
	\frac{(10n+2)}{(6n+3)}
	\frac{
		\binom{2n}{n}^3
	}{
		\binom{6n+2}{3n+1}
	}
\right)
\prod_{i=1}^{2n}\prod_{j=1}^{2n+2}\prod_{k=1}^{2n}\frac{i+j+k-1}{i+j+k-2}.
\end{equation}
The number of 
rhombus tilings of a hexagon with side lengths $2n+1,2n-1,2n+1,2n+1,2n-1,2n+1$,
which contain the rhombus above and next to the center of the hexagon, equals
\begin{equation}
\label{eq:ProppOdd1}
\left(
	\frac{1}{3}+
	\frac{(10n+3)}{24n}
	\frac{
		\binom{2n}{n}^3
	}{
		\binom{6n}{3n}
	}
\right)
\prod_{i=1}^{2n+1}\prod_{j=1}^{2n-1}\prod_{k=1}^{2n+1}\frac{i+j+k-1}{i+j+k-2}.
\end{equation}
The number of 
rhombus tilings of a hexagon with side lengths $2n+2,2n,2n+2,2n+2,2n,2n+2$,
which contain the rhombus above and next to the center of the hexagon, equals
\begin{equation}
\label{eq:ProppEven2}
\left(
	\frac{1}{3}+
	4\frac{
		\binom{2n}{n}^3
	}{
		\binom{6n+4}{3n+2}
	}
\right)
\prod_{i=1}^{2n+2}\prod_{j=1}^{2n}\prod_{k=1}^{2n+2}\frac{i+j+k-1}{i+j+k-2}.
\end{equation}
The number of 
rhombus tilings of a hexagon with side lengths $2n+3,2n-1,2n+3,2n+3,2n-1,2n+3$,
which contain the rhombus above and next to the center of the hexagon, equals
\begin{equation}
\label{eq:ProppOdd2}
\left(
	\frac{1}{3}+
	\frac{2(6n^2+9n+2)}{(n+1)^2}
	\frac{
		\binom{2n}{n}^3
	}{
		\binom{6n+4}{3n+2}
	}
\right)
\prod_{i=1}^{2n+3}\prod_{j=1}^{2n-1}\prod_{k=1}^{2n+3}\frac{i+j+k-1}{i+j+k-2}.
\end{equation}
\end{cor}

If $m$ is not $n+1$, $n$ or $n-1$, then one may still try to obtain at
least estimates for the number of rhombus tilings that contain this
particular rhombus.
Indeed,
from Theorems~\ref{thm:MEven}--\ref{thm:NOdd2}, we are able to derive
an ``arcsine law" for this kind of enumeration, which is analogous
to the ones in \cite{CiucKratAB} and \cite{FuKrAC}.

\begin{cor}\label{thm:asy}
Let $a$ be any nonnegative real number, and consider a hexagon with
side lengths $2n,2m,2n,2n,2m,2n$.
For $m\sim an$, the proportion
of the rhombus tilings that contain the
rhombus above and next to the center of the hexagon
in the total number of rhombus tilings is
$\sim\ \frac {2} {\pi}\arcsin(1/(a+1))$
as $n$ tends to infinity.
The same is true for a
hexagon with side lengths $2n+1,2m-1,2n+1,2n+1,2m-1,2n+1$. Moreover,
for hexagons 
with side lengths $2n,2m-1,2n,2n,2m-1,2n$ or hexagons
with side lengths $2n-1,2m,2n-1,2n-1,2m,2n-1$,
the same result holds for the proportion 
of the rhombus tilings that contain 
the rhombus above and next to the central rhombus
in the total number of rhombus tilings, as well as 
for the proportion of the rhombus tilings that contain 
the rhombus above and next to the rhombus which is adjacent to the
center of the hexagon
in the total number of rhombus tilings of hexagons 
with side lengths $2n,2m,2n,2n,2m,2n$ or hexagons
with side lengths $2n+1,2m-1,2n+1,2n+1,2m-1,2n+1$.
\end{cor}
Also this result is (as well as Corollary~4 in \cite{CiucKratAB}
and Theorem~1.3 in \cite{FuKrAC}; see the
respective comments in \cite{FuKrAC}) 
in accordance with Conjecture~1 in
\cite{CoLPAA}, to which it adds evidence in further special instances.

In the next sections we describe proofs of 
Theorems~\ref{thm:MEven}--\ref{thm:NOdd2},
and of 
Corollaries~\ref{thm:Propp} and \ref{thm:asy}.
In Section~2 we provide proofs of Corollaries~\ref{thm:Propp} and
\ref{thm:asy},
and we outline the proofs of
Theorems~\ref{thm:MEven}--\ref{thm:NOdd2},
the latter consisting 
of two basic steps. In the first step we
build on the approach of Helfgott and Gessel in 
\cite{HeGeAA}, a short summary of which is the contents of Section~3. 
It allows to write the
number that we are interested in 
in form of a determinant. The evaluation
of this determinant is not easy and is carried out in detail in
Section~4. For the evaluation we follow a ``method" that was
first introduced in \cite{KratBD} (see the tutorial description in
\cite[Sec.~2.4]{KratBN} or \cite[Sec.~2]{KratBI}). 
For accomplishing the required computations, we need
to evaluate certain Hankel determinants featuring Bernoulli numbers,
which are, in fact, of independent interest.
As it turns out, some of the evaluations of these Hankel determinants are
already known, provided certain results about
orthogonal polynomials, in particular, about
continuous Hahn polynomials, and continued fractions are properly combined. 
For the convenience of the reader,
we collect these facts, and their implications, in Section~5.
In particular, the evaluation of the relevant Hankel determinants
featuring Bernoulli numbers is given in
Theorem~\ref{thm:HankelBernoulli}.
However, in the proofs of Theorems~\ref{thm:NEven2} and \ref{thm:NOdd2}
(more precisely, in the proof of the subordinate Lemma~\ref{lem:Det3})
we encounter a certain Hankel determinant of Bernoulli numbers
(see \eqref{eq:extra-Bernoulli}), the evaluation of which requires
considerable effort. (This is one of the added difficulties mentioned
earlier in comparison to the proofs of Theorems~\ref{thm:MEven} and
\ref{thm:MOdd}.) We evaluate this Hankel determinant by combining the
knowledge about continuous Hahn polynomials with a recent theorem on
orthogonal polynomials due to Leclerc \cite{Leclerc} (restated here
as Theorem~\ref{thm:Leclerc}), and applying some integral calculus
(see the proof of Lemma~\ref{lem:weird}). Section~6 is devoted to
provide the details of these calculations. In
Section~\ref{sec:further} we make
explicit a few unusual evaluations of Hankel determinants of
Bernoulli polynomials, which are implicit in the proofs
of our enumeration results.
Finally, in Section~\ref{sec:open}, we point to further directions in this
research, and propose a few open problems.


\section{Outline of proofs}
Here we outline the proofs of Theorems~\ref{thm:MEven}--\ref{thm:NOdd2},
and we deduce Corollaries~\ref{thm:Propp} and \ref{thm:asy}. 
We fill in the details in the subsequent sections.

\smallskip
{\sc Proof of Theorems~\ref{thm:MEven}--\ref{thm:NOdd2}}.
Following the approach of Helfgott and Gessel \cite{HeGeAA} (see
Section~3), we may
write the number of rhombus tilings of a hexagon with side lengths
$N,M,N,N,M,N$, which
contain an {\em arbitrary\/} rhombus on the $(N+M)$-long vertical symmetry
axis (see Figure~\ref{figure-thm6}),
in form of a determinant. This determinant is given by
Proposition~\ref{thm:HeGe-det}.
That is to say, in
order to prove Theorems~\ref{thm:MEven} and \ref{thm:MOdd}, we need
to evaluate the determinant in \eqref{eq:HeGe-det} with
$l=\frac{N+M}{2}$, and in order to prove Theorems~\ref{thm:NEven} and \ref{thm:NOdd}, we need to evaluate the same determinant with
$l=\frac{N+M+1}{2}$, and in order to prove Theorems~\ref{thm:NEven2} and \ref{thm:NOdd2}, we need to evaluate the same determinant with
$l=\frac{N+M+2}{2}$. Modulo replacement of parameters, we may thus
concentrate on the determinants $D(n,n,N)$, $D(n,n-1,N)$ and
$D(n,n-2,N)$, where
\begin{equation}
\label{eq:gen-det}
D(a,b,N):=\det_{1\le i,j\le  N}\bigg(\sum_{s=-a}^{b-1} s^{i+j}\bigg).
\end{equation}

The evaluation of 
determinants $D(n,n,N)$, $D(n,n-1,N)$ and $D(n,n-2,N)$ is then carried out in
Lemmas~\ref{lem:Det1}, \ref{lem:Det2} and \ref{lem:Det3}, respectively.

Theorem~\ref{thm:MEven} follows upon combining
Proposition~\ref{thm:HeGe-det} with $N=2n$, $M=2m$, $l=n+m$ and 
Lemma~\ref{lem:Det1} with
$n$ replaced by $n+m$ and $N=2n-1$, and some rearrangement of terms.
Similarly, Theorem~\ref{thm:MOdd} follows upon combining
Proposition~\ref{thm:HeGe-det} with $N=2n+1$, $M=2m-1$, $l=n+m$ and 
Lemma~\ref{lem:Det1} with $n$ replaced by $n+m$ and $N=2n$.

Likewise, Theorem~\ref{thm:NEven} follows upon combining
Proposition~\ref{thm:HeGe-det} with $N=2n$, $M=2m-1$, $l=n+m$ and 
\eqref{eq:Det2odd} with
$n$ replaced by $n+m$ and $m$ replaced by $n$, and
Theorem~\ref{thm:NOdd} follows upon combining
Proposition~\ref{thm:HeGe-det} with $N=2n-1$, $M=2m$, $l=n+m$ and 
\eqref{eq:Det2even} with $n$ replaced by $n+m$ and $m$ replaced by $n-1$.

Likewise, Theorem~\ref{thm:NEven2} follows upon combining
Proposition~\ref{thm:HeGe-det} with $N=2n$, $M=2m$, $l=n+m+1$ and 
\eqref{eq:Det3odd} with
$n$ replaced by $n+m+1$ and $m$ replaced by $n-1$, and
Theorem~\ref{thm:NOdd2} follows upon combining
Proposition~\ref{thm:HeGe-det} with $N=2n+1$, $M=2m-1$, $l=n+m+1$ and 
\eqref{eq:Det3even} with $n$ replaced by $n+m+1$ and $m$ replaced by $n$.
\hfil\qed

\bigskip
{\sc Proof of Corollary~\ref{thm:Propp}}.
We have to compute the value of the expressions \eqref{eq:MEven}
and \eqref{eq:MOdd} for $m=n+1$ (in order to establish
\eqref{eq:ProppEven1} and \eqref{eq:ProppOdd}), 
$m=n$ (in order to establish
\eqref{eq:ProppEven} and \eqref{eq:ProppOdd1}), 
and $m=n-1$ (in order to establish
\eqref{eq:ProppEven2} and \eqref{eq:ProppOdd2}).
Clearly, except for
trivial manipulations, we will be done once we are able to evaluate
the sum in \eqref{eq:MEven} and \eqref{eq:MOdd} (it is indeed the
same sum!) for $m=n+1$, $m=n$, respectively $m=n-1$.

We treat the case $m=n$ first.
We claim that 
\begin{equation} \label{eq:Summe}
	\sum_{h=0}^{n-1}
		\frac{
		\po{2}{h}\,\po{3/2+n}{h}\,\po{1-n}{h}\,\po{1+2n}{h}\,\po{1-2n}{h}
		}{
		\po{1}{h}\,\po{3/2-n}{h}\,\po{2+n}{h}\,\po{2-2n}{h}\,\po{2+2n}{h}
		}
=\frac {(n+1)(2n-1)^2} {2n^2}\left(\frac {1} {6}+\frac {1} {3}
\frac {\binom {6n}{3n}} {{\binom {2n}{n}}^3}\right).
\end{equation}
Let us denote the sum by $S(n)$ and its summand by $F(n,i)$.
We use the Gosper--Zeilberger algorithm \cite{PeWZAA,ZeilAM,ZeilAV}
to obtain the relation
\begin{multline} \label{eq:WZ}
6 n^2  (n+2) ( 6 n+1) ( 6 n+5) \,F(n,i) - 
     6 (n+1) (2 n-1)^2  (3 n+1) ( 3 n+2) \,F(n+1,i) \\
=G(n,i+1)-G(n,i),
\end{multline}
with
\begin{multline*} 
G(n,i)=\frac {n (n+2) (2n-2h-1) (2n-h-1) 
    }
  {(h+1) (h - n) (2n-h) (2n - h + 1) (2n + h + 2)}\\
\times(
 144 n^5
- 432 h^2 n^4 
- 432 h n^4 
+ 312 n^4 
- 936 h^2 n^3 
- 936 h n^3 
+ 236 n^3 
+ 108 h^4 n^2 \\
\quad \quad + 216 h^3 n^2 
- 588 h^2 n^2 
- 696 h n^2 
+ 70 n^2 
+ 117 h^4 n 
+ 234 h^3 n 
- 83 h^2 n 
- 200 h n \\
+ 6 n 
+ 24 h^4 
+ 48 h^3 
+ 6 h^2 
-18 h 
)F(n,i).
\end{multline*}
Summation of the relation \eqref{eq:WZ} from $i=0$ to $i=n-1$, and little
rearrangement, leads to the recurrence
\begin{multline*}
6 n^2  (n+2) ( 6 n+1) ( 6 n+5) \,S(n) - 
     6 (n+1) (2 n-1)^2  (3 n+1) ( 3 n+2) \,S(n+1) \\
=\frac{( n+2) (2 n-1)^2 (36 n^3 + 60 n^2 + 29 n + 3)}{2 (n+1)}.
\end{multline*}
for the sum in \eqref{eq:Summe}. (Paule and Schorn's \cite{PaScAA} {\sl Mathematica}
implementation of the Gosper--Zeilberger algorithm,
which is the one we used, gives this recurrence directly.) Since
$S(1)=1$, and since the right-hand side of \eqref{eq:Summe} satisfies the same
recurrence, equation \eqref{eq:Summe} is proved.

The procedure in the other two cases is analogous. The respective
evaluations that we need to prove are
\begin{equation} \label{eq:Summe1}
	\sum_{h=0}^{n-1}
		\frac{
		\po{2}{h}\,\po{3/2+n}{h}\,\po{1-n}{h}\,\po{2+2n}{h}\,\po{-2n}{h}
		}{
	\po{1}{h}\,\po{3/2-n}{h}\,\po{2+n}{h}\,\po{1-2n}{h}\,\po{3+2n}{h}
		}
=\frac {(n+1)^2(2n-1)} {(2n+1)^2}\left(-\frac {2} {3}+\frac {1} {3}
\frac {\binom {6n+2}{3n+1}} {{\binom {2n}{n}}^3}\right)
\end{equation}
in the case that $m=n+1$, and
\begin{equation} \label{eq:Summe2}
	\sum_{h=0}^{n-1}
		\frac{
	\po{2}{h}\,\po{3/2+n}{h}\,\po{1-n}{h}\,\po{2n}{h}\,\po{2-2n}{h}
		}{
	\po{1}{h}\,\po{3/2-n}{h}\,\po{2+n}{h}\,\po{3-2n}{h}\,\po{1+2n}{h}
		}
=\frac {2n(n+1)} {2n+1}\left({1}+\frac {n} {12(2n-1)}
\frac {\binom {6n-2}{3n-1}} {{\binom {2n-2}{n-1}}^3}\right)
\end{equation}
in the case that $m=n-1$. We leave it to the reader to fill in the
details.
\hfil\qed

\bigskip
{\sc Proof of Corollary~\ref{thm:asy}}.
We concentrate first on the case of a hexagon with side lengths 
$2n,2m,2n,2n,2m,2n$ and the rhombus above and next to the center 
of the hexagon. 

{}From MacMahon's formula \eqref{eq:MacMahon} for the total number of rhombus
tilings together with Theorem~\ref{thm:MEven} we infer
that the proportion of the rhombus tilings that contain the
rhombus above and next to the center of the hexagon with sides
lengths $2n,2m,2n,2n,2m,2n$
in the total number of rhombus tilings is given by
\begin{multline*}
\frac{
n m \binom{2n}{n}\binom{2n-1}{n}\binom{2m}{m}
}{
\binom{4n+2m-1}{2n+m}
}
\Biggl(
	-\frac{1}{(n+m)^2}
	+\frac{4n+2}{(n+1)(2n-1)(n+m-1)(n+m+1)}\\
	\times\sum_{h=0}^{n-1}
		\frac{
			\po{2}{h}\,\po{3/2+n}{h}\,\po{1-n}{h}\,\po{1+n+m}{h}\,\po{1-n-m}{h}
		}{
			\po{1}{h}\,\po{3/2-n}{h}\,\po{2+n}{h}\,\po{2-n-m}{h}\,\po{2+n+m}{h}
		}
\Biggr).
\end{multline*}
Using the standard hypergeometric notation
\begin{equation} \label{eq:hypergeom}
\HypF{r}{s}{a_1,\dots,a_r}{b_1,\dots,b_s}{z} =
	\sum_{k=0}^\infty
		\frac{
			\po{a_1}{k}\cdots\po{a_r}{k}
		}{
			k!\,\po{b_1}{k}\cdots\po{b_s}{k}
		}z^k,
\end{equation}
we may write the above expression as
\begin{multline} \label{eq:very-well-poised}
\frac{
n m \binom{2n}{n}\binom{2n-1}{n}\binom{2m}{m}
}{
\binom{4n+2m-1}{2n+m}
}
\Biggl(
	-\frac{1}{(n+m)^2}
	+\frac{4n+2}{(n+1)(2n-1)(n+m-1)(n+m+1)}\\
	\times
{}_7F_6\!\left[\begin{matrix} 2,2,1 + n + m, 1 - n - m, 1, {\frac 3 2} + n, 1 - n\\
1,2-n-m,2+n+m,2,{\frac 3 2}-n,2+n\end{matrix}; 1\right]
\Biggr).
\end{multline}
Now we observe that the sum in this expression is in fact a
very-well-poised hypergeometric series, and therefore we can transform it
into a more convenient form by applying 
Whipple's transformation (see \cite[(2.4.1.1)]{SlatAC})
\begin{multline}
\label{eq:Whipple}
\HypF{7}{6}{a,1+\frac{a}{2},b,c,d,e,-N}%
{\frac{a}{2},1+a-b,1+a-c,1+a-d,1+a-e,1+a+N}{1}\\
=\frac{\po{a+1}{N}\,\po{a-d-e+1}{N}}{\po{a-d+1}{N}\,\po{a-e+1}{N}}
\HypF{4}{3}{a-b-c+1,d,e,-N}{a-b+1,a-c+1,-a+d+e-N}{1}
\end{multline}
(where $N$ is a nonnegative integer)
with $a=2$, $b=1+n+m$, $c=1-n-m$, $d=1$, $e=3/2+n$, and $N=n-1$ to it.
Thus we obtain the expression
\begin{multline}\label{eq:4F3}
\frac{
n m \binom{2n}{n}\binom{2n-1}{n}\binom{2m}{m}
}{
\binom{4n+2m-1}{2n+m}
}
\Biggl(
	-\frac{1}{(n+m)^2}
	+\frac{2n+1}{(n+m-1)(n+m+1)}\\
	\times
   {} _{4} F _{3} \!\left [ \begin{matrix} { 1, 1, {\frac 3 2} + n, 1 - n}\\ 
   { 2 - n - m, 2 + n + m, {\frac 3 2}}\end{matrix} ; {\displaystyle 1}\right ] 
\Biggr).
\end{multline}
Now we substitute $m\sim a n$ and perform the limit
$n\rightarrow\infty$. The asymptotics of the binomials appearing in
front of the expression \eqref{eq:4F3} is easily determined by means of Stirling's
formula. For the $\HypFsimple{4}{3}$-series itself,
we may exchange limit and summation
by uniform convergence,
\begin{equation}\label{eq:2F1}
\lim_{n\rightarrow\infty}\HypF{4}{3}{1,1,\frac{3}{2}+n,1-n}%
{2-n-an,2+n+an,\frac{3}{2}}{1}=
\HypF{2}{1}{1,1}{\frac{3}{2}}{ {\frac {1} {(a+1)^2}} }.
\end{equation}
Combining all this, and making use of the identity 
(see \cite[p.~463, (133)]{Prudnikov})
\begin{equation}\label{eq:arcsin}
\HypF{2}{1}{1,1}{\frac{3}{2}}{z}=\frac{\arcsin\sqrt{z}}{\sqrt{z(1-z)}},
\end{equation}
we obtain exactly $\frac {2} {\pi}\arcsin (1/(a+1))$ as the
asymptotic form of \eqref{eq:4F3}, and, hence, as the asymptotic form
of the proportion of rhombus tilings in the statement of the
corollary, as desired.

\smallskip
The case of a hexagon with side lengths 
$2n+1,2m-1,2n+1,2n+1,2m-1,2n+1$ can be handled in (almost) the same
way because the sums in \eqref{eq:MEven} and \eqref{eq:MOdd} are
exactly the same.

\smallskip
For the next two cases, i.e., in order to estimate
\eqref{eq:NEven} and \eqref{eq:NOdd}, we proceed in a similar way.
Again we write the sum which appears in \eqref{eq:NEven} and
\eqref{eq:NOdd} as a $_7F_6$-series,
$${}_7F_6\!\left[\begin{matrix} 3,\frac {5} {2},1+n+m,2-n-m,2,\frac {3} {2}+n,
2-n\\\frac {3} {2},3-n-m,2+n+m,2,\frac {5} {2}-n,2+n\end{matrix}; 1\right],$$
apply Whipple's transformation \eqref{eq:Whipple}, and then let $n$ tend
to infinity. Here, the $_2F_1$-series which is obtained is a slightly
different one as before (compare \eqref{eq:2F1}),
$${}_2F_1\!\left[\begin{matrix} 1,2\\\frac {5} {2}\end{matrix}; \frac {1}
{(a+1)^2}\right].$$
In order to be able to use \eqref{eq:arcsin}, we use the relation
$${}_2F_1\!\left[\begin{matrix} 1,2\\\frac {5} {2}\end{matrix}; \frac {1}
{(a+1)^2}\right]=
-\frac {3(a+1)^2} {2} + 
  {\frac {{{3\left( a+1 \right) }^2}\,
 } 2}      {} _{2} F _{1} \!\left [ \begin{matrix} { 1, 1}\\
{{\frac 3 2}}\end{matrix} ; {\displaystyle \frac {1} {(a+1)^2}}\right
       ].
$$
The computation is then completed by straightforward use of Stirling's
formula, and subsequent simplification.

The remaining two cases, i.e., the estimations of \eqref{eq:NEven2}
and \eqref{eq:NOdd2}, can be dealt with in just the same manner.
The $_2F_1$-series which is obtained here is again slightly
different. It reads
$${}_2F_1\!\left[\begin{matrix} 1,2\\\frac {3} {2}\end{matrix}; \frac {1}
{(a+1)^2}\right]=
-{\frac {{{\left( a+1 \right) }^2}} 2} + 
  {\frac {{{\left( a+1 \right) }^2}
 } 2}
      {} _{2} F _{1} \!\left [ \begin{matrix} { 1, 1}\\ { {\frac 1
2}}\end{matrix} ;
       {\displaystyle \frac {1} {(a+1)^2}}\right ].
$$
To the $_2F_1$-series on the right-hand side 
we apply the formula (see \cite[p.~464, (138)]{Prudnikov})
$${} _{2} F _{1} \!\left [ \matrix { 1, 1}\\ { {\frac 1 2}}\endmatrix ;
   {\displaystyle z}\right ]  = 
  {\frac 1 {1 - z}} + {\frac {{\sqrt{z}}\,\arcsin {\sqrt{z}}} 
     {{{\left( 1 - z \right) }^{{\frac 3 2}}}}}.
$$
Again, the computation is then completed by straightforward use of Stirling's
formula, and subsequent simplification.
\hfil\qed


\section{From rhombus tilings to determinants}
\label{sec:rhombdet}

This section is entirely based on ideas by Helfgott and Gessel 
\cite{HeGeAA}.
These allow us to find a determinantal expression for the number 
of rhombus tilings
of a hexagon with side lengths $N,M,N,N,M,N$ that contain an {\em arbitrary\/}
fixed rhombus on the $(N+M)$-long symmetry axis (the vertical
symmetry axis in Figure~\ref{figure-thm6}). We shall state
two auxiliary results (Propositions~\ref{prop:HeGe1} and \ref{prop:HeGe}) 
without proof (the reader can find the details in
\cite{HeGeAA}), and then derive Helfgott and Gessel's determinant
(Proposition~\ref{thm:HeGe-det}). It is the
specialization $l=(N+M)/2$ of Proposition~\ref{thm:HeGe-det} (compare
the paragraph above \eqref{eq:gen-det}) which in the long run leads 
to a proof of our Theorems~\ref{thm:MEven} and \ref{thm:MOdd}, it is the
specialization $l=(N+M+1)/2$ of Proposition~\ref{thm:HeGe-det}
which in the long run leads 
to a proof of our Theorems~\ref{thm:NEven} and \ref{thm:NOdd}, and it is the
specialization $l=(N+M+2)/2$ of Proposition~\ref{thm:HeGe-det}
which in the long run leads 
to a proof of our Theorems~\ref{thm:NEven2} and \ref{thm:NOdd2}.
We do want to alert the reader that we use a different convention
in our figures of how to draw the hexagons than Helfgott and Gessel. 
To be precise, our figures turn into those in \cite{HeGeAA} by
a rotation by $90^\circ$.

\begin{figure}
\bigskip
\centertexdraw{
\drawdim truecm
  \linewd.02
  \rhombus \rhombus \rhombus \rhombus \ldreieck
  \move (-0.866025403784439 -.5)
  \rhombus \rhombus \rhombus \rhombus \rhombus \ldreieck
  \move (-1.732050807568877 -1)
  \rhombus \rhombus \rhombus \rhombus \rhombus \rhombus \ldreieck
  \move (-2.598076211353316 -1.5)
  \rhombus \rhombus \rhombus \rhombus \rhombus \rhombus \rhombus
  \move (-2.598076211353316 -1.5)
  \rdreieck \rhombus \rhombus \rhombus \rhombus \rhombus \rhombus
  \move (0.866025403784439 -3.5)
  \bsegment
    \rlvec(0.866025403784439 -.5) \rlvec(-0.866025403784439 -.5)
	   \rlvec(-0.866025403784439 .5)
    \lfill f:0.4
  \esegment
  \move (-2.598076211353316 -2.5)
  \rdreieck \rhombus \rhombus \rhombus \rhombus \rhombus
  \move (-2.598076211353316 -3.5)
  \rdreieck \rhombus \rhombus \rhombus \rhombus
  \htext (-1.4 -0.2){$N$} \htext (2.6 -0.2){$N$}
  \htext (-3 -3){$M$}     \htext (4.6 -3){$M$}
  \htext (-1.4 -6){$N$}   \htext (2.6 -6){$N$}
  \htext (0.80 -4.1){$2$}
}
\caption{}
\label{figure-thm6}
\end{figure}

The first observation is that for any rhombus tiling of the hexagon
with side lengths $N,M,N,N,M,N$, 
there are exactly $N$ rhombi of the tiling that are
cut in two by the $(N+M)$-long symmetry axis. Removing these rhombi, and cutting 
the hexagon in two along the symmetry axis, leaves two symmetric halves of trapezoidal
shape with $N$ ``dents''. The following statement
counts rhombus tilings of such ``dented trapezoids''. (This is Lemma~2 in
\cite{HeGeAA}.)

\begin{figure}
\bigskip
\centertexdraw{
  \drawdim truecm  \linewd.02
  \ldreieck
  \move (-0.866025403784439 -.5)
  \rhombus 
  \move (-1.732050807568877 -1)
  \rhombus \rhombus 
  \move (-2.598076211353316 -1.5)
  \rhombus \rhombus \rhombus
  \ldreieck
  \move (-2.598076211353316 -1.5)
  \rdreieck
  \rhombus \rhombus \rhombus
  \ldreieck
  \move (-2.598076211353316 -2.5)
  \rdreieck
  \rhombus \rhombus \rhombus 
  \move (-2.598076211353316 -3.5)
  \rdreieck
  \rhombus \rhombus \rhombus
\htext (-1.4 -0.2){$N$}
\htext (-3.5 -3){$M$}
\htext (-1.4 -6){$N$}
\htext (1 -0.1){$6$}
\htext (1 -1.1){$5$}
\htext (1 -2.1){$4$}
\htext (1 -3.1){$3$}
\htext (1 -4.1){$2$}
\htext (1 -5.1){$1$}
\htext (1 -6.1){$0$}
}
\caption{}
\label{figure-prop4}
\end{figure}

\begin{prop}\label{prop:HeGe1}
The number of rhombus tilings of
a semi-hexagon with side lengths $N,\break M,N$ (i.e., the ``half'' of a
hexagon with side lengths 
$N,M,N,N,M,N$) and $N$ ``dents'' at positions $0\leq r_0<\dots<r_{N-1}<N+M$ 
(see Figure~\ref{figure-prop4}, where $N=4$, $M=3$, 
and the ``dents'' are at positions $0$, $1$, $4$ and $5$) is
\begin{equation*}
\left(\prod_{i=1}^{N-1}\frac{1}{i!}\right)\prod_{0\leq i<j\le N-1}(r_i-r_j)
=\left(\prod_{i=1}^{N-1}\frac{1}{i!}\right)
	\det_{0\le i,j\le N-1}\left(p_i(r_j)\right) ,
\end{equation*}
where $p_i(x)$ is an {\em arbitrary monic\/} polynomial of degree $i$
in $x$.
(The ``standard'' case would be  $p_i(x)=x^i$, which corresponds to the
Vandermonde determinant.)
\end{prop}

{From} this proposition, Helfgott and Gessel 
deduce another enumeration result. (This is Proposition~4 in \cite{HeGeAA}.)

\begin{figure}
\bigskip
\centertexdraw{
  \drawdim truecm  \linewd.02
  \rhombus \rhombus \rhombus \rhombus \ldreieck
  \move (-0.866025403784439 -.5)
  \rhombus
  \move (0.866025403784439 -1.5)
  \rhombus \rhombus \rhombus \ldreieck
  \move (0.866025403784439 -0.5)
  \bsegment
    \rlvec(0.866025403784439 -.5) \rlvec(-0.866025403784439 -.5)
	   \rlvec(-0.866025403784439 .5)
    \lfill f:0.8
  \esegment
  \move (-1.732050807568877 -1)
  \rhombus \rhombus
  \move (0.866025403784439 -2.5)
  \rhombus \rhombus \rhombus \ldreieck
  \move (0.866025403784439 -1.5)
  \bsegment
    \rlvec(0.866025403784439 -.5) \rlvec(-0.866025403784439 -.5)
	   \rlvec(-0.866025403784439 .5)
    \lfill f:0.8
  \esegment
  \move (-2.598076211353316 -1.5)
  \rhombus \rhombus \rhombus
  \move (0.866025403784439 -3.5)
  \rhombus \rhombus \rhombus
  \move (0.866025403784439 -2.5)
  \bsegment
    \rlvec(0.866025403784439 -.5) \rlvec(-0.866025403784439 -.5)
	   \rlvec(-0.866025403784439 .5)
    \lfill f:0.8
  \esegment
  \move (-2.598076211353316 -1.5)
  \rdreieck \rhombus \rhombus \rhombus \rhombus \rhombus \rhombus
  \move (-2.598076211353316 -2.5)
  \rdreieck \rhombus \rhombus \rhombus 
  \move (0.866025403784439 -5.5)
  \rhombus
  \move (0.866025403784439 -4.5)
  \bsegment
    \rlvec(0.866025403784439 -.5) \rlvec(-0.866025403784439 -.5)
	   \rlvec(-0.866025403784439 .5)
    \lfill f:0.8
  \esegment
  \move (-2.598076211353316 -3.5)
  \rdreieck \rhombus \rhombus \rhombus
  \move (0.866025403784439 -5.5)
  \bsegment
    \rlvec(0.866025403784439 -.5) \rlvec(-0.866025403784439 -.5)
	   \rlvec(-0.866025403784439 .5)
    \lfill f:0.8
  \esegment
  \linewd.1
  \move (0.866025403784439 -4.5)
  \rlvec(0 1)
  \move (0.866025403784439 -0.5)
  \rlvec(0 1)
  \htext (-1.4 -0.2){$N$} \htext (2.6 -0.2){$N$}
  \htext (-3.5 -3){$M$}     \htext (4.6 -3){$M$}
  \htext (-1.4 -6){$N$}   \htext (2.6 -6){$N$}
  \htext (0.80 -1.1){$5$}
  \htext (0.80 -2.1){$4$}
  \htext (0.80 -3.1){$3$}
  \htext (0.80 -5.1){$1$}
  \htext (0.80 -6.1){$0$}
}
\caption{}
\label{figure-prop5}
\end{figure}

\begin{prop}
\label{prop:HeGe}
Let $L$ be a subset of $0,1,\dots,N+M-1$ of cardinality at least $N$.
Then the number of rhombus tilings of
a hexagon with side lengths $N,M,N,N,M,N$, in which the set of 
rhombi on the $(N+M)$-long symmetry axis is a subset of $L$
(given in terms of the numbers of the rhombi, where the rhombi on the
symmetry axis are numbered from bottom to
top as $0,1,\dots,N+M-1$; see Figure~\ref{figure-prop5}, where $N=4$,
$M=3$, $L=\{0,1,3,4,5\}$; the set $L$ consists of the shaded rhombi), is
\begin{equation}\label{eq:HeGe}
\left(\prod_{i=1}^{N-1}\frac{1}{i!^2}\right)
\det_{0\le i,j\le  N-1}\bigg(\sum_{s\in L}p_i(s)p_j(s)\bigg),
\end{equation}
where, again, $p_i(x)$ is an {\em arbitrary monic\/} polynomial of degree
$i$ in $x$.
\end{prop}

{From} this proposition 
we can derive the following determinantal expression for
the number of rhombus tilings
of a hexagon with side lengths $N,M,N,N,M,N$ which contain an arbitrary
fixed rhombus on the $(N+M)$-long symmetry axis (see \cite{HeGeAA}).

\begin{prop}
\label{thm:HeGe-det}
The number of rhombus tilings of
a hexagon with side lengths $N,M,N,\break N,M,N$ which contain the $l$-th 
rhombus, $0\leq l\leq N+M-1$, on the $(N+M)$-long symmetry axis
(see Figure~\ref{figure-thm6}; the rhombus which is contained in every
tiling is shaded, i.e.,
$l=2$), is
\begin{equation}
\label{eq:HeGe-det}
\left(\prod_{i=1}^{N-1}\frac{1}{i!^2}\right)
\det_{1\le i,j\le  N-1}\bigg(\sum_{s=-l}^{N+M-l-1}s^{i+j}\bigg).
\end{equation}
\end{prop}

\begin{proof}
Let us first count the complementary set, i.e., the rhombus tilings
which do {\em not\/} contain rhombus $l$.
Obviously, we get this number from Proposition~\ref{prop:HeGe} with
$L=\{0,1,\dots,N+M-1\}\setminus\{l\}$. Then, if in addition we choose
$p_i(x)=(x-l)^i$, by formula \eqref{eq:HeGe} this number is
\begin{equation}
\label{eq:aux}
\left(\prod_{i=1}^{N-1}\frac{1}{i!^2}\right)
	\det_{0\le i,j\le  N-1}\Bigg(\underset{s\ne l}{\sum_{s=0}^{N+M-1}}(s-l)^{i+j}\Bigg)
=
\left(\prod_{i=1}^{N-1}\frac{1}{i!^2}\right)
	\det_{0\le i,j\le  N-1}\Bigg(\underset{s\ne
0}{\sum_{s=-l}^{N+M-l-1}}s^{i+j}\Bigg).
\end{equation}
We have to subtract this number from the total number of possible rhombus
tilings. This can be again expressed by making use of 
Proposition~\ref{prop:HeGe}, this time with $L=\{0,1,\dots,N+M-1\}$.
So, choosing $p_i(x)=(x-l)^i$ in \eqref{eq:HeGe} again, we obtain
\begin{equation}
\label{eq:auxtotal}
\left(\prod_{i=1}^{N-1}\frac{1}{i!^2}\right)
	\det_{0\le i,j\le N-1}\bigg({\sum_{s=0}^{N+M-1}}(s-l)^{i+j}\bigg)
=
\left(\prod_{i=1}^{N-1}\frac{1}{i!^2}\right)
	\det_{0\le i,j\le N-1}\bigg({\sum_{s=-l}^{N+M-l-1}}s^{i+j}\bigg)
\end{equation}
for the total number of rhombus tilings.

It should be observed that the determinants in \eqref{eq:aux} and
\eqref{eq:auxtotal} are almost the same. The only difference
is the $(0,0)$-entry, which is by 1 less in \eqref{eq:aux} than in
\eqref{eq:auxtotal}. Therefore, if we expand the determinant in
\eqref{eq:aux} with respect to the first row, we may rewrite \eqref{eq:aux} 
as
\begin{equation}
\label{eq:aux2}
\left(\prod_{i=1}^{N-1}\frac{1}{i!^2}\right)
\left(\det_{0\le i,j\le  N-1}\bigg(\sum_{s=-l}^{N+M-l-1}s^{i+j}\bigg) -
\det_{1\le i,j\le N-1}\bigg(\sum_{s=-l}^{N+M-l-1}s^{i+j}\bigg)\right).
\end{equation}
The number of rhombus tilings that we are interested in is the
difference of \eqref{eq:auxtotal} and \eqref{eq:aux2}, which is
exactly \eqref{eq:HeGe-det}.
\end{proof}


\section{Determinant evaluations}
\label{sec:proofs}

\begin{lem} \label{lem:Det1}
Let $n$ and $N$ be positive integers. Then the determinant $D(n,n,N)$,
as defined in {\em\eqref{eq:gen-det}}, is equal to
\begin{multline} \label{eq:Det1}
n^N\prod_{i=1}^{\floor{N/2}}\big((n^2-i^2)^{N-2i+1}
(n^2-(i-1/2)^2)^{N-2i+1}\big)\\
\times
\frac{
	2^{N^2}N!\,(N+1)!\,\po{n-\ceil{N/2}}{2\ceil{N/2}+1}\prod_{i=1}^{N} i!^4	
}{
	n\floor{N/2}!^2\ceil{N/2}!^2 \prod_{i=1}^{2N+1}i!
}\\
\times\Biggl(	
	(-1)^N +
	\frac{
		4(N+2)n^2\floor{N/2}\ceil{N/2}
	}{
		(N-1)N(n^2-1)(\floor{N/2}+1)(\ceil{N/2}+1)
	}\\
\times	\sum_{h=0}^{\ceil{N/2}-1}
		\frac{
		\po{2}{h}\,\po{1-\ceil{N/2}}{h}\,\po{3/2+\ceil{N/2}}{h}\,
		\po{1-n}{h}\,\po{1+n}{h}
	}{
		\po{1}{h}\,\po{2+\ceil{N/2}}{h}\,\po{3/2-\ceil{N/2}}{h}\,
	\po{2+n}{h}\,\po{2-n}{h}
	}
\Biggr).
\end{multline}
\end{lem}

\begin{proof}
We proceed in several steps. An outline is as
follows.
In the first step we make the obvious observation that
$D(n,n,N)$ is actually a polynomial in $n$, of degree at most
$N(N+2)$. 
Next, we show that $D(n,n,N)$, as polynomial in $n$, has a lot of linear factors.
More precisely, in the second step, we show that $n^N$ is a factor of $D(n,n,N)$.
Then, in the third step, we show that $\prod _{i=1}
^{\floor{N/2}}(n^2-i^2)^{N-2i+1}$ is a factor of $D(n,n,N)$. Moreover,
in the fourth step, we show that $\prod _{i=1}
^{\floor{N/2}}(n^2-(i-1/2)^2)^{N-2i+1}$ is a factor of $D(n,n,N)$. 
 From a
combination of these four steps we are forced to conclude that
\begin{equation} \label{eq:polydef}
D(n,n,N)=n^N\bigg(\prod_{i=1}^{\floor{N/2}}(n^2-i^2)^{N-2i+1}
(n^2-(i-1/2)^2)^{N-2i+1}\bigg)	\cdot P(n,N),
\end{equation}
where $P(n,N)$ is a polynomial in $n$ of degree at most $2\ceil{N/2}$.
Finally, in the fifth step, we evaluate $P(n,N)$ at
$n=-\ceil{N/2},-\ceil{N/2}+1,\dots,\ceil{N/2}-1,\ceil{N/2}$.
Namely, we show that 
\begin{equation}
\label{eq:poly-at-0-ganz}
P(0,N)=\frac {(-1)^{N(N-1)/2}
\,3^{-N}}
{\prod_{i=1}^{\floor{N/2}}i^{2(N-2i+1)}\,(i-1/2)^{2(N-2i+1)}}
\prod_{i = 1}^{N-1}\left(\frac{i\,{{\left( i+1 \right) }^4}\,
             \left( i+2 \right) }{\left( 2 i+1 \right) \,
             {{\left(  2 i+2 \right) }^2}\,\left(  2 i+3 \right) }
         \right)^{N-i},
\end{equation}
and, for $1\le e\le\ceil{N/2}$, that
\begin{multline}
\label{eq:poly-at-eplus-ganz}
P(\pm e,N) =\dfrac {(-1)^{(N-2e+1)(N-2e)/2}\,2^{N-2e+1}
\prod_{i=1}^{2e-1}i!^2}
{\displaystyle e^N\underset{i\ne e}{\prod _{i=1} ^{\floor{N/2}}}(e-i)^{N-2i+1}
\prod_{i=1}^{\floor{N/2}}(e+i)^{N-2i+1}
(e^2-(i-1/2)^2)^{N-2i+1}}
\\
\times
\left(\frac{\left( 2 e \right) !^4}
  {\left( 4 e+1 \right) !}\right)^{N-2e+1}
   \prod_{i = 1}^{N-2e}\left(\frac{i\,\left( 2 e + i \right) ^4\,
            \left( 4 e + i \right) }{\left(  4 e + 2 i -1\right) \,
            {{\left( 4 e + 2 i \right) }^2}\,
            \left(  4 e + 2 i +1\right) }\right)^{N-2e+1-i}.
\end{multline}
Clearly, this determines a polynomial of maximal
degree $2\ceil{N/2}$ uniquely. In fact, an explicit expression for $P(n,N)$ can
immediately be written down using Lagrange interpolation.
As it turns out, the resulting expression for $P(n,N)$ is exactly
the expression covering the last three lines of \eqref{eq:Det1}.
In view of \eqref{eq:polydef},
this would finish the proof of the Lemma.

\newstep{Step 1. $D(n,n,N)$ is a polynomial in $n$ of degree at most $N(N+2)$}
It is standard that sums of powers, such as the entries of $D(n,n,N)$,
can be expressed using Bernoulli numbers. More precisely, we have
(cf.~\cite[p. 269ff]{Knuth}),
\begin{equation}
\label{eq:BernoulliP}
\sum _{s=-n} ^{n-1}s^{m}=\sum_{\ell =0}^m\frac {1}
{\ell +1}\binom{m}{\ell }B_{m-\ell }\left( n^{\ell +1}-(-n)^{\ell +1}\right),
\end{equation}
where $B_\ell $ denotes the $\ell $-th Bernoulli number.
Hence, the $(i,j)$-entry of $D(n,n,N)$ is a polynomial in $n$ of degree
$i+j+1$. Thus, by expanding the determinant $D(n,n,N)$ according to the
definition of a determinant and determining the degree of each term, 
it follows that $D(n,n,N)$ is a polynomial
in $n$ of degree at most $\sum _{1\le i,j\le N} ^{}(i+j+1)=2\binom
{N+1}2 +N=N(N+2)$.

{From} \eqref{eq:BernoulliP} we may read off another property of
$D(n,n,N)$, which we record here for later use,
\begin{equation}
D(-n,-n,N) 	= (-1)^ND(n,n,N).\label{eq:symD}
\end{equation}

\newstep{Step 2. $n^N$ is a factor of $D(n,n,N)$}
{From} the definition \eqref{eq:gen-det} of the entries of $D(n,n,N)$ and
\eqref{eq:BernoulliP} it is immediate that $n$ divides each entry of
$D(n,n,N)$. Hence, $n^N$ divides $D(n,n,N)$.

\newstep{Step 3. $\prod _{i=1}
^{\floor{N/2}}(n^2-i^2)^{N-2i+1}$ is a factor of $D(n,n,N)$}
In view of \eqref{eq:symD}, it suffices to prove that $(n-e)^{N-2e+1}$
divides $D(n,n,N)$ for $e=1,2,\dots,\floor{N/2}$. 
In order to do so, we claim that for each such $e$ there are $N-2e+1$
linear combinations of the columns, which are themselves linearly
independent, that vanish for $n=e$. 
Define the coefficients $a_j$ by 
\begin{equation}
\label{eq:lincomb-coefficients}
\underset{z\ne0}{\prod_{z=-e}^{e-1}}(x-z)=
\po{x-e+1}{e-1}\,\po{x+1}{e}=
\sum_{j=0}^{2e-1} a_jx^j.
\end{equation}
Then we claim that for $k=2e,2e+1,\dots,N$ we have
\begin{equation}
\label{eq:lincomb-definition}
\sum_{j=k-2e+1}^{k}a_{j-k+2e-1}\cdot\left(\text{column $j$ 
in $D(e,e,N)$}\right)=0.
\end{equation}
To verify \eqref{eq:lincomb-definition}, for $i=1,2,\dots,N$ we have
to prove that 
\begin{equation}
\label{eq:lincomb-compute}
\sum_{j=k-2e+1}^{k}a_{j-k+2e-1}\sum_{s=-e}^{e-1}s^{i+j}=0.
\end{equation}
To do so, we interchange sums in \eqref{eq:lincomb-compute},
and then use \eqref{eq:lincomb-coefficients}, to obtain
\begin{equation*}
\sum_{s=-e}^{e-1}
	\left(s^{i+k-2e+1}\sum_{j=0}^{2e-1}a_{j}s^{j}\right)
=\sum_{s=-e}^{e-1}
\bigg(s^{i+k-2e+1}\underset{z\ne0}{\prod_{z=-e}^{e-1}}(s-z)\bigg)
\end{equation*}
for the left-hand side in \eqref{eq:lincomb-compute},
which clearly vanishes for any integer $e>0$.

\newstep{Step 4. $\prod _{i=1}
^{\floor{N/2}}(n^2-(i-1/2)^2)^{N-2i+1}$ is a factor of $D(n,n,N)$}
In view of \eqref{eq:symD}, it suffices to prove that $(n-(e-1/2))^{N-2e+1}$
divides $D(n,n,N)$ for $e=1,2,\dots,\floor{N/2}$. In order to do so,
we claim, in a manner analogous
to Step~3, that for each such $e$ there are $N-2e+1$
linear combinations of the columns, which are themselves linearly
independent, that vanish for $n=e-1/2$. 
Define the coefficients $\tilde a_j$ by 
\begin{equation}
\label{eq:lincomb-coefficients-halb}
\prod_{z=-e+1}^{e-1}(x-z+1/2)=
\po{x-e+3/2}{2e-1}=
\sum_{j=0}^{2e-1} \tilde a_jx^j.
\end{equation}
Then we claim that for $k=2e,2e+1,\dots,N$ we have
\begin{equation}
\label{eq:lincomb-definition-halb}
\sum_{j=k-2e+1}^{k}\tilde a_{j-k+2e-1}\cdot\left(\text{column $j$ 
in $D(e-1/2,e-1/2,N)$}\right)=0.
\end{equation}
To verify \eqref{eq:lincomb-definition-halb}, for $i=1,2,\dots,N$ we have
to prove that 
\begin{equation*}
\sum_{j=k-2e+1}^{k}\tilde a_{j-k+2e-1}
\sum_{\ell =0}^{i+j}\frac {1}
{\ell +1}\binom{i+j}{\ell }B_{i+j-\ell }\left( (e-1/2)^{\ell
+1}-(-e+1/2)^{\ell +1}\right)=0,
\end{equation*}
where we used the right-hand side of formula \eqref{eq:BernoulliP} to
express the entries of $D(n,n,N)$. By a variation of
\eqref{eq:BernoulliP}, the last equation can be rewritten as
\begin{equation}
\label{eq:lincomb-compute-halb}
\sum_{j=k-2e+1}^{k}\tilde a_{j-k+2e-1}
\sum _{s=-e+1} ^{e-1}(s-1/2)^{i+j}=0.
\end{equation}
Again, to show that \eqref{eq:lincomb-compute-halb} is true,
we interchange sums,
and then use \eqref{eq:lincomb-coefficients-halb}, to obtain
\begin{equation*}
\sum_{s=-e+1}^{e-1}
	\left((s-1/2)^{i+k-2e+1}\sum_{j=0}^{2e-1}\tilde a_{j}(s-1/2)^{j}\right)
=\sum_{s=-e+1}^{e-1}
	\left((s-1/2)^{i+k-2e+1}{\prod_{z=-e+1}^{e-1}}(s-z)\right)
\end{equation*}
for the left-hand side in \eqref{eq:lincomb-compute-halb},
which clearly vanishes for any integer $e>0$.

\newstep{Step 5. Evaluation of the polynomial $P(n,N)$ at
$n=-\ceil{N/2},-\ceil{N/2}+1,\dots,\ceil{N/2}$}
We start by observing that the symmetry relation \eqref{eq:symD} for
$D(n,n,N)$ is ``inherited" by $P(n,N)$. To be precise, we have
\begin{equation}
P(-n,N) 	= P(n,N).\label{eq:symP}
\end{equation}
In view of \eqref{eq:symP} it suffices to determine the evaluations of
$P(n,N)$ at $n=0,1,\dots,\ceil{N/2}$. 

What we would like to do is, for any $e$ with $0\le 
e\le \ceil{N/2}$,
to set $n=e$ in
\eqref{eq:polydef}, compute $D(e,e,N)$,
and then express $P(e,N)$ as
the ratio of $D(e,e,N)$ and the right-hand side
product evaluated 
at $n=e$. Unfortunately, this is typically a ratio $0/0$ and, hence,
undetermined. So, we have to first divide both sides of \eqref{eq:polydef}
by the appropriate power of $(n-e)$, and only then set $n=e$.

This program is easily carried out for $e=0$.
As we observed in Step~2, each entry of $D(n,n,N)$ is divisible by $n$.
Hence, division of both sides of \eqref{eq:polydef} by $n^N$ and then
specializing to $0$, transforms
\eqref{eq:polydef} into the equation
\begin{equation*}
P(0,N)=\frac{
	\det_{1\le i,j\le N}\left(2B_{i+j}\right)
}{\prod_{i=1}^{\floor{N/2}}i^{2(N-2i+1)}\,(i-1/2)^{2(N-2i+1)}}.
\end{equation*}
The Hankel determinant of Bernoulli numbers which appears in the
right-hand side expression can be evaluated by 
Theorem~\ref{thm:HankelBernoulli} with $n$ replaced by $N$ and $a=b=c=d=1$.
Thus, we obtain \eqref{eq:poly-at-0-ganz}.
\smallskip

The case $1\leq e \leq\ceil{N/2}$ requires more work. First, we have
to ``preprocess'' the determinant $D(n,n,N)$. Define the coefficients
$a_j$ as before in \eqref{eq:lincomb-coefficients},
\begin{equation*}
\underset{z\ne0}{\prod_{z=-e}^{e-1}}(x-z)=
\po{x-e+1}{e-1}\,\po{x+1}{e}=
\sum_{j=0}^{2e-1} a_jx^j.
\end{equation*}
Then, for $k=N,N-1,\dots,2e$,
in this order, add 
\begin{equation*}
\sum_{j=k-2e+1}^{k-1}a_{j-k+2e-1}\cdot\left(\text{column $j$ 
in $D(n,n,N)$}\right)
\end{equation*}
to column $k$. Thus, by \eqref{eq:lincomb-definition}, each entry in
column $k$, $k=2e,2e+1,\dots,N$, will be divisible by $(n-e)$ after
performing these operations. Next, apply the analogous row operations.
I.e., for $k=N,N-1,\dots,2e$,
in this order, add 
\begin{equation*}
\sum_{j=k-2e+1}^{k-1}a_{j-k+2e-1}\cdot\left(\text{row $j$ 
in $D(n,n,N)$}\right)
\end{equation*}
to row $k$.

Now we divide $(n-e)^{N-2e+1}$ on both sides of \eqref{eq:polydef},
and only then set $n=e$. As a result, from equation \eqref{eq:polydef} 
we obtain
\begin{align} \notag
P(e,N)&=\dfrac {\det_{1\le i,j,\le N}\left(\begin{matrix}
X&*\\0&Y\end{matrix}\right)} 
{\displaystyle e^N\underset{i\ne e}{\prod _{i=1} ^{\floor{N/2}}}(e-i)^{N-2i+1}
\prod_{i=1}^{\floor{N/2}}(e+i)^{N-2i+1}
(e^2-(i-1/2)^2)^{N-2i+1}} \\
&=\dfrac {\det X\,\det Y}
{\displaystyle e^N\underset{i\ne e}{\prod _{i=1} ^{\floor{N/2}}}(e-i)^{N-2i+1}
\prod_{i=1}^{\floor{N/2}}(e+i)^{N-2i+1}
(e^2-(i-1/2)^2)^{N-2i+1}},
\label{eq:matXY}
\end{align}
where $X$ is the $(2e-1)\times (2e-1)$ matrix
$X=\big(\sum_{s=-e}^{e-1} s^{i+j}\big)_{i,j=1,\dots, 2e-1}$, and
where $Y$ is the $(N-2e+1)\times (N-2e+1)$ matrix
$(Y_{ij})_{i,j=1,\dots, N-2e+1}$ whose entries are given by
$$Y_{ij}=\bigg(\frac {1}{n-e}
\sum_{s=-n}^{n-1}s^{i+j}\left(\po{s-e+1}{e-1}\,\po{s+1}{e}\right)^2\bigg)
\bigg\vert_{n=e}\ .$$
Alternatively, 
\begin{equation} \label{eq:Yij} 
Y_{ij}=\frac {d} {d n}\bigg(
\sum_{s=-n}^{n-1}s^{i+j}\left(\po{s-e+1}{e-1}\,\po{s+1}{e}\right)^2\bigg)
\bigg\vert_{n=e}\ .
\end{equation}
The expression on the right-hand side is a certain linear combination
of expressions of the form
$$H_m:=\frac {d} {d n}\bigg(
\sum_{s=-n}^{n-1}s^m\bigg)
\bigg\vert_{n=e}\ .
$$
In order to compute $H_m$, we use \eqref{eq:BernoulliP} (and
variations thereof) to rewrite it as
\begin{align}
\notag
H_m&=\sum _{\ell =0} ^{m}\binom m\ell  B_{m-\ell }(e^\ell
+(-e)^\ell )\\
&=2B_m+ \sum _{\ell =1} ^{m}m\binom {m-1}{\ell -1}\frac {1} {\ell }
B_{m-\ell }(e^\ell +(-e)^\ell )\notag\\
&=2B_m+ \sum _{s=0} ^{e-1}ms^{m-1}-\sum _{s=-e} ^{-1}ms^{m-1}\notag\\
&=2B_m+ \sum _{s=0} ^{e-1}\left(\frac {d} {ds}s^m\right)
-\sum _{s=-e} ^{-1}\left(\frac {d} {ds}s^m\right).\label{eq:H_m}
\end{align}
Thus, using the symbolic notation $B^k\equiv B_k$ for the
Bernoulli numbers, a combination of \eqref{eq:Yij} and \eqref{eq:H_m}
leads to
\begin{align*} 
Y_{ij}&=2B^{i+j}\left(\po{B-e+1}{e-1}\,\po{B+1}{e}\right)^2+
\sum _{s=0} ^{e-1}\frac {d} {ds}\left(s^{i+j}
\left(\po{s-e+1}{e-1}\,\po{s+1}{e}\right)^2\right)\\
&\quad \quad 
-\sum _{s=-e} ^{-1}\frac {d} {ds}\left(s^{i+j}
\left(\po{s-e+1}{e-1}\,\po{s+1}{e}\right)^2\right)\\
&=2B^{i+j}\left(\po{B-e+1}{e-1}\,\po{B+1}{e}\right)^2,
\end{align*}
as each summand in either sum on the right-hand side vanishes.

Substituting all this into \eqref{eq:matXY}, we arrive at
\begin{equation}
\label{eq:poly-at-eplus}
P(e,N) = \dfrac{D(e,e,2e-1)
	\det_{1\le i,j\le N-2e+1}\left(
		2B^{i+j}\left(\po{B-e+1}{e-1}\,\po{B+1}{e}\right)^2
	\right)}
{\displaystyle e^N\underset{i\ne e}{\prod _{i=1} ^{\floor{N/2}}}(e-i)^{N-2i+1}
\prod_{i=1}^{\floor{N/2}}(e+i)^{N-2i+1}
(e^2-(i-1/2)^2)^{N-2i+1}},
\end{equation}
which is valid for $1\leq e\leq \ceil{N/2}$.
The determinant $D(e,e,2e-1)$ is evaluated separately in
Lemma~\ref{lem:detklein}. 
The Hankel determinant of (linear combinations of) Bernoulli numbers
which appears in the 
right-hand side expression can be evaluated by 
Theorem~\ref{thm:HankelBernoulli} with $n$ replaced by $N-2e+1$, 
$a=b=e+1$ and $c=d=e$.
Thus, we obtain \eqref{eq:poly-at-eplus-ganz}.

\smallskip
This completes the proof of the Lemma.
\end{proof}

\begin{lem} \label{lem:Det2}
Let $n$ and $m$ be positive integers. Then the determinant
$D(n,n-1,2m-1)$,
as defined in \eqref{eq:gen-det}, is equal to
\begin{multline}\label{eq:Det2odd}
\frac{{2^{- 4 m + 5}}}{\po{
       {n- m + \frac{1}{2}}}{2 m - 1}\,
      \left( \prod_{i = 1}^{4 m - 1}i! \right) }\,\left(2 m + 1\right) \,
      {{\left(2 n - 1\right) }^2}\,{\binom{2 m}{m}}\,
      {\binom{2 m - 4}{m - 2}}
		\\\times
		\po{n-m}{m - 2}\,
      \po{n + 2}{m - 2}
      \left( \prod_{i = 1}^{2 m - 1}{{i!}^4} \right) \,
      \left( \prod_{i = 1}^{2 m - 2}
         \po{2 n -i-2}{2i + 3}\right)
\\\times
\Biggl(
{\frac {m \left( 2 m-3 \right)  \left( 2 m-1 \right)  
      \left( n-2 \right)  \left( n+1 \right)  
      \left( n^2-n + 2 m+1 \right) } 
    {\left( m-1 \right)  \left( 2 m+1 \right) }}\kern5cm\\
+
  {\frac {6(n-1)n} {(m+1)}
		\sum_{h=0}^{m-2}
			\frac{
				\po{3}{h}\,\po{\frac{5}{2}}{h}\,
                \po{2 - m}{h}\,\po{\frac{3}{2} + m}{h}\,
					\po{2 - n}{h}\,\po{1 + n}{h}
			}{
				\po{1}{h}\,\po{\frac{3}{2}}{h}\,
                \po{2 + m}{h}\,\po{\frac{5}{2} - m}{h}\,
					\po{2 + n}{h}\,\po{3 - n}{h}
			}
	}
\Biggr).
\end{multline}
Likewise, the determinant $D(n,n-1,2m)$ is equal to
\begin{multline}\label{eq:Det2even}
\frac{4}{\prod_{i = 1}^{4 m + 1}
      i!}\,\left(2 m + 3\right) \,{{\left(2 n - 1\right) }^2}\,
     {\binom{2 m - 2}{m - 1}}\,{\binom{2 m + 2}{m + 1}}
	  \\\times
     \po{n-m}{m - 2}\,
     \po{n + 2}{m - 2}\,
     \po{n- m - 1}{2 m + 2}
     \left( \prod_{i = 1}^{2 m}{{i!}^4} \right) \,
     \left( \prod_{i = 1}^{2 m - 2}
        \po{2 n -i-2}{2i + 3} \right)
\\\times
\Biggl(
{\frac {\left( 2 m-1 \right)  \left( m+1 \right)  \left( 2 m+1 \right)  
     \left(  n-2 \right)  \left( n+1 \right)  
     \left( n^2-n-2 m  \right) } {m \left( 2 m+3 \right)
}}\kern5cm\\
+
  {\frac {6(n-1)n} {(m+2)}
		\sum_{h=0}^{m-1}
			\frac{
				\po{3}{h}\,\po{\frac{5}{2}}{h}\,
                \po{1 - m}{h}\,\po{\frac{5}{2} + m}{h}\,
					\po{2 - n}{h}\,\po{1 + n}{h}
			}{
				\po{1}{h}\,\po{\frac{3}{2}}{h}\,
                \po{3 + m}{h}\,\po{\frac{3}{2} - m}{h}\,
					\po{2 + n}{h}\,\po{3 - n}{h}
			}
	}
\Biggr)
\end{multline}
\end{lem}

\begin{proof}
Basically, the proof proceeds in the same way as above. By 
considerations which parallel Steps~1--4 of the previous proof, 
we deduce that
\begin{multline} \label{eq:polydef2}
D(n,n-1,N)=
	\left(n-\frac{1}{2}\right)^N
	\left((n-1)n\right)^{N-1}\\
\times
	\Bigg(\prod_{i=2}^{\ceil{(N+1)/2}}\left(
		\left(n-i+\frac{1}{2}\right)
        \left(n+i-\frac{3}{2}\right)(n-i)(n+i-1)
	\right)^{N-2i+2}\Bigg)	\cdot Q(n,N),
\end{multline}
where $Q(n,N)$ is a polynomial in $n$ of degree at most $2\ceil{(N+1)/2}$.

Also, Step~5 of the previous proof has a parallel here. Eventually, 
this yields evaluations of $Q(n,N)$ at
$n=-\ceil{(N+1)/2}+1,-\ceil{(N+1)/2}+2,
\dots,\ceil{(N+1)/2}-1,\ceil{(N+1)/2}$. In particular, the symmetry
relation which plays the role of \eqref{eq:symP} in this new context is
\begin{equation}
Q(-n+1,N) 	= Q(n,N),\label{eq:symQ}
\end{equation}
while the evaluation of determinant $D(e,e-1,2e-2)$ which is 
needed here (and replaces the evaluation of $D(e,e,2e-1)$ 
in this new context; compare \eqref{eq:poly-at-eplus}) 
is evaluated separately in Lemma~\ref{lem:detklein2}.

Unfortunately, this is not good enough. The polynomial $Q(n,N)$
is a polynomial of maximal degree $2\ceil{(N+1)/2}$, but by now we have
found only $2\ceil{(N+1)/2}$ explicit special evaluations of
$Q(n,N)$. Hence, we need one more information about $Q(n,N)$.
We get this missing piece of information by computing the 
leading coefficient of $Q(n,N)$. This is easily done. By the
definition of $D(n,n-1,N)$, given by \eqref{eq:gen-det}, by the
analogue of \eqref{eq:BernoulliP},
\begin{equation}
\label{eq:BernoulliP2}
\sum _{s=-n} ^{n-2}s^{m}=\sum_{\ell =0}^m\frac {1}
{\ell +1}\binom{m}{\ell }B_{m-\ell }\left( (n-1)^{\ell +1}-(-n)^{\ell +1}\right),
\end{equation}
and by \eqref{eq:polydef2}, 
the leading coefficient is given by
\begin{equation*}
\det_{1\leq i,j\leq N}\left(
	\begin{cases}
		\frac{2}{i+j+1} & \text{if $i+j$ is even,} \\
		0 & \text{otherwise.}
	\end{cases}
\right)  .
\end{equation*}
It is easy to see, by sorting rows and columns with odd indices
to the beginning, that this determinant equals the product of the
minor consisting of the odd rows and columns times the minor 
consisting of the even rows and columns, explicitly
\begin{equation*}
\det_{1\leq i,j\leq \ceil{N/2}}\left(
		\frac{2}{2i+2j-1} 
\right)  
\det_{1\leq i,j\leq \floor{N/2}}\left(
		\frac{2}{2i+2j+1} 
\right)  .
\end{equation*}
Both of these determinants are easily evaluated by means of the
Cauchy determinant evaluation (see \cite[vol.~III, p.~311]{MuirAB})
\begin{equation*}
\det_{1\leq i,j\leq K}\left(
		\frac{1}{x_i+y_j} 
\right)  =
	\frac{
		\prod_{1\leq i<j\leq K}(x_i-x_j)(y_i-y_j)
	}{
		\prod_{1\leq i, j\leq K}(x_i+y_j)
	}.
\end{equation*}

For the simplifications of the resulting expression for $Q(n,N)$ 
it turns out to be convenient to separate the cases $N=2m$ and $N=2m-1$.

\smallskip
This concludes the proof of the Lemma.
\end{proof}

\begin{lem}
\label{lem:Det3}
Let $n$ and $m$ be positive integers. 
Then the determinant $D(n,n-2,2m+1)$,
as defined in {\em\eqref{eq:gen-det}}, is equal to
\begin{multline}
\label{eq:Det3odd}
\frac{
	(2m-2)!\,(2m+2)!\,
	\po{n-m-1}{2m+1}\,\po{n-m-2}{2m+3}
}{
	(m-1)!\,m!\,(m+1)!\,(m+4)!\,(4m+3)!\,n(n-2)
}
\\
\times
\prod_{i=0}^{2m}
\frac{
	(i+1)!^3\,\po{2n -i- 2}{2i+1}
}{
	(2m+i+2)!
}
\Biggl( {\frac {2   \left( m+3 \right)  \left( m+4 \right)  
     \left( 2 m-1 \right) \left( 2 m +1\right)  x(m,n) }
   {\left( n-2 \right)  {{\left( n-1 \right) }^2} n}}\\
	-
	\frac{
		24(2m + 3)(2m + 5)(n - m - 3)(n + m + 1)
	}{
		(n - 3)(n + 1)
	}\\ 
\cdot
		\sum_{h = 0}^{m}
			\frac{
				\po{4}{h}\,
				\po{- m}{h}\,
				\po{{\frac{7}{2}} + m}{h}\,
				\po{3 - n}{h}\,
				\po{1 + n}{h}
			}{
				\po{1}{h}\,
				\po{5 + m}{h}\,
				\po{{\frac{3}{2}} - m}{h}\,
				\po{2 + n}{h}\,
				\po{4 - n}{h}
			}
\Biggr),
\end{multline}
where 

\vskip-10pt
{\small
\begin{multline*} 
x(m,n) = -72 m - 204 {m^2} - 212 {m^3} - 
       96 {m^4} - 16 {m^5} - 30 n - 86 m n - 68 {m^2} n - 
       16 {m^3} n \\
- 29 {n^2} + 3 m {n^2} + 26 {m^2} {n^2} + 
       8 {m^3} {n^2} + 44 {n^3} + 40 m {n^3} + 8 {m^2} {n^3} - 
       11 {n^4} - 10 m {n^4} - 2 {m^2} {n^4} .
\end{multline*}
}
\vskip-10pt

\noindent
Likewise, the determinant $D(n,n-2,2m)$ is equal to
\begin{multline}
\label{eq:Det3even}
\frac{
	(2m-2)!\,(2m+2)!\,
	\po{n-m-1}{2m+1}^2
}{
	(2m-3)\,(m-1)!\,m!\,(m+1)!\,(m+3)!\,(4m+1)!\,n(n-2)
}
\\
\times
\prod_{i=0}^{2m - 1}
\frac{
	(i+1)!^3\,\po{2n -i- 2}{2i+1}
}{
	(2m+i+1)!
}
\Biggl( {\frac {  \left( m+2 \right)  \left( m+3 \right)  
     \left( 2 m-3 \right) \left( 2 m-1 \right)  y(m,n) }
   {\left( n-2 \right)  {{\left( n-1 \right) }^2} n}}\\
	-
	\frac{
		24m(2m + 3)(n - m - 2)(n + m)
	}{
		(n - 3)(n + 1)
	}
		\sum_{h = 0}^{m - 1}
			\frac{
				\po{4}{h}\,
				\po{1 - m}{h}\,
				\po{{\frac{5}{2}} + m}{h}\,
				\po{3 - n}{h}\,
				\po{1 + n}{h}
			}{
				\po{1}{h}\,
				\po{4 + m}{h}\,
				\po{{\frac{5}{2}} - m}{h}\,
				\po{2 + n}{h}\,
				\po{4 - n}{h}
			}
\Biggr),
\end{multline}
where 

\vskip-10pt
{\small
\begin{multline*} 
y(m,n)=     12 {m^2} + 44 {m^3} + 48 {m^4} + 16 {m^5} + 16 m n + 
       28 {m^2} n + 16 {m^3} n + 4 {n^2} - 20 m {n^2}\\
 - 
       22 {m^2} {n^2} - 8 {m^3} {n^2} - 4 {n^3} + 12 m {n^3} + 
       8 {m^2} {n^3} + {n^4} - 3 m {n^4} - 2 {m^2} {n^4} .
\end{multline*}
}
\end{lem}

\begin{proof}
Basically, the proof proceeds in the same was as in the preceding lemmas.
By considerations which parallel Steps 1--4 of the proof of Lemma
\ref{lem:Det1}, we deduce that
\begin{equation} \label{eq:polydef3}
D(n,n-2,N)=
	(n-1)^N\Bigg(\prod_{i=1}^{N-1}\left((n-1)^2-\frac{i^2}{4}\right)^{N-i}
	\prod_{i=2}^{\floor{N/2}}\left((n-1)^2-i^2\right)\Bigg)\cdot R(n,N),
\end{equation}
where $R(n,N)$ is a polynomial in $n$ of degree at most $2\ceil{(N)/2}+2$
if $N\geq 2$.

Step~5 of the proof of Lemma \ref{lem:Det1} has a parallel here, too, which
yields evaluations of $R(n,N)$ at
$n=-\ceil{N/2}+1,\dots,0,2,\dots,\ceil{N/2}+1$.
In particular, the symmetry relation which plays the role of \eqref{eq:symP}
in this context is
\begin{equation}
\label{eq:symR}
R(-n+2,N)=R(n,N).
\end{equation}
While the evaluation of the determinant $D(e,e-2,2e-3)$ which is needed here
is done separately in Lemma~\ref{lem:detklein3} for $e>2$, 
the case $e=2$ leads to the $2\times2$ determinant
$\det_{1\le i,j\le 2}\left((-1)^{i+j}+(-2)^{i+j}\right)=4$.

The evaluation of $R(n,N)$ at $n=1$ requires extra treatment.
First, factoring $(n-1)$ out of each 
column of $D(n,n-2,N)$ and then setting $n=1$ yields 
the determinant $\det_{1\le i,j\le N}\big(2(B-1)^{i+j}\big)$
(using again the symbolic notation $B^k\equiv B_k$). 
By row and column operations,
this determinant can be transformed into the form
\begin{equation*}
2^N\det_{1\le i,j\le N}\big(B^{i+j-2}(B-1)^2\big).
\end{equation*}
This determinant agrees with the determinant on the left-hand side of
\eqref{eq:HankelBernoulli} with $n=N$, $a=b=0$ and $c=d=2$ (adopting
the usual convention $(B+1)_{-1}:=1/B$). 
Unfortunately, formula \eqref{eq:HankelBernoulli} does not hold for
this choice of parameters. In fact, the evaluation of this
determinant is rather tedious
and therefore given separately in Lemma~\ref{lem:extra-Bernoulli}.

Altogether, we have $2\ceil{N/2}+1$ evaluations for our polynomial
$R(n,N)$ at special values of $n$ so far. 
So we are short of exactly two informations on $R(n,N)$.
We can get these missing
informations by computing the leading coefficient of the polynomial,
which is done in the same way as in the proof of the preceding lemma,
and by exploiting the symmetry \eqref{eq:symR} once again.

For the simplifications of the resulting expression for $R(n,N)$ 
it turns out to be again convenient to separate the cases $N=2m$ and
$N=2m+1$.

\smallskip
This concludes the proof of the Lemma.
\end{proof}



\begin{lem}
\label{lem:detklein}
For $e\geq 1$, we have
\begin{equation}
\label{eq:detklein}
D(e,e,2e-1)=
	\det_{1\le i,j\le  2e-1}\left(\sum_{s=-e}^{e-1} s^{i+j}\right) =
\prod_{i=1}^{2e-1}i!^2.
\end{equation}
{\em(}The determinant $D(a,b,N)$ was defined in
\eqref{eq:gen-det}.{\em)}
\end{lem}

\begin{proof}
For $\ell \ge1$ define the polynomials $p_\ell (x)$ by
$p_\ell (x):=(x-\ceil{\ell /2}+1)_\ell $.
Furthermore, define coefficients $a_{\ell ,k}$, $k=1,2,\dots,\ell $, by
$p_\ell (x)=\sum _{k=1} ^{\ell }a_{\ell ,k}x^k$.
Now, for $\ell =2e-1,2e-2,\dots,2$, in this order, 
we replace column $c_\ell $ of the matrix underlying $D(e,e,2e-1)$ by
the linear combination of columns
$\sum _{k=1} ^{\ell }a_{\ell ,k}c_k,$
and afterwards do the same sort of replacement in the rows. These
operations yield a matrix with $(i,j)$-entry equal to
\begin{equation} \label{eq:Sum1}
\sum _{s=-e} ^{e-1}p_i(s)p_j(s)=
\sum _{s=-e} ^{e-1}(s-\ceil{i/2}+1)_i\,(s-\ceil{j/2}+1)_j.
\end{equation}

The next sequence of operations will turn this matrix into upper
triangular form, so that the determinant is easily obtained by
forming the product of the diagonal entries.

To begin with, it should be noted that for $j=2e-1$ the sum \eqref{eq:Sum1}
consists of just the term corresponding to $s=-e$. We subtract
$p_j(-e)/p_{2e-1}(-e)$ times column $2e-1$ from column $j$,
$j=1,2,\dots,2e-2$. The previous observation tells us that these
operations have the effect that the entries that are in one of the
columns $j=1,2,\dots,2e-2$ become
\begin{equation} \label{eq:Sum2}
\sum _{s=-e+1} ^{e-1}p_i(s)p_j(s)=
\sum _{s=-e+1} ^{e-1}(s-\ceil{i/2}+1)_i\,(s-\ceil{j/2}+1)_j,
\end{equation}
i.e., the summand corresponding to $s=-e$ has been eliminated. In
particular, all the entries in the last row, except the rightmost
entry, of course, are 0.

Next consider column $2e-2$. After the above column operations, the
sum \eqref{eq:Sum2} with $j=2e-2$ 
which defines the entries collapses to just the term
corresponding to $s=e-1$. We subtract
$p_j(e-1)/p_{2e-2}(e-1)$ times column $2e-2$ from column $j$,
$j=1,2,\dots,2e-3$. The previous observation tells us that these
operations have the effect that now the entries that are in one of the
columns $j=1,2,\dots,2e-3$ become
\begin{equation} \label{eq:Sum3}
\sum _{s=-e+1} ^{e-2}p_i(s)p_j(s)=
\sum _{s=-e+1} ^{e-2}(s-\ceil{i/2}+1)_i\,(s-\ceil{j/2}+1)_j,
\end{equation}
i.e., the summand corresponding to $s=e-1$ has been eliminated as well. In
particular, all the entries in the next-to-last row, except the two rightmost
entries, of course, are 0.

If we continue in the same manner, then eventually we arrive at an
upper triangular matrix whose $i$-th diagonal entry is $i!^2$. 
Thus, the result \eqref{eq:detklein} follows.
\end{proof}

\begin{lem}
\label{lem:detklein2}
For $e\geq 2$ we have
\begin{equation}
\label{eq:detklein2}
D(e,e-1,2e-2)=
	\det_{1\le i,j\le  2e-2}\left(\sum_{s=-e}^{e-2} s^{i+j}\right) =
	\prod_{i=1}^{2e-2}i!^2.
\end{equation}
{\em(}The determinant $D(a,b,N)$ was defined in
\eqref{eq:gen-det}.{\em)}
\end{lem}
\begin{proof}
This identity can be established in essentially the same manner as the
previous lemma. The basic difference is that one has to replace the
polynomials $p_\ell (x)$ in the previous proof by the polynomials
\begin{equation*}
q_\ell (x):=
\begin{cases}
	x			&\text{if $\ell =1$,}\\
	\po{x-\floor{\ell /2}+1}{\ell }	&\text{if $\ell >1$.}
\end{cases}
\end{equation*}
Everything else is completely analogous. We leave the details to the
reader.
\end{proof}

\begin{lem}
\label{lem:detklein3}
For $e\geq 3$ we have
\begin{equation}
\label{eq:detklein3}
D(e,e-2,2e-3)=
	\det_{1\le i,j\le  2e-3}\left(\sum_{s=-e}^{e-3} s^{i+j}\right) =
	\prod_{i=1}^{2e-3}i!^2.
\end{equation}
{\em(}The determinant $D(a,b,N)$ was defined in
\eqref{eq:gen-det}.{\em)}
\end{lem}
\begin{proof}
This identity can be established in essentially the same manner as the
previous lemmas. The basic difference is that one has to replace the
polynomials in the previous proofs by the polynomials
\begin{equation*}
r_\ell(x):=
\begin{cases}
	x			&\text{if $\ell =1$,}\\
	x(x+1)	&\text{if $\ell =2$,}\\
	\po{x-\ceil{\ell /2}+2}{\ell }	&\text{if $\ell >2$.}
\end{cases}
\end{equation*}
Everything else is completely analogous. We leave the details to the
reader.
\end{proof}


\section{Orthogonal polynomials, continued fractions, and Hankel
determinants of Bernoulli numbers}
\label{sec:orthpol}
In this section we review some facts about the interrelations between
orthogonal polynomials, continued fractions, and Hankel determinants.
Good sources for information about these topics are
\cite{JoThAA,KoSwAA,PerronKB,SzegoOP,VienAE,WallCF}.

To begin with, we recall Favard's Theorem.

\begin{thm} \label{thm:Favard}
{\em (Cf\@.  \cite[Th\'eor\`eme~9 on
p.~I-4]{VienAE} or \cite[Theorem 50.1]{WallCF}).}
Let $(p_n(x))_{n\ge0}$ be a sequence of monic polynomials, the
polynomial $p_n(x)$ having degree $n$. Then the sequence $(p_n(x))$ is
(formally) orthogonal with respect to a linear functional $L$, i.e.,
$L(p_n(x)p_m(x))=\delta_{mn}c_n$ for some sequence $(c_n)_{n\ge0}$ of
nonzero numbers, with $\delta_{m,n}$ denoting the Kronecker delta (i.e.,
$\delta_{m,n}=1$ if $m=n$ and $\delta_{m,n}=0$ otherwise) if and only
if there exist sequences $(a_n)_{n\ge1}$ and $(b_n)_{n\ge1}$, with
$b_n\ne0$ for all $n\ge1$, such that the three-term recurrence
\begin{equation}
p_{n+1}(x)=(a_{n}+x)p_{n}(x)-b_{n}p_{n-1}(x),\quad \quad 
\text{ for } n\geq 1,
\label{eq:three-term}
\end{equation}
holds, with initial conditions $p_0(x)=1$ and $p_1(x)=x+a_0$.
\end{thm}

It is a simple fact that, given a linear functional, the
corresponding orthogonal polynomials can be expressed in form of
certain determinants.

\begin{thm} \label{thm:ortho}
{\em (Cf\@. \cite[Theorem~7.15, (7.2.41)]{WallCF} or \cite[(36)]{VienAE}).}
Let $L$ be a linear functional defined on polynomials with moments
$\mu_n=L(x^n)$. Then the corresponding set 
$(p_n(x))_{n\ge0}$ of orthogonal polynomials is given by
\begin{equation}
\label{eq:ortho}
p_{n}(x)=d_n^{-1}\det\begin{pmatrix}
		\mu_0	&	\mu_1	&	\mu_2	&	\dots &	\mu_n \\
		\mu_1	&	\mu_2	&	\dots	&	\mu_n	&	\mu_{n+1} \\
		\mu_2	&	\dots	&	\mu_n	&	\mu_{n+1}&\mu_{n+2} \\
		\hdotsfor{5} \\
		\mu_{n-1} &	\mu_n	&	\mu_{n+1}&\dots	&\mu_{2n-1} \\
		1	&	x&\dots&	x^{n-1}&x^n 
	\end{pmatrix},
\end{equation}
where $d_n=\det_{0\le i,j\le n-1}(\mu_{i+j})$.
\end{thm}

The next theorem addresses the relation between
orthogonal polynomials and continued fractions, the link being the
generating function of the moments.
\begin{thm} \label{thm:momentgf}
{\em (Cf\@. \cite[Theorem 51.1]{WallCF} or \cite[Proposition~1, (7), on p.~V-5]{VienAE})}.
Let $(p_n(x))_{n\ge0}$ be a sequence of monic polynomials, the
polynomial $p_n(x)$ having degree $n$, which is orthogonal with
respect to some functional $L$. Let
\begin{equation}
p_{n+1}(x)=(a_{n}+x)p_{n}(x)-b_{n}p_{n-1}(x)
\label{eq:three-term2}
\end{equation}
be the corresponding three-term recurrence which is guaranteed by
Favard's theorem. Then the generating function for the moments
$\mu_n=L(x^n)$ satisfies
\begin{equation}
\label{eq:momentgf}
\sum_{n=0}^\infty{\mu_n}x^n=\cfrac{
	\mu_0}
		{1+a_0x-\cfrac{
			b_1x^2}
				{1+a_1x-\cfrac{
					b_2x^2}
						{1+a_2x-\cdots}}}\quad .
\end{equation}
\end{thm}

We remark that a continued fraction of the type as in
\eqref{eq:momentgf} is called a {\em J-fraction}.

The next theorem addresses the relation between J-fractions and
Hankel determinants.
\begin{thm}
\label{cor:cfracHankel}
{\em (Cf\@. \cite[Theorem 51.1]{WallCF} or
\cite[Corollaire~6, (19), on p.~IV-17]{VienAE})}.
Let $(\mu_n)_{n\ge0}$ be a sequence of numbers whose generating
function, when written in terms of a J-fraction, is given by
\eqref{eq:momentgf}. 
Then the Hankel determinant $\det_{0\le i,j\le n-1}(\mu_{i+j})$
equals $\mu_0^nb_1^{n-1}b_2^{n-2}\cdots b_{n-2}^2b_{n-1}$.
\end{thm}

The family of orthogonal polynomials that is of interest to us is the
{\em continuous Hahn polynomials}, which were first studied by
Atakishiyev and Suslov \cite{AtSuAA} and Askey \cite{AskeyCHP}. 
The theorem below lists the
relevant facts about these polynomials.
\begin{thm} \label{thm:cHahn}
{\em (Cf\@. \cite[Section 1.4]{KoSwAA})}.
The continuous Hahn polynomials $(p_n(a,b,c,d;x))_{n\ge0}$ are
the monic polynomials defined as the terminating hypergeometric series
\begin{multline} \label{eq:cHahn}
p_n(a,b,c,d;x)\\
=\frac{\complexikl^n \,\po{a+c}{n}\,\po{a+d}{n}}{\po{a+b+c+d+n-1}{n}}
\HypF{3}{2}{-n,n+a+b+c+d-1,a+x\complexi }{a+c,a+d}{1}.
\end{multline}
These polynomials satisfy the three-term recurrence
\begin{equation} \label{eq:cHahnrec}
p_{n+1}(a,b,c,d;x)  
=	(x-A_n(a,b,c,d))p_{n}(a,b,c,d;x)-B_n(a,b,c,d)p_{n-1}(a,b,c,d;x),
\end{equation}
where
\begin{multline}
A_n(a,b,c,d) = \complexi\,\Biggl( a + {\frac{n\,\left(  b + c + n-1 \right) \,
        \left(  b + d + n-1 \right) }{\left(  a + b + c + d + 
          2 n-2 \right) \,\left(  a + b + c + d + 2 n-1 \right) }} + \\
    {\frac{\left( 1 - a - b - c - d - n \right) \,
        \left( a + c + n \right) \,\left( a + d + n \right) }{
        \left( -1 + a + b + c + d + 2 n \right) \,
        \left( a + b + c + d + 2 n \right) }} \Biggr),
\end{multline}
and
\begin{multline}
B_n(a,b,c,d) = -{\frac{n\,\left(a + c + n - 1\right) \,
     \left(b + c + n - 1\right) \,\left(a + d + n - 1\right)}
	  {\left(a + b + c + d + 
       2 n - 3\right) \,{{\left(a + b + c + d + 2 n - 2\right) }^2}}}\\
\times{\frac{\left(b + d + n - 1\right) \,
     \left(a + b + c + d + n - 2\right) }{
     \left(a + b + c + d + 2 n - 1\right) }}.
\end{multline}
They are orthogonal with respect to the functional $L$ which is given
by
\begin{equation} \label{eq:cHahnortho}
L(p(x))=\frac{1}{2\pi}\int_{-\infty}^{\infty}
\Gamma(a+x\complexi )\,\Gamma(b+x\complexi )
\,\Gamma(c-x\complexi )\,\Gamma(d-x\complexi )\,p(x)\,dx\,.
\end{equation}
In this integral, if necessary, the path of integration has to be 
deformed, so that it separates the increasing sequences
$\big((a+k)\complexi\big)_{k\ge0}$ and 
$\big((b+k)\complexi\big)_{k\ge0}$ of poles from
the decreasing sequences 
$\big(-(c+k)\complexi\big)_{k\ge0}$ and
$\big(-(d+k)\complexi\big)_{k\ge0}$ of poles (see \cite{AskeyCHP}).

Explicitly, the orthogonality relation is
\begin{multline}
L(p_m(a,b,c,d;x)p_n(a,b,c,d;x)) \\
=\frac{
\Gamma(n+a+c)\,\Gamma(n+a+d)\,\Gamma(n+b+c)\,\Gamma(n+b+d)\,n!
}{
\po{a+b+c+d+n-1}{n}\,\po{a+b+c+d+n-1}{n+1}\Gamma(n+a+b+c+d-1)
}\delta_{m,n}.
\end{multline}
In particular, 
\begin{equation}
L(1)={\frac{\left(a + c -1\right) !\,\left(b + c - 1\right) !\,
     \left(a + d - 1\right) !\,\left(b + d - 1\right) !}{\left(
       a + b + c + d - 1\right) !}}.
\end{equation}
\end{thm}
\begin{rem}
(1)
The reader should be aware that Theorem~\ref{thm:cHahn} is formulated for 
the {\em monic} 
form of the continuous Hahn polynomials, so that our polynomials are
not the polynomials which are denoted by $p_n(x;a,b,c,d)$ in 
\cite[Section 1.4]{KoSwAA} but those which are denoted there, slightly
confusingly, by $p_n(x)$.

(2) Although the definition \eqref{eq:cHahn} does not show it, the
continuous Hahn polynomials are symmetric in $a$ and $b$, and in $c$
and $d$, because the orthogonality
measure \eqref{eq:cHahnortho} has
these symmetries. In addition, there is another symmetry when the
roles of $a,b$ are interchanged with the roles of $c,d$. Namely, we
have 
\begin{equation} \label{eq:cHahn-symmetry}
p_n(a,b,c,d;x)= (-1)^np_n(c,d,a,b;-x).
\end{equation}
This follows from the fact that, when performing the substitution
$x\to-x$ in the
integral \eqref{eq:cHahnortho}, we obtain almost the same integrand,
the only differences being that the roles of $a,b$ and $c,d$ are
interchanged, and that $p(x)$ is replaced by $p(-x)$. 
The path of integration also does not change, at least as long as
all of $a,b,c,d$ are positive. This is, however, sufficient to
conclude that \eqref{eq:cHahn-symmetry} holds for all $a,b,c,d$
because, for fixed $n$, both sides of \eqref{eq:cHahn-symmetry} are
rational of bounded degree in $a,b,c,d$.
(The same conclusion could also be less elegantly derived by applying some
$_3F_2$-transformation formulas.)
\end{rem}

Now, by combining Theorems~\ref{thm:momentgf}, \ref{cor:cfracHankel},
and \ref{thm:cHahn} we are able to derive without difficulty the
determinant evaluation that we need in Step~5 of the proof of
Lemma~\ref{lem:Det1} and the analogous places in Lemmas~\ref{lem:Det2}
and \ref{lem:Det3} (which, in turn, are essential for the proofs of
Theorems~\ref{thm:MEven}--\ref{thm:NOdd2}). 

\begin{thm} \label{thm:HankelBernoulli}
For positive integers $a,b$ and nonnegative integers $c,d$ there holds
\begin{multline} \label{eq:HankelBernoulli}
\det_{1\le i,j\le n}
\left(
	B^{i+j}\,\po{B+1}{a-1}\,\po{B+1}{b-1}\,\po{-B+1}{c-1}\,\po{-B+1}{d-1}
\right)\\
=(-1)^{n(n-1)/2}\left({\frac{\left(a + c -1\right) !\,\left(b + c - 1\right) !\,
     \left(a + d - 1\right) !\,\left(b + d - 1\right) !}{\left(
       a + b + c + d - 1\right) !}}\right)^n\\
\times\prod_{i=1}^{n-1}\Biggl({\frac{i\,\left(a + c + i - 1\right) \,
     \left(b + c + i - 1\right) \,\left(a + d + i - 1\right)}
	  {\left(a + b + c + d + 
       2 i - 3\right) \,{{\left(a + b + c + d + 2 i - 2\right) }^2}}}\\
\times{\frac{\left(b + d + i - 1\right) \,
     \left(a + b + c + d + i - 2\right) }{
     \left(a + b + c + d + 2 i - 1\right) }}
\Biggr)^{n-i},
\end{multline}
where we have again used the symbolic notation $B^k\equiv B_k$. (In
case that $c$ or $d$ should be zero, we have to interpret $(-B+1)_{-1}$
as $1/(-B)$.)
\end{thm}
\begin{proof}
We claim that the moments for the continuous Hahn
polynomials, i.e., for the linear functional as defined in
\eqref{eq:cHahnortho}, are
\begin{equation}
\label{eq:CHmoments}
\mu_n=(B/\complexi)^n B^2
\,\po{B+1}{a-1}\,\po{B+1}{b-1}\,\po{-B+1}{c-1}\,\po{-B+1}{d-1}.
\end{equation}
Let us suppose
for the moment that this is already established. Then, by
Theorem~\ref{thm:momentgf} and \eqref{eq:cHahnrec} of Theorem~\ref{thm:cHahn},
we have
\begin{multline*}
\sum_{n\geq 0}{(B/\complexi)^n B^2
\,\po{B+1}{a-1}\,\po{B+1}{b-1}\,\po{-B+1}{c-1}\,\po{-B+1}{d-1}}x^n\\ 
=\cfrac{B^2
\,\po{B+1}{a-1}\,\po{B+1}{b-1}\,\po{-B+1}{c-1}\,\po{-B+1}{d-1}	}
		{1+A_0(a,b,c,d)x-\cfrac{
			B_1(a,b,c,d)x^2}
				{1+A_1(a,b,c,d)x-\cfrac{
					B_2(a,b,c,d)x^2}
						{1+A_2(a,b,c,d)x-\cdots}}}\quad .
\end{multline*}
But, by Theorem~\ref{cor:cfracHankel},
this immediately implies the truth of \eqref{eq:HankelBernoulli}.

It remains to verify \eqref{eq:CHmoments}. Using
\eqref{eq:cHahnortho}, this is a rather straightforward computation:
\begin{align}
\mu_n &=
	\frac{1}{2\pi}\int_{-\infty}^{\infty}
	\Gamma(a+x\complexi )\,\Gamma(b+x\complexi )
	\,\Gamma(c-x\complexi )\,\Gamma(d-x\complexi )\,x^n\,dx\notag\\
&\kern-5pt= \frac{1}{2\pi}\int_{-\infty\complexi }^{\infty\complexi }
	\Gamma(a+z)\,\Gamma(b+z)\,\Gamma(c-z)\,\Gamma(d-z)
	\frac{z^n}{\complexikl ^{n+1}}\,dz\notag\\
&\kern-5pt= \frac{1}{2\pi \complexikl^{n+1} }
	\int_{-\infty\complexi }^{\infty\complexi }
	\po{z}{a}\,\po{z}{b}\,\po{-z}{c}\,\po{-z}{d}\left(\Gamma(z)\,\Gamma(-z)\right)^2
	z^n\,dz\notag\\
&\kern-5pt= \frac{1}{2\pi \complexikl^{n+1} }
	\int_{-\infty\complexi }^{\infty\complexi }
	\po{z}{a}\,\po{z}{b}\,\po{-z}{c}\,\po{-z}{d}
		\left(-\frac{\pi}{z\sin \pi z}\right)^2 z^n\,dz\notag\\
&\kern-5pt= \frac{1}{2\pi \complexikl^{n+1} }
	\int_{-\infty\complexi }^{\infty\complexi }
	\po{z+1}{a-1}\,\po{z+1}{b-1}\,\po{-z+1}{c-1}\,\po{-z+1}{d-1}
		\left(\frac{\pi}{\sin\pi z}\right)^2
	z^{n+2}\,dz.\label{eq:momentn}
\end{align}
In the third line we used the relation $\Gamma(a+z)=\po{z}{a}\,\Gamma(z)$
for $a\in {\mathbb N}$ 
(see, e.g., \cite[1.2(2)]{Erdelyi}), and in the fourth line we used 
the formula
$\Gamma(z)\,\Gamma(-z)=-\frac{\pi}{z\sin\pi z}$ (see \cite[1.2(5)]{Erdelyi}).
The reader should notice that, because of the convention regarding
the path of integration in \eqref{eq:cHahnortho}, in case that $c$ or
$d$ are zero the path of
integration in the last line of
\eqref{eq:momentn} is deformed so that it crosses the real axis
between the singularities $z=-1$ and $z=0$.

To finish the calculation, we appeal to the following integral
representation of Bernoulli numbers (see
\cite[p.~75]{Noerlund}, and let $\alpha\to0^+$ in the relevant
identity in the middle of the page)
\begin{equation}
B_\nu=\frac{1}{2\pi \complexi }
\int_{-\infty\complexi }^{\infty\complexi }z^\nu
\left(\frac{\pi}{\sin\pi z}\right)^2 dz.
\end{equation}
(If $\nu=0$ or $\nu=1$ then the path of integration is indented so
that it avoids the singularity $z=0$, passing it on the negative
side, but to the right of the singularity $z=-1$.)
If we use this formula in \eqref{eq:momentn}, 
we obtain \eqref{eq:CHmoments} immediately. 
\end{proof}

At least two special cases of Theorem~\ref{thm:HankelBernoulli} have
explicitly appeared in the literature before. The case $a=b=1$,
$c=d=0$ appears for example in \cite[(3.1)]{AlCaAA}. The case
$a=b=c=d=1$ appears in \cite[App.~A.5, p.~322]{MehtAB}.

\medskip
In the proof of Theorems~\ref{thm:NEven2} and \ref{thm:NOdd2} (to be
precise, in the proof of Lemma~\ref{lem:extra-Bernoulli}), 
we make use of a rather recent result on (formal)
orthogonal polynomials, due to Leclerc \cite[Theorem~1]{Leclerc}.

\begin{thm} \label{thm:Leclerc}
For an arbitrary sequence of numbers $(\mu_n)_{n\ge0}$ let
$(P_n(x))_{n\ge0}$ be the sequence of polynomials defined by
\begin{equation}
\label{eq:defineP}
P_{n}(x):=\det\begin{pmatrix}
		\mu_0	&	\mu_1	&	\mu_2	&	\dots &	\mu_n \\
		\mu_1	&	\mu_2	&	\dots	&	\mu_n	&	\mu_{n+1} \\
		\mu_2	&	\dots	&	\mu_n	&	\mu_{n+1}&\mu_{n+2} \\
		\hdotsfor{5} \\
		\mu_{n-1} &	\mu_n	&	\mu_{n+1}&\dots	&\mu_{2n-1} \\
		1	&	x&\dots&	x^{n-1}&x^n 
	\end{pmatrix}.
\end{equation}
(These are, up to normalization, the orthogonal polynomials
with respect to the linear functional with moments $(\mu_n)$; 
compare Theorem~\ref{thm:ortho}.)
Furthermore, let
$(Q_n(x))_{n\ge0}$ be the sequence of polynomials defined by
\begin{equation*}
Q_{n}(x):=\sum_{k=0}^n\mu_k\binom{n}{k}(-x)^{n-k}.
\end{equation*}
Then, for arbitrary
integers $l\geq1$, $m\geq1$, there holds
\begin{multline}
\label{eq:Leclerc}
\det\begin{pmatrix}
P_{l}(x) & P_{l+1}(x) &\dots & P_{l+m-1}(x) \\
P_{l}^\prime(x) & P_{l+1}^\prime(x) 
	&\dots& P_{l+m-1}^\prime(x) \\
\vdots & \vdots & & \vdots \\
P_{l}^{(m-1)}(x) & P_{l+1}^{(m-1)}(x) &
	\dots& P_{l+m-1}^{(m-1)}(x)
\end{pmatrix} \\
= C_{l,m}
\det\begin{pmatrix}
Q_{m}(x) & Q_{m+1}(x) &\dots& Q_{m+l-1}(x) \\
Q_{m+1}(x) & Q_{m+2}(x) &\dots& Q_{m+l}(x) \\
\vdots & \vdots & & \vdots \\
Q_{m+l-1}(x) & Q_{m+l}(x) &\dots& Q_{m+2l-2}(x)
\end{pmatrix},
\end{multline}
where
\begin{equation}
\label{eq:Cln}
C_{l,m}=(-1)^{lm}\prod_{k=1}^{m-1}k!
\det\begin{pmatrix}
\mu_0 & \mu_1 &\dots& \mu_{k+l-1} \\
\mu_1 & \mu_2 &\dots& \mu_{k+l} \\
\vdots & \vdots & & \vdots \\
\mu_{k+l-1} & \mu_{k+l} &\dots& \mu_{2k+2l-2} \\ 
\end{pmatrix}.
\end{equation}
\end{thm}


\section{Auxiliary lemmas}
\begin{lem}
\label{lem:extra-Bernoulli}
For any integer $n\geq 1$, we have
\begin{multline}
\label{eq:extra-Bernoulli}
\det_{1\leq i,j\leq n}\left(B^{i+j-2}(B-1)^2\right)\\
=
{(-1)^{{n(n-1)}/{2}}} {6^{-n}}(1 + 2n + 5n^2 + 4n^3 + n^4)
\prod_{i=1}^n\left(\frac{i(i+1)^4(i+2)}{(2i+1)(2i+2)^2(2i+3)}\right)^{n-i}
\end{multline}
(again, using the symbolic notation $B^k\equiv B_k$).
\end{lem}
\begin{proof}
The determinant in \eqref{eq:extra-Bernoulli} is equal to the determinant in
\eqref{eq:HankelBernoulli} with $a=b=0$ and $c=d=2$ (again, with the
convention that $(B+1)_{-1}$ is interpreted as $1/B$). Unfortunately,
for this choice of parameters, formula \eqref{eq:HankelBernoulli} is
not valid. However, the determinant in \eqref{eq:extra-Bernoulli} is
very close to the determinant in \eqref{eq:HankelBernoulli} with 
$a=b=2$ and $c=d=0$. In fact, because of the well-known property of
Bernoulli numbers that $B_{2k+1}=0$ for all positive integers $k$,
multiplication of all even numbered rows and columns of the
determinant in \eqref{eq:HankelBernoulli} by $-1$ (which does not change the
value of the determinant) with this choice of
parameters transforms the latter 
into the determinant $\Delta_n:=\det_{1\le i,j\le n}(\la_{i+j-2})$, where
\begin{align} \notag
\la_0	&=B^0(B+1)^2= B^{0}(B-1)^2-2,\\
\notag
\la_1	&=-B^1(B+1)^2= B^{1}(B-1)^2+1,\\
\label{eq:la}
\la_n	&=(-1)^nB^n(B+1)^2= B^{n}(B-1)^2\text{ for }n\geq2.
\end{align}
So, because of the deviating definitions of $\la_0$ and $\la_1$, 
the only difference between $\Delta_n$ and the determinant in
\eqref{eq:extra-Bernoulli} is in the top-left entry and its right and
bottom neighbour. Linearity of the determinant $\Delta_n$ 
in the first row and column then implies that
\begin{equation} \label{eq:extra-expansion}
\det_{1\leq i,j\leq n}\left(B^{i+j-2}(B-1)^2\right)=
\Delta_{n} -
\Delta_{n}^{\{1,2\};\{1,2\}} +
2\Delta_{n}^{\{1\};\{1\}} +
2\Delta_{n}^{\{1\};\{2\}}.
\end{equation}
Here, $A^{\{i_1,i_2,\dots\};\{j_1,j_2\dots\}}$ denotes the minor of $A$ with
rows $i_1,i_2,\dots$ and columns $j_1,j_2\dots$ deleted. (Empty minors
are defined to be zero).

Now observe that from Theorem~\ref{thm:HankelBernoulli} with $a=b=2$
and $c=d=0$ we obtain immediately that
\begin{equation} \label{eq:Delta}
\Delta_{n}=(-1)^{{n(n-1)}/2}
\prod_{i=1}^{n}\frac{(i-1)!\,i!^4\,(i+1)!}{(2i)!\,(2i+1)!}.
\end{equation}
Likewise, from Theorem~\ref{thm:HankelBernoulli} with $a=b=1$ and
$c=d=2$ we have
\begin{equation} \label{eq:Delta-1,1}
\Delta_{n}^{\{1\};\{1\}}=(-1)^{{(n-1)(n-2)}/2}
\prod_{i=1}^{n-1}\frac{(i-1)!\,(i+1)!^4\,(i+3)!}{(2i+2)!\,(2i+3)!}.
\end{equation}

Our next observation is that $\Delta_{n}^{\{1\};\{2\}}$ is,
essentially, the coefficient of $x$ in the continuous Hahn polynomial
$p_{n-1}(0,1,2,2;x)$. To make a more precise statement, consider
\eqref{eq:ortho} with $\mu_k=\la_{k+1}$, $k=0,1,\dots$. Then,
obviously, $\Delta_{n}^{\{1\};\{2\}}$ equals
$$(-1)^n\det_{1\le i,j\le n-1}(\la_{i+j-1})
\cdot\big(\text {coefficient of $x$ in } p_{n-1}(x)\big).$$
On the other hand, Theorem~\ref{thm:ortho} says that, with this
choice of the $\mu_k$, the polynomials $p_n(x)$ are orthogonal with
associated moments $\la_{k+1}$, $k=0,1,\dots$. By comparing
\eqref{eq:la} and \eqref{eq:CHmoments} with $a=b=2$, $c=0$, and
$d=1$, we see that we must have
$$p_n(x)=(-\complexi)^n\,p_n(2,2,0,1;-x/\complexi).$$
Therefore, if we remember \eqref{eq:cHahn-symmetry},
we obtain that $\Delta_{n}^{\{1\};\{2\}}$ equals
\begin{equation} \label{eq:0122}
(-1)^n(\complexi)^{n-1}\det_{1\le i,j\le n-1}(\la_{i+j-1})
\cdot\big(\text {coefficient of $x$ in }p_{n-1}(0,1,2,2;x/\complexi)\big).
\end{equation}
The Hankel determinant in this expression can be evaluated 
by using Theorem~\ref{thm:HankelBernoulli} with $a=b=2$, $c=0$, $d=1$, 
and $n$ replaced by $n-1$. By substituting the result in
\eqref{eq:0122} and by using the definition \eqref{eq:cHahn} of
continuous Hahn polynomials with $a=0$, $b=1$, $c=d=2$, we obtain
\begin{multline} \label{eq:Delta-1,2}
\Delta_{n}^{\{1\};\{2\}}=
(-1)^{n(n+1)/2}\frac{n!^2}{(n+3)_{n-1}}
\\ \times
\bigg(\prod_{i=1}^{n-1}\frac{(i-1)!\,i!^2\,(i+1)!^2\,(i+2)!}
{(2i+1)!\,(2i+2)!}\bigg)
\sum_{k=1}^{n-1}
	\frac{\po{1-n}{k}\,\po{n+3}{k}}{k\,(k+1)!^2}.
\end{multline}

\def\a{a}

The remaining minor $\Delta_{n}^{\{1,2\};\{1,2\}}$ requires additional
work. We employ Theorem~\ref{thm:Leclerc} with
$l=n-2$, $m=2$, $x=0$, and $\mu_k=\la_{k+2}$, $k=0,1,\dots$.
With this choice of parameters, the determinant on the
right-hand side of \eqref{eq:Leclerc} is precisely our remaining minor
$\Delta_{n}^{\{1,2\};\{1,2\}}$. The single determinant occurring in the
definition \eqref{eq:Cln} of $C_{n-2,2}$ can be evaluated by
using Theorem~\ref{thm:HankelBernoulli} with $a=b=1$, $c=d=2$, and $n$
replaced by $n-2$, so that we have
\begin{equation}
\label{eq:Cl2}
C_{n-2,2}=(-1)^{(n-1)(n-2)/2}
\prod_{i=1}^{n-1}\frac{(i-1)!\,(i+1)!^4\,(i+3)!}{(2i+2)!\,(2i+3)!}.
\end{equation}
The determinant on the left-hand side \eqref{eq:Leclerc} 
is a $2\times2$-determinant
whose entries are the constant term and the coefficient of $x$, respectively,
of $P_{n-2}(x)$ and $P_{n-1}(x)$. The polynomials $P_n(x)$,
$n=0,1,\dots$, on the other hand, are orthogonal polynomials with
associated moments $\la_{k+2}$, $k=0,1,\dots$ (compare \eqref{eq:defineP} 
and Theorem~\ref{thm:ortho}). By comparing the definition
\eqref{eq:la} of the $\la_i$'s with \eqref{eq:CHmoments}, it is seen
that the polynomials $P_n(x)$ must agree, up to
normalization, with the continuous Hahn polynomials with parameters
$a=b=1$ and $c=d=2$. To be precise, we have
$$P_n(x)=(\complexi)^n\det_{1\le i,j\le
n}(\la_{i+j})\,p_n(1,1,2,2;x/\complexi) .$$
Clearly, the Hankel
determinant on the right-hand side can be once again evaluated
by means of Theorem~\ref{thm:HankelBernoulli} with $a=b=1$ and
$c=d=2$.
In summary, from \eqref{eq:Leclerc} with the above choice of
parameters we infer
\begin{multline} \label{eq:Delta-12,12}
\Delta_{n}^{\{1,2\};\{1,2\}}=
	(-1)^{n(n-1)/2}\frac{\po{3}{n-2}^2\,\po{3}{n-1}^2}
        {\po{n+3}{n-2}\,\po{n+4}{n-1}}
	\left(\prod_{i=1}^{n-2}\frac{(i-1)!\,(i+1)!^4\,(i+3)!}
        {(2i+2)!\,(2i+3)!}\right)
	\\ \times
	\left(c_{0,n-2}c_{1,n-1}-c_{0,n-1}c_{1,n-2}\right),
\end{multline}
where $c_{0,n}$ and $c_{1,n}$ denote the coefficient of $x^0$ and $x^1$,
respectively, in
\begin{equation*}
\sum_{k=0}^{n}\frac{\po{-n}{k}\,\po{n+5}{k}\,\po{1+x}{k}}{k!\,\po{3}{k}^2}.
\end{equation*}
Using hypergeometric notation \eqref{eq:hypergeom}, 
the first of these two, $c_{0,n}$, can be expressed as
$$
  {\frac {4} 
    {\left( 1 + n \right)  \left( 2 + n \right)  \left( 3 + n \right)  
      \left( 4 + n \right) }}
\left(
{} _{2} F _{1} \!\left [ \begin{matrix} { -2 - n, 3 + n}\\ {1}\end{matrix} ;
       {\displaystyle 1}\right ]
-   1 + (n+2)(n+3)\right).
$$
The $_2F_1$-series can be evaluated by means of the Chu--Vandermonde
summation formula (see \cite[(1.7.7); Appendix (III.4)]{SlatAC})
\begin{equation} \label{eq:Chu-Van}
{} _{2} F _{1} \!\left [ \matrix { a, -N}\\ { c}\endmatrix ; {\displaystyle
   1}\right ]  = {\frac {({ \textstyle c-a}) _{N} } 
    {({ \textstyle c}) _{N} }},
\end{equation} 
where $N$ is a nonnegative integer. This yields
\begin{equation} \label{eq:c}
c_{0,n}=\begin{cases}
	\frac{4}{(n+1)(n+4)} &\text{ if $n$ is even,}\\
	\frac{4}{(n+2)(n+3)} &\text{ if $n$ odd.}
\end{cases}
\end{equation}
By combining \eqref{eq:Delta-1,2}, \eqref{eq:Delta-12,12} and
\eqref{eq:c}, we get, after a considerable amount of simplification,
\begin{multline*}
2\Delta_{n}^{\{1\};\{2\}}-\Delta_{n}^{\{1,2\};\{1,2\}}=
        (-1)^{\binom {n-1}2}
	\prod_{i=1}^{n}\frac{(i-1)!\,i!^4\,(i+1)!}{(2i)!\,(2i+1)!}\\
	\times\Biggl(\sum_{k=0}^{n-2}
	\frac{(-1)^k\,(n+k+3)!\,(1+(-1)^n(n+1))}{(k+1)\,(k+2)!^2\,(n-k-2)!}\\
+       \frac {(n+2)!} {(n-1)!}
	\sum_{k=0}^{n-2}
	\frac{(-1)^{n+k-1}(n+k+3)!}{(k+2)!\,(k+3)!\,(n-k-2)!}
        \sum_{j=0}^{k}\frac {1} {j+1}
\Biggr).
\end{multline*}
By Lemma~\ref{lem:weird}, with $n$ replaced by $n-1$, 
this expression reduces to
\begin{equation} \label{eq:Delta-diff}
2\Delta_{n}^{\{1\};\{2\}}-\Delta_{n}^{\{1,2\};\{1,2\}}=
(-1)^{\binom {n+1}2}
((-1)^n(n+1)+2)\frac {(n+2)!} {(n-1)!}
\prod_{i=1}^{n}\frac{(i-1)!\,i!^4\,(i+1)!}{(2i)!\,(2i+1)!}.
\end{equation}
Substituting \eqref{eq:Delta}, \eqref{eq:Delta-1,1}, and
\eqref{eq:Delta-diff} in \eqref{eq:extra-expansion}, and simplifying the
resulting expression, we eventually arrive at \eqref{eq:extra-Bernoulli}.
\end{proof}


\begin{lem}
\label{lem:weird}
For $n\geq 0$, we have
\begin{multline}
\label{eq:weird}
\frac{n!}{(n+3)!}\sum_{k=0}^{n-1}
	\frac{(-1)^k\,(n+k+4)!\,(1-(-1)^n(n+2))}{(k+1)\,(k+2)!^2\,(n-k-1)!}\\
+        \sum_{k=0}^{n-1}
	\frac{(-1)^{n+k}(n+k+4)!}{(k+2)!\,(k+3)!\,(n-k-1)!}
        \sum_{j=0}^{k}\frac {1} {j+1}=(-1)^n(n+2)-2.
\end{multline}
\end{lem}

\begin{proof}
We shall treat the two sums (in the first and the second line in
\eqref{eq:weird}, respectively) separately.

First we consider the sum in the first line of \eqref{eq:weird}.
We replace the term $1/(k+1)$ by $\int_0^1x^k\dx$, interchange
summation and integration, and write the sum
in hypergeometric notation \eqref{eq:hypergeom}. This gives
\begin{equation*}
\frac{1-(-1)^n(n+2)}{\po{n+1}{3}}\int_0^1\left(\frac{n+2}{x^2}
\HypF{3}{2}{n+3,-n-1}{1}{x}-\frac{n+2}{x^2}+\frac{\po{n+1}{3}}{x}
\right)\,{dx}.
\end{equation*}
Using the transformation formula (see \cite[(1.7.1.3), sum
reversed at the right-hand side]{SlatAC})
\begin{equation} \label{eq:T2106}
{} _{2} F _{1} \!\left [ \matrix { a, -N}\\ { c}\endmatrix ; {\displaystyle
   z}\right ]  = {z^N}
  {\frac {  ({ \textstyle c-a}) _{N} } {({ \textstyle c}) _{N} }}
  {} _{2} F _{1} \!\left [ \matrix { -N, 1 - c - N}\\ {
       1 + a - c - N}\endmatrix ; {\displaystyle {\frac {z-1} z}}\right ] 
\end{equation} 
(where $N$ is a nonnegative integer) with $a=n+3$, $N=n+1$, 
$c=1$, and $z=x$, this is transformed into
\begin{multline} \label{eq:ex1}
\frac{1-(-1)^n(n+2)}{\po{n+1}{3}}
\int_0^1\Biggl((n+2)^2(-1)^{n+1}
\sum_{k=0}^{n+1}\frac{(-1)^k\,\po{-n-1}{k}^2}{k!\,(k+1)!}(1-x)^k\,x^{n-k-1}\\
-\frac{n+2}{x^2}+\frac{\po{n+1}{3}}{x}
\Biggr)\,{dx}.
\end{multline}
Recall that for $\Re(\alpha)>-1$ and $\Re(\beta)>-1$ we have the following
integral representation for the {\em Euler beta function},
\begin{equation}
\label{eq:beta-integral}
\int_0^1(1-x)^\alpha x^\beta {\, dx}=
\frac{\Gamma(\alpha+1)\Gamma(\beta+1)}{\Gamma(\alpha+\beta+2)}.
\end{equation}
Use of this identity in \eqref{eq:ex1} 
wherever possible (i.e., it is applied to the summands with 
$k=0,\dots,n-1$) yields
after some simplification the expression
\begin{multline} \label{eq:ex2}
\frac{1-(-1)^n(n+2)}{\po{n+1}{3}}\Biggl(
(-1)^{n+1}\left(1-\sum_{k=0}^n\frac{\po{-n-2}{k}^2}{\po{-n}{k}\,k!}\right)
\\+\int_0^1\left(	
	\frac{\po{n+1}{2}\left((n+3)-(n+2)(1-x)^n\right)}{x}+
	\frac{(n+2)\left((1-x)^{n+1}-1\right)}{x^2}
\right){\, dx}\Biggr).
\end{multline}
We would like to write the sum in the first line as a hypergeometric
series. Unfortunately, this cannot be done by just straighforwardly
extending the summation over all nonnegative $k$ because of the term
$(-n)_k$ in the denominator, which is 0 for $k=n+1$. The way to
overcome this problem is to rewrite the sum as a limit,
\begin{multline*} 
\sum_{k=0}^n\frac{\po{-n-2}{k}^2}{\po{-n}{k}\,k!}\\
=\lim_{\ep\to0}\Bigg({}_2F_1\!\left[\matrix -n-2,n-2\\-n-\ep\endmatrix;
1\right]
-\frac {(-n-2)_{n+1}^2} {(n+1)!\,(-n-\ep)_{n+1}}
-\frac {(-n-2)_{n+2}^2} {(n+2)!\,(-n-\ep)_{n+2}}\Bigg).
\end{multline*}
The $_2F_1$-series can be evaluated by means of the Chu--Vandermonde
summation formula \eqref{eq:Chu-Van}. Substituing the result in
\eqref{eq:ex2}, we obtain
\begin{multline*}
\hskip-8pt
\frac{1-(-1)^n(n+2)}{\po{n+1}{3}}\Biggl(
(-1)^{n+1}\left(
	1-\lim_{\varepsilon\rightarrow0}\left(
	\frac{\po{2-\varepsilon}{n+2}-(n+2)\,(n+2)!\,(1-\varepsilon)-(n+2)!}{
		\po{-n-\varepsilon}{n+2}
	}
\right)\right)\\
+\int_0^1\left(	
	\frac{\po{n+1}{2}\left((n+3)-(n+2)\,(1-x)^n\right)}{x}+
	\frac{(n+2)\left((1-x)^{n+1}-1\right)}{x^2}
\right){\, dx}
\Biggr).
\end{multline*}
Of course, the limit can be computed by de l'Hospital's rule, so that
we get
\begin{multline*}
\frac{1-(-1)^n(n+2)}{\po{n+1}{3}}\Biggl(
(-1)^{n+1}
	+\frac{1}{n!}\left(
	(n+3)!\sum_{k=1}^{n+2}\frac{1}{k+1}-(n+2)!\,(n+2)
\right)\\
+\int_0^1\left(	
	\frac{\po{n+1}{2}\left((n+3)-(n+2)(1-x)^n\right)}{x}+
	\frac{(n+2)\left((1-x)^{n+1}-1\right)}{x^2}
\right){\, dx}
\Biggr).
\end{multline*}
Now, let us turn to the integral. Expanding the integrand by the binomial
theorem and simplifying the result, we obtain
\begin{multline} 
\frac{1-(-1)^{n}(n+2)}{\po{n+1}{3}}\left(
	(-1)^{n+1}-{(n+1)\,(n+2)^2}+
        {(n+1)_3}\sum_{k=1}^{n+2}\frac{1}{k+1}
	\right)\\
\label{eq:ex3}
+\frac{1-(-1)^{n}(n+2)}{\po{n+1}{3}}(n+2)\int_0^1\left(
	\sum_{k=0}^{n-1}\frac{\po{-n-1}{k+2}-(k+2)\po{-n-2}{k+3}}{(k+2)!}x^k
\right){\, dx}.
\end{multline}
Now the integration can be performed without any difficulty. 
We rewrite the resulting sum over $k$ in the second line as a limit,
\begin{multline*} 
\sum_{k=0}^{n-1}\frac{\po{-n-1}{k+2}-(k+2)\po{-n-2}{k+3}}{(k+1)\,(k+2)!}
\\
=\lim _{\ep\to0}\Bigg(\frac {n+1} {\ep}-\frac {n+1} {\ep}
{}_2F_1\!\left[\matrix -n,\ep\\2\endmatrix; 1\right]
+\frac {(n+1)(n+2)} {\ep}-\frac {(n+1)(n+2)} {\ep}
{}_2F_1\!\left[\matrix -n,\ep\\1\endmatrix; 1\right]
\Bigg).
\end{multline*}
Again, Chu--Vandermonde summation \eqref{eq:Chu-Van} can be applied
to evaluate the two $_2F_1$-series. If the result is substituted in
\eqref{eq:ex3}, we get
\begin{multline*}
\frac{1-(-1)^{n}(n+2)}{\po{n+1}{3}}\left(
	(-1)^{n+1}
        - {(n+1)\,(n+2)^2} 
+{(n+1)_3\sum_{k=1}^{n+2}\frac{1}{k+1}}
	\right)\\
+\frac{1-(-1)^{n}(n+2)}{n+3}
	\lim_{\varepsilon\rightarrow0}\left(
	\frac{(n+1)!-\po{2-\varepsilon}{n}}{(n+1)!\,\varepsilon}+
	(n+2)\frac{n!-\po{1-\varepsilon}{n}}{n!\,\varepsilon}
\right).
\end{multline*}
Using de l'Hospital's rule once more, we can compute the limit and obtain after
some simplification the expression
\begin{equation}
\label{eq:simplesumresult}
(1-(-1)^{n}(n+2))
\Biggl(\sum_{k=0}^{n}\frac{2}{k+1}
-	\frac{(-1)^{n}+2n^3+11n^2+19n+11
	}{(n+1)(n+2)(n+3)
	}
\Biggr).
\end{equation}

Now we turn our attention to the double sum in the second line of
\eqref{eq:weird}. Analogously to before,
we replace the term $1/(j+1)$ by $\int_0^1x^j\dx$. This enables
us to evaluate the inner harmonic sum $\sum _{j=0} ^{k}1/(j+1)$ to $\int
_{0} ^{1}(1-x^{k+1})/(1-x)\,dx$. We substitute this in the double sum
in the second line of \eqref{eq:weird}. Using hypergeometric
notation, the result is
\begin{equation*}
(-1)^n(n+2)\int_0^1\left(\frac{1}{1-x}\left(
	\HypF{2}{1}{n+3,-n-1}{2}{1}-
	\frac{1}{x}\HypF{2}{1}{n+3,-n-1}{2}{x}
	\right)+\frac{1}{x}\right){\, dx}.
\end{equation*}
The first hypergeometric series can simply be computed by 
Chu--Vandermonde summation \eqref{eq:Chu-Van}. To the second
hypergeometric series we apply 
the transformation formula (see \cite[(1.8.10),
terminating form]{SlatAC}) 
$$
{} _{2} F _{1} \!\left [ \matrix { a, -N}\\ { c}\endmatrix ; {\displaystyle
   z}\right ]  =
{\frac {    ({ \textstyle c-a}) _{N} } {({ \textstyle c}) _{N} }}
{} _{2} F _{1} \!\left [ \matrix { a,-N}\\ { 1 + a - c -
       N}\endmatrix ; {\displaystyle 1 - z}\right ]  ,
$$
where $N$ is a nonnegative integer.
These operations yield
\begin{equation} \label{eq:ex4}
(-1)^n\int_0^1\left(\frac{(-1)^{n+1}}{(1-x)}
-	\frac{(-1)^{n+1}}{x(1-x)}\HypF{2}{1}{n+3,-n-1}{1}{1-x}
	+\frac{n+2}{x}\right){\, dx}.
\end{equation}
Now we would like to apply the Euler beta integral formula
\eqref{eq:beta-integral} once more. However, this is not possible
just straightforwardly, because the beta integral on the left-hand
side of \eqref{eq:beta-integral} is not defined for $\beta=-1$.
In order to overcome this problem, we first rewrite the term $1/x(1-x)$
(which appears in the second term of the integrand in \eqref{eq:ex4})
as $1/x+1/(1-x)$, and then replace all occurrences of $1/x$ by
$\lim_{\ep\to0^+}x^{1-\ep}$, so that \eqref{eq:ex4} becomes
\begin{multline*}
\int_0^1\lim_{\varepsilon\rightarrow0^+}\Biggl(
\sum_{k=0}^{n+1}\frac{\po{n+3}{k}\,\po{-n-1}{k}}{k!^2}
x^{\varepsilon-1}\,(1-x)^k\\
+\sum_{k=1}^{n+1}\frac{\po{n+3}{k}\,\po{-n-1}{k}}{k!^2}
x^{\varepsilon}\,(1-x)^{k-1}
+(-1)^n(n+2)x^{\varepsilon-1}
\Biggr){\, dx}.
\end{multline*}
Next we interchange limit and integration, and apply the Euler beta
integral formula
\eqref{eq:beta-integral} wherever possible. In the result, the first
sum can easily be evaluated by Chu--Vandermonde summation
\eqref{eq:Chu-Van}. Subsequently, we compute the limit 
by using de l'Ho\-spi\-tal's rule again. This yields the expression
\begin{equation*}
(-1)^n(n+2)\left(\sum_{k=0}^n\frac{1}{k+1}+
	\sum_{k=1}^{n+1}\frac{1}{k+1}\right)+
\sum_{k=1}^{n+1}\frac{\po{n+3}{k}\,\po{-n-1}{k}}{k!^2\, k}.
\end{equation*}
Replacing once more the term ${1}/{k}$ in the
rightmost sum by $\int_0^1x^{k-1}\dx$, we
obtain the expression
\begin{equation*}
(-1)^n(n+2)\left(\sum_{k=0}^n\frac{1}{k+1}+
	\sum_{k=1}^{n+1}\frac{1}{k+1}\right)+
\int_0^1\frac {1} {x}\left(-1+{\HypF{2}{1}{n+3,-n-1}{1}{x}}\right){\, dx}.
\end{equation*}
To the $_2F_1$-series we apply
the transformation formula \eqref{eq:T2106}. In the result, we
replace any occurrence of $x^m$ by $\lim _{\ep\to0^+}x^{m+\ep}$, so that
we arrive at the expression
\begin{multline*}
(-1)^n(n+2)\left(\sum_{k=0}^n\frac{1}{k+1}+
	\sum_{k=1}^{n+1}\frac{1}{k+1}\right)\\
+\int_0^1\lim_{\varepsilon\rightarrow0^+}
	\left(-x^{\varepsilon-1}+(-1)^{n+1}(n+2)
	\sum_{k=0}^{n+1}\frac{\po{-n-1}{k}^2}{k!\,(k+1)!}(x-1)^kx^{n-k+\varepsilon}
\right){\, dx}.
\end{multline*}
Again, we interchange limit and integration and apply the Euler beta
integral formula
\eqref{eq:beta-integral} once more. Writing the result in
hypergeometric notation, we obtain
\begin{multline*}
(-1)^n(n+2)\left(\sum_{k=0}^n\frac{1}{k+1}+
	\sum_{k=1}^{n+1}\frac{1}{k+1}\right)\\
+\lim_{\varepsilon\rightarrow0^+}
\left(
	-\frac{1}{\varepsilon}
        +(-1)^{n+1}\frac{1}{n+2}
	\left(
		1-
		\HypF{2}{1}{-n-2,-n-2}{-n-1-\varepsilon}{1}
	\right)
\right).
\end{multline*}
The hypergeometric series can be computed by 
Chu--Vandermonde summation \eqref{eq:Chu-Van}. Once more we make use
of de l'Hospital's rule for the limit and obtain after some simplifications
the expression
\begin{equation}
\label{eq:doublesumresult}
-\frac{1+(-1)^n\left(n^2+3n+3\right)}{n+2}-
(1-(-1)^n(n+2))\sum_{k=0}^{n}\frac{2}{k+1}.
\end{equation}

When adding together \eqref{eq:simplesumresult}, the result of our
computation for the sum in the first line of \eqref{eq:weird}, and
\eqref{eq:doublesumresult}, the result of our computation for the
double sum in the second line of \eqref{eq:weird},
the harmonic sums cancel, and it is easy to verify that, magically, 
the remaining
terms simplify to the right-hand side of \eqref{eq:weird}.
\end{proof}

\begin{rem}
Peter Paule demonstrated to us, that the identity \eqref{eq:weird}
can also be proved algorithmically. Clearly, the 
Gosper--Zeilberger algorithm \cite{PeWZAA,ZeilAM,ZeilAV} finds a
recurrence for the sum in the first line of \eqref{eq:weird}. Carsten
Schneider's extension of Karr's algorithm \cite{KarrAA}, implemented
by Schneider, finds a recurrence for the double sum in the second
line of \eqref{eq:weird}. Finally, Mallinger's {\sl Mathematica}
package {\tt GeneratingFunctions} \cite{MallAA} 
or Salvy and Zimmermann's {\sl
Maple} package {\tt gfun} \cite{SaZiAA} can be used to combine these two
recurrences into one, a recurrence of order 10. It is then routine to
check (preferably on the computer) that the right-hand side of
\eqref{eq:weird} satisfies
this same recurrence. However, in the present implementation, these algorithms
are not able to find the explicit evaluations, in terms of harmonic
numbers, of the sums in the first and second line of
\eqref{eq:weird}, given in \eqref{eq:simplesumresult} and
\eqref{eq:doublesumresult}, respectively.
\end{rem}

%
%


\section{Evaluations of Hankel determinants featuring Bernoulli
polynomials}
\label{sec:further}

There are several theorems hidden in the body of this paper. 
Among these are evaluations of Hankel determinants of
{\em Bernoulli polynomials} evaluated at special values. 
Recall that the $l$-th Bernoulli polynomial is defined by
\begin{equation*} 
B_l(x):=\sum_{k=0}^l\binom{l}{k}B_{l-k}x^k.
\end{equation*}
The Hankel determinants of Bernoulli polynomials of which we are talking
are special cases of the determinant $B(N;x)$ given by
\begin{equation} \label{eq:BPol}
B(N;x):=\det_{1\leq i,j\leq N}\left(B_{i+j}(x)\right).
\end{equation}
This is in a fundamental way different from the Hankel determinant
$$\det_{0\le i,j\le N}(B_{i+j}(x)),$$
which has been considered earlier (see \cite[Sec.~5]{AlCaAA}). 
(Note that the
difference is that, in the latter determinant, indices start
already with 0.) As is not difficult to see (cf\@.
\cite[p.~419]{MuirAD} or \cite[Lemma~15]{KratBN}), the latter 
determinant does in fact not depend on $x$
(i.e., the powers of $x$ cancel in the expansion of the
determinant), so that its value is equal to its value at $x=0$,
which, in turn,
is given by Theorem~\ref{thm:HankelBernoulli} with $a=b=1$, $c=d=0$.
This is in sharp contrast to the Hankel determinant \eqref{eq:BPol},
where the powers of $x$
do not cancel, so that \eqref{eq:BPol} is a nontrivial
polynomial in $x$. As such, the evaluation of the determinant
\eqref{eq:BPol} is much more difficult. Below, we provide evaluations
of \eqref{eq:BPol} for $x=-1$, $x=-1/2$ and $x=1/2$. Needless to say that 
the evaluation in the special case $x=0$ (and as well in the special
case $x=1$) is given by 
Theorem~\ref{thm:HankelBernoulli} with $a=b=c=d=1$.

First of all, in the proof of Lemma~\ref{lem:Det3}, we observed that
(in symbolic notation $B^k\equiv B_k$)
$$\det_{1\le i,j\le N}\big((B-1)^{i+j}\big)=
\det_{1\le i,j\le N}\big(B^{i+j-2}(B-1)^2\big).
$$
The determinant on the right-hand side was then evaluated in
Lemma~\ref{lem:extra-Bernoulli}. The linear combination $(B-1)^l$
of Bernoulli numbers is nothing else but $B_l(-1)$, the $l$-th
Bernoulli polynomial evaluated at $-1$. Thus, we obtain the following
corollary.

\begin{cor}
Let $N$ be a positive integer. Then there holds
\begin{multline}
B(N;-1)=
\det_{1\le i,j\le N}\left(B_{i+j}+(-1)^{i+j}(i+j)\right)
\\=
{(-1)^{{N(N-1)}/{2}}} {6^{-N}}(1 + 2N + 5N^2 + 4N^3 + N^4)
\prod_{i=1}^N\left(\frac{i(i+1)^4(i+2)}{(2i+1)(2i+2)^2(2i+3)}\right)^{N-i}.
\end{multline}
\end{cor}

Our next evaluation results from the determinant evaluation in
Lemma~\ref{lem:Det2}.

\begin{thm}
\label{prop:extrad-det}
Let $m$ be a positive integer. Then for
\begin{equation*}
B(N;-1/2)=
\det_{1\le i,j\le N}\left((2^{1-i-j}-1)B_{i+j}-(-1/2)^{i+j-1}(i+j)\right)
\end{equation*}
there hold
\begin{multline} \label{eq:B1}
B(2m-1;-1/2)=\frac{
	\left( -1 \right)^{m - 1}\,\left(2 m - 1\right)!^2
	\left( \prod_{i = 1}^{2 m - 1}i! \right)^4
	\left( \prod_{i = 1}^{m}\left(2 i - 1\right) ! \right)^4
}{
	2^{6 \left(m-1\right)}\,\left(m - 1\right)!^6
	\left( \prod_{i = 1}^{4 m - 1}i! \right)
}\\\times
\left( 3 + 8 m + 
	\frac{
		8 \left(m - 1\right)  \left(2 m + 1\right)
	}{
		3 m \left(m + 1\right)  \left(2 m - 3\right)  
      \left(2 m - 1\right)
	}
   \HypF{4}{3}{
		3,\frac{3}{2},2-m,\frac{3}{2}+m
	}{
		\frac{5}{2},2+m,\frac{5}{2}-m
	}{
		1
	}
\right)
\end{multline}
and
\begin{multline} \label{eq:B2}
B(2m;-1/2)=\frac{
	\left( -1 \right)^{m}
	\left( \prod_{i = 1}^{2 m}i! \right)^4
	\left( \prod_{i = 1}^{m + 1}\left(2 i - 1\right) ! \right)^4
}{
	2^{6 m}\,m!^6\,\left( \prod_{i = 1}^{4 m + 1}i! \right)
}\\\times
\left(
	1 + 8 m - 
	\frac{
		8 m \left(2 m + 3\right)
	}{
		3 \left(m + 1\right)  \left(m + 2\right)  
      \left(2 m - 1\right)  \left(2 m + 1\right)
	} 
	\HypF{4}{3}{
		{3,\frac{3}{2}},1- m,{\frac{5}{2}+m}
	}{
		{\frac{5}{2}},3+m,{\frac{3}{2}-m} 
	}{
		1
	}
\right).
\end{multline}
\end{thm}
\begin{proof} 
Consider the determinant $D(n,n-1,N)$ (see
\eqref{eq:gen-det} for definition). Factor
$2(n-1/2)$ out of each column of $D(n,n-1,N)$, and then
set $n=1/2$. By the appropriate variant of \eqref{eq:BernoulliP} and
de l'Hospital's rule,
this yields the Hankel determinant $B(N;-1/2)$. 
On the other
hand, we have evaluated $D(n,n-1,N)$ in Lemma~\ref{lem:Det2}.
Then by dividing the results in \eqref{eq:Det2odd} and
\eqref{eq:Det2even} by $2^N(n-1/2)^N$ and then setting $n=1/2$,
we obtain the expressions on the right-hand sides of \eqref{eq:B1}
and \eqref{eq:B2}.
\end{proof}

When solving the enumeration of rhombus tilings of a hexagon which
contain the central rhombus, the determinant $D(n,n+1,N)$ 
was explicitly \cite{HeGeAA} (see the proof of Lemma~10 and
Proposition~14 in \cite{HeGeAA}) or
implicitly \cite{CiucKratAB,FuKrAC} evaluated (compare
Proposition~\ref{thm:HeGe-det}). If one adapts the
preceding proof to this case, one obtains the following result.

\begin{thm}
\label{prop:extrad-det2}
Let $m$ be a positive integer. Then for
\begin{equation*}
B(N;1/2)=
\det_{1\le i,j\le N}\left((2^{1-i-j}-1)B_{i+j}\right)
\end{equation*}
there hold
\begin{multline} \label{eq:B3}
B(2m-1;1/2)=
{\frac {{{\left( 2 m \right) !}^2} 
   }
   {{2^{4 m-2}}\,{{\left( m-1 \right) !}^2}\,{{m!}^2}}}
     \bigg( \prod_{i = 1}^{2 m}
        {\frac {{{\left( i-1 \right) !}^5}}
          {\left( 2m+i-1 \right) !}} \bigg)
 \sum_{i = 0}^{m-1}
     \frac {{{\left( -1 \right) }^{m-i}}} {\left( 2m- 2 i-1 \right)}
        {\frac {  ({ \textstyle {\frac 1 2} - i}) _{2 i} }
          {  {{i!}^2}}} 
\end{multline}
and
\begin{equation} \label{eq:B4}
B(2m;1/2)=
{\frac {{{\left( 2m + 2  \right) !}^2} 
   }
   {{2^{4 m+2}}\,{{m!}^2}\,{{\left( m+1 \right) !}^2}}}
     \bigg( \prod_{i = 1}^{2 m+1}
        {\frac {{{\left( i-1 \right) !}^5}} {\left( 2 m+i \right) !}}
        \bigg)
 \sum_{i = 0}^{m}
     \frac {{{\left( -1 \right) }^{m-i}}} {\left( 2 m-2i+1 \right)}
        {\frac {  ({ \textstyle {\frac 1 2} - i}) _{2 i} }
          {  {{i!}^2}}} .
\end{equation}
\end{thm}
\begin{proof} 
Consider the determinant $D(n,n+1,N)$. Factor
$2(n+1/2)$ out of each column of $D(n,n+1,N)$, and then
set $n=-1/2$. By the appropriate variant of \eqref{eq:BernoulliP} and
de l'Hospital's rule,
this yields the Hankel determinant $B(N;1/2)$. 
On the other
hand, by Proposition~\ref{thm:HeGe-det} (with $N$ replaced by $N+1$,
$l=n$, and $M=2n-N$) we know that $D(n,n+1,N)$ counts, up to a
multiplicative constant, the number of rhombus tilings of a hexagon
with sides $N+1,2n-N,N+1,N+1,2n-N,N+1$, which contain the central
rhombus. This enumeration problem was solved in 
\cite{CiucKratAB,FuKrAC,HeGeAA}. If in the result we perform
the according manipulations and then set $n=-1/2$,
we obtain the expressions on the right-hand sides of \eqref{eq:B3}
and \eqref{eq:B4}.
\end{proof} 

\section{Concluding comments and open problems}
\label{sec:open}

We conclude this article by pointing to open questions
which are raised by this work.

\smallskip
(1)
In Corollary~\ref{thm:Propp} it was demonstrated that, for $M$ close
to $N$, the number of 
rhombus tilings of a hexagon with side lengths $N,M,N,N,M,N$,
which contain the rhombus above and next to the center of the hexagon, equals
$(\frac {1} {3}+r(N,M))\,T(N,M)$, where $T(N,M)$ is the total number of
rhombus tilings of the hexagon, and where $r(N,M)$ is a ``closed form"
expression. (Magically,
the value of $1/3$ which appears here is, 
according to \cite[Corollary~3]{CiucKratAB}, the {\em exact\/} 
proportion of the rhombus tilings that contain the
{\em central\/} rhombus in the total number of rhombus tilings of a
hexagon with side lengths $2n-1,2n,2n-1,2n-1,2n,2n-1$ or with side
lengths $2n,2n-1,2n,2n,2n-1,2n$.) 
As we mentioned in the Introduction, it is easy to derive
many more such results, also for the central rhombus and the other
two cases that were considered in
Theorems~\ref{thm:NEven}--\ref{thm:NOdd2}.
Our proof, given in Section~2, consisted of starting with
the expressions \eqref{eq:MEven} and \eqref{eq:MOdd} and applying
Zeilberger's algorithm to establish the simplification of 
the sum in these expressions when
$m$ and $n$ are close. This is, unfortunately, not conceptual (as it
just {\em verifies}, but does not {\em derive} the result), and
therefore does not explain why these simplifications take place.
The fact that apparently many more such results exist indicates
that there must be a hypergeometric transformation formula lurking in
the background, which we were, however, unable to discover. (It is obvious
that the sums in \eqref{eq:MEven}--\eqref{eq:NOdd2} can be written as
very-well-poised $_7F_6$-series --- see e.g\@.
\eqref{eq:very-well-poised} --- and, by means of Whipple's
transformation formula \eqref{eq:Whipple}, can therefore be transformed
into balanced $_4F_3$-series, to which, in turn, we could apply
Sears' $_4F_3$ transformation formulas. However, it seems that this
does not suffice to find the desired identity which would ``explain"
Corollary~\ref{thm:Propp}.)

\smallskip
(2)
Is it possible to find a uniform formula for the number of
rhombus tilings of a hexagon with side lengths $N,M,N,N,M,N$, which
contain an {\em arbitrary} (but fixed) rhombus on the ``vertical"
symmetry axis (i.e., the symmetry axis
which runs in parallel to the sides of length $M$)?
Recall that (as we mentioned already in the Introduction) 
in \cite{FuKrAC} such a formula was found for the ``horizontal" symmetry
axis (i.e., the symmetry axis which cuts through the sides of length
$M$). In contrast, here we encountered increasing difficulties in the
proofs of our enumerations the farther we moved the rhombus which is
contained in every tiling from the center. Recall that for solving
our enumeration problems we needed to compute the determinants
$D(n,n-t,N)$ (see \eqref{eq:gen-det} for definition) for $t=0,1,2$.
For the case of a rhombus which is even farther away from the center, we
would have to evaluate this determinant for even larger values of $t$. The
increasing difficulties in doing this arise in Step~5 (compare the
proof of Lemma~\ref{lem:Det1}) of the computation. The previous
steps, Steps~1--4, would prove that
$$D(n,n-t,N)=\big(\text {product of linear factors in $n$}\big)\cdot
S(n,N,t),$$
where $S(n,N,t)$ is a polynomial of degree $2\ceil{(N+t)/2}$ (compare
\eqref{eq:polydef}, \eqref{eq:polydef2}, and \eqref{eq:polydef3}).
Thus, in order to determine $S(n,N,t)$, the larger $t$ becomes, the
more evaluations of $S(n,N,t)$ at special values of $n$ 
(or other informations about $S(n,N,t)$) we need. (The computations in
\cite{FuKrAC} have exactly the opposite behaviour: The farther the
rhombus which is contained in every tiling is moved away from the
center, the {\em smaller} in degree becomes the irreducible
polynomial in the result.) Even worse, the larger $t$ becomes,
the more difficult it becomes to obtain these special values.
(Remember, for example, the difficulty of evaluation of $R(n,N)$ at
$n=1$ via Lemmas~\ref{lem:extra-Bernoulli} and \ref{lem:weird}.)

That the problem that we considered here is at a different level of
complexity than the problem in \cite{FuKrAC} is also indicated by the
(partially conjectural) form of the asymptotic behaviour of the
proportion of the rhombus tilings that contain this particular
rhombus in the total number of rhombus tilings. While the asymptotic
behaviour is totally smooth when the rhombus which is contained in
every tiling is moved along the ``horizontal" symmetry axis (see
\cite[Theorem~1.3]{FuKrAC}), the conjectured form
\cite[Conjecture~1]{CoLPAA} of the asymptotics when
the rhombus which is contained in
every tiling is moved along the ``vertical" symmetry axis behaves
nonsmoothly. It is increasing for some time when the rhombus is moved
away from the center, but at some point, when the rhombus enters the ``arctic
region" near the (top or bottom) corner, it becomes 1 and stays
1 from thereon. Thus, a formula for exact enumeration must,
somehow, reflect this nonsmooth asymptotic behaviour.

Is there a way to overcome these difficulties?

\smallskip
(3)
In Theorem~\ref{thm:HankelBernoulli} only $c$ or $d$ may be 0, but
not $a$ or $b$. In fact, Theorem~\ref{thm:HankelBernoulli} is wrong
if $a=0$ or $b=0$. But, apparently, not terribly wrong. 
Lemma~\ref{lem:extra-Bernoulli}
shows the evaluation of the determinant in \eqref{eq:HankelBernoulli}
with $a=b=0$, $c=d=2$. Remarkably, the result is almost identical
with the right-hand side in \eqref{eq:HankelBernoulli}, the only
difference being the polynomial in $n$ of fourth degree in
\eqref{eq:extra-Bernoulli}. In fact, computer experiments suggest
that a much more general result holds.

\begin{conj}
For positive integers integers $c,d$ there holds
\begin{multline} \label{eq:HankelBernoulli-extra1}
\det_{1\le i,j\le n}
\left(
	B^{i+j-2}\,\po{-B+1}{c-1}\,\po{-B+1}{d-1}
\right)\\
=(-1)^{n(n-1)/2}\left({\frac{\left( c -1\right) !\,\left( c - 1\right) !\,
     \left( d - 1\right) !\,\left( d - 1\right) !}{\left(
        c + d - 1\right) !}}\right)^n\kern6cm\\
\times\prod_{i=1}^{n-1}\Biggl({\frac{i\,\left( c + i - 1\right) \,
     \left( c + i - 1\right) \,\left( d + i - 1\right)\,
\left( d + i - 1\right) \,
     \left( c + d + i - 2\right) }
	  {\left( c + d + 
       2 i - 3\right) \,{{\left( c + d + 2 i - 2\right) }^2\,
     \left( c + d + 2 i - 1\right) }}}
\Biggr)^{n-i} P(n;c,d),
\end{multline}
where $P(n;c,d)$ is a certain polynomial in $n$ of degree
$2(c+d-2)$.

Furthermore, for positive integers integers $b,c,d$ there holds
\begin{multline} \label{eq:HankelBernoulli-extra2}
\det_{1\le i,j\le n}
\left(
	B^{i+j-1}\,\po{B+1}{b-1}\,\po{-B+1}{c-1}\,\po{-B+1}{d-1}
\right)\\
=(-1)^{n(n-1)/2}\left({\frac{\left( c -1\right) !\,\left(b + c - 1\right) !\,
     \left(d - 1\right) !\,\left(b + d - 1\right) !}{\left(
        b + c + d - 1\right) !}}\right)^n\\
\times\prod_{i=1}^{n-1}\Biggl({\frac{i\,\left( c + i - 1\right) \,
     \left(b + c + i - 1\right) \,\left( d + i - 1\right)}
	  {\left( b + c + d + 
       2 i - 3\right) \,{{\left( b + c + d + 2 i - 2\right) }^2}}}\\
\times{\frac{\left(b + d + i - 1\right) \,
     \left( b + c + d + i - 2\right) }{
     \left( b + c + d + 2 i - 1\right) }}
\Biggr)^{n-i} R(n;b,c,d),
\end{multline}
where $R(n;b,c,d)$ is a certain rational function in $n$, which can be
written with a
numerator of degree $c+d-2$ and a denominator of degree $b-1$.
\end{conj}

In principle, our approach of proving Lemma~\ref{lem:extra-Bernoulli}
(the special case $c=d=2$ of \eqref{eq:HankelBernoulli-extra1}),
which consisted of using linearity of the determinant in order to
break it into several pieces, to each of which we could either 
apply Theorem~\ref{thm:HankelBernoulli} or Theorem~\ref{thm:Leclerc},
should make a proof of the above conjecture possible. However,
serious difficulties have to be expected in actually doing the
calculations, in particular, when working through a generalized form
of Lemma~\ref{lem:weird}. We believe that, in view of the simplicity
of the result \eqref{eq:extra-Bernoulli} and of the conjectured
results \eqref{eq:HankelBernoulli-extra1} and
\eqref{eq:HankelBernoulli-extra2}, there must be a more elegant way
to attack these Hankel determinant evaluations, in particular, if one
also desires to obtain explicit forms for the polynomial $P(n;c,d)$
and the rational function $R(n;b,c,d)$.


\bigskip

{\small
\noindent
{\sc Note.} \scriptsize
Since first versions of this article were distributed,
Ilse Fischer (``Enumeration of rhombus tilings of a hexagon which
contain a fixed rhombus in the centre", preprint,
{\tt math/9906102}) generalized
Theorems~\ref{thm:MEven} and \ref{thm:MOdd} to arbitrary semiregular
hexagons.
}

\ifx\undefined\bysame
\newcommand{\bysame}{\leavevmode\hbox to3em{\hrulefill}\,}
\fi


\end{document}